\documentclass[titlepage,11pt]{article}
\usepackage{latexsym}
\usepackage{amsfonts}
\usepackage{amsmath}
\newcommand\blackslug{\hbox{\hskip 1pt \vrule width 4pt height 8pt depth 1.5pt
        \hskip 1pt}}
\newcommand\bbox{\hfill \quad \blackslug \medbreak}
\def\d{\hbox{-}}
\def\c{\hbox{-}\cdots\hbox{-}}
\def\l{,\ldots,}

\usepackage[dvips]{graphicx}
\usepackage{color}
\usepackage{epsfig}

\title{Excluded minors in cubic graphs}
\author{Neil Robertson\thanks{Research partially supported by DIMACS,
and by ONR grant N00014-92-J-1965, and by NSF grant DMS-8903132,
and partially performed under a consulting agreement with Bellcore.} \\
Ohio State University, Columbus, Ohio 43210 
\and
Paul Seymour\thanks{This research was partially performed while Seymour was 
employed at Bellcore in Morristown, New Jersey, and partially supported by
ONR grant N00014-10-1-0680 and
NSF grant DMS-1265563.}\\
Princeton University, Princeton, New Jersey 08544
\and
Robin Thomas\thanks{Research partially supported by DIMACS,
by ONR grant N00014-93-1-0325, and by NSF grants DMS-9303761 and DMS-1202640, 
and partially performed under a consulting agreement with Bellcore.} \\
Georgia Institute of Technology, Atlanta, Georgia 30332}
\date{February 1, 1995; revised \today}

\newtheorem{thm}{}[section]

\newcommand{\Proof}{\noindent{\bf Proof.}\ \ }
\newcommand{\Subproof}{\noindent{\em Subproof.}\ \ }
\newcommand{\he}{homeomorphic embedding }

\begin{document}
\maketitle
\begin{abstract}
Let $G$ be a cubic graph, with girth at least five, such that for 
every partition
$X,Y$ of its vertex set with $|X|, |Y| \geq 7$ there are at least six
edges between $X$ and $Y$.
We prove that if there is no \he of the Petersen graph in $G$, and $G$ is
not one particular 20-vertex graph, then either 
\begin{itemize}
\item $G \setminus v$
is planar for some vertex $v$, or 
\item $G$
can be drawn with crossings in the plane, but with only two crossings,
both on the infinite region.
\end{itemize}
We also prove several other theorems of the same kind.
\end{abstract}

\section{Introduction}


All graphs in this paper are simple and finite.
Circuits have no repeated vertices or edges; the {\em girth} of a graph is the length of the shortest circuit.
If $G$ is a graph and $X \subseteq V(G), \delta_{G} (X)$ or
$\delta (X)$ denotes the set of edges with one end in $X$ and the other in
$V(G)  \setminus X$.
We say a cubic graph $G$ is {\em cyclically} $k$-{\em connected}, for
$k \geq 1$ an integer, if $G$ has girth $\geq k$, and
$| \delta_{G} (X) | \geq k$ for every $X \subseteq V(G)$
such that both $X$ and $V(G) \setminus X$ include the vertex set of a circuit of $G$.

A {\em homeomorphic embedding} of a graph $G$ in a graph $H$ is a function $\eta$
such that
\begin{itemize}
\item
for each $v \in V(G)$, $\eta (v)$ is a vertex of $H$, and
$\eta (v_{1} ) \neq \eta (v_{2} )$ for all distinct $v_{1},v_{2} \in V(G)$
\item
for each $e \in E(G)$, $\eta (e)$ is a path of $H$ with ends
$\eta (v_{1})$ and $\eta (v_{2})$, where $e$ has ends
$v_{1},v_{2}$ in $G$; and no edge or internal vertex of
$\eta (e_{1})$ belongs to $\eta (e_{2})$, for all distinct
$e_{1},e_{2} \in E(G)$
\item
for all $v \in V(G)$ and $e \in E(G)$, $\eta (v)$ belongs to
$\eta (e)$ if and only if $v$ is an end of $e$ in $G$.
\end{itemize}
We denote by $\eta(G)$ the subgraph of $H$ consisting of all the vertices
$\eta(v)\;(v\in V(G))$ and all the paths $\eta(e)\;(e\in E(G))$.
We say that $H$ {\em contains} $G$ if there is a
\he of $G$ in $H$.

\begin{figure} [h!]
\centering
\includegraphics{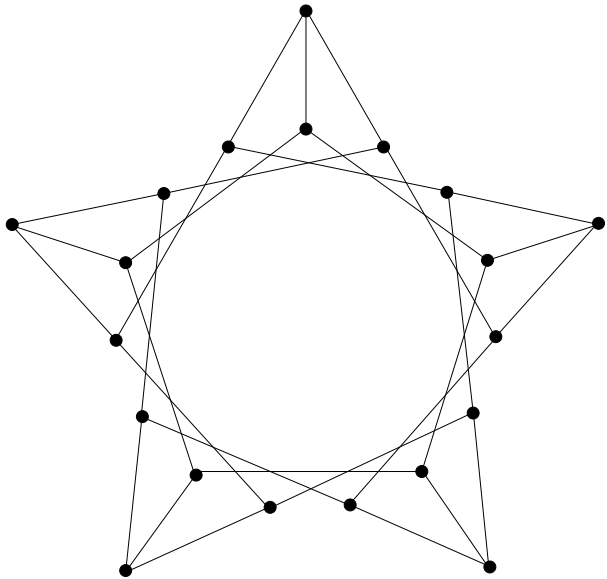}
\caption{Starfish}
\label{starfishfig}
\end{figure}

Let us say that $G$ is {\em theta-connected} if $G$
is cubic and cyclically five-connected, and
$| \delta_{G} (X) | \geq 6$ for all
$X \subseteq V(G)$ with $|X| ,|V(G)  \setminus X| \geq 7$.
We say $G$ is {\em apex} if $G \setminus v$ is planar for some vertex
$v$ (we use $\setminus$ to denote deletion); and $G$ is
{\em doublecross} if it can be drawn in the plane with only two
crossings, both on the infinite region.
Our goal in this paper is to give a construction for all theta-connected graphs 
not containing Petersen (we define {\em Petersen} to be the Petersen graph.)
This is motivated by a result of a previous paper~\cite{RST2}, where we showed that to 
prove Tutte's conjecture~\cite{Tutte} that every two-edge-connected cubic graph not 
containing Petersen is three-edge-colourable, it is enough to prove the same
for theta-connected graphs not containing Petersen, and for apex graphs.

The graph {\em Starfish} is shown in Figure 1.
Our main result is the following.

\begin{thm}\label{main}
Let $G$ be theta-connected.
Then $G$ does not contain Petersen if and only if either $G$ is apex, or
$G$ is doublecross, or $G$ is isomorphic to Starfish.
\end{thm}

The ``if'' part of \ref{main} is easy and we omit it. (It is enough to check 
that Petersen itself is not apex or doublecross, and is not contained 
in Starfish.)
The ``only if'' part is an immediate consequence of the following
three theorems. The graph {\em Jaws} is defined in Figure 2.

\begin{figure} [h!]
\centering
\includegraphics{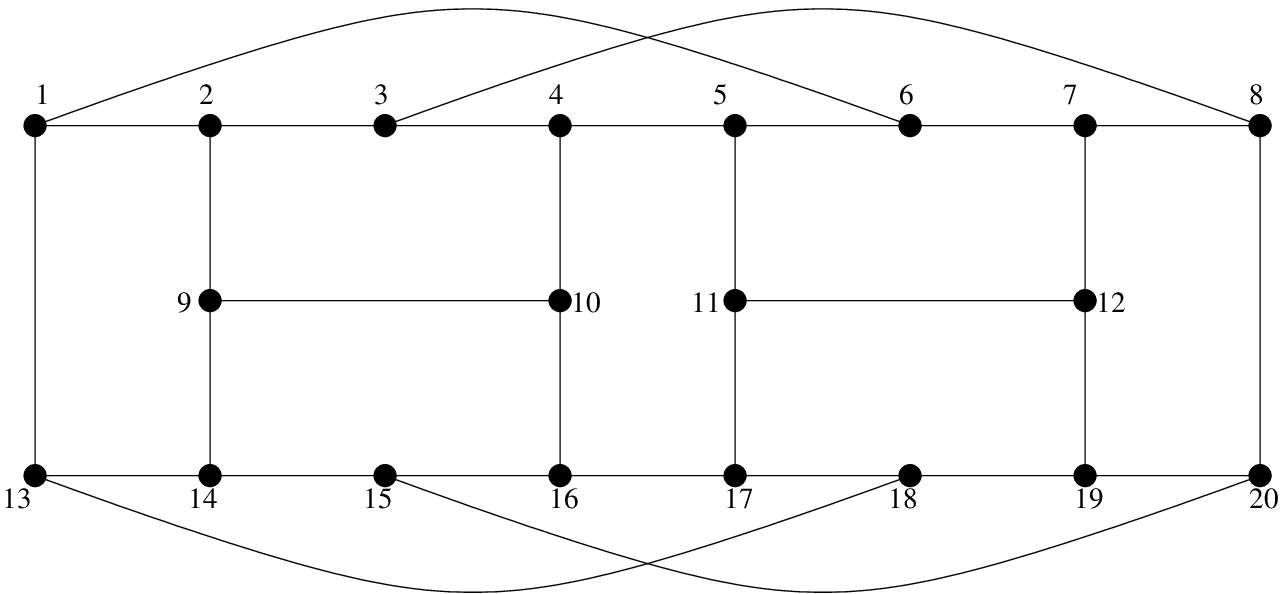}
\caption{Jaws}
\label{jawsfig}
\end{figure}

\begin{thm}\label{starfish0}
Let $G$ be theta-connected, and not contain Petersen.
If $G$ contains Starfish then $G$ is isomorphic to Starfish.
\end{thm}

\begin{thm}\label{main0}
Let $G$ be theta-connected, and not contain Petersen.
If $G$ contains Jaws then $G$ is doublecross.
\end{thm}

\begin{thm}\label{none0}
Let $G$ be theta-connected, and not contain Petersen.
If $G$ contains neither Jaws nor Starfish, then $G$ is apex.
\end{thm}
\begin{figure} [h!]
\centering
\includegraphics{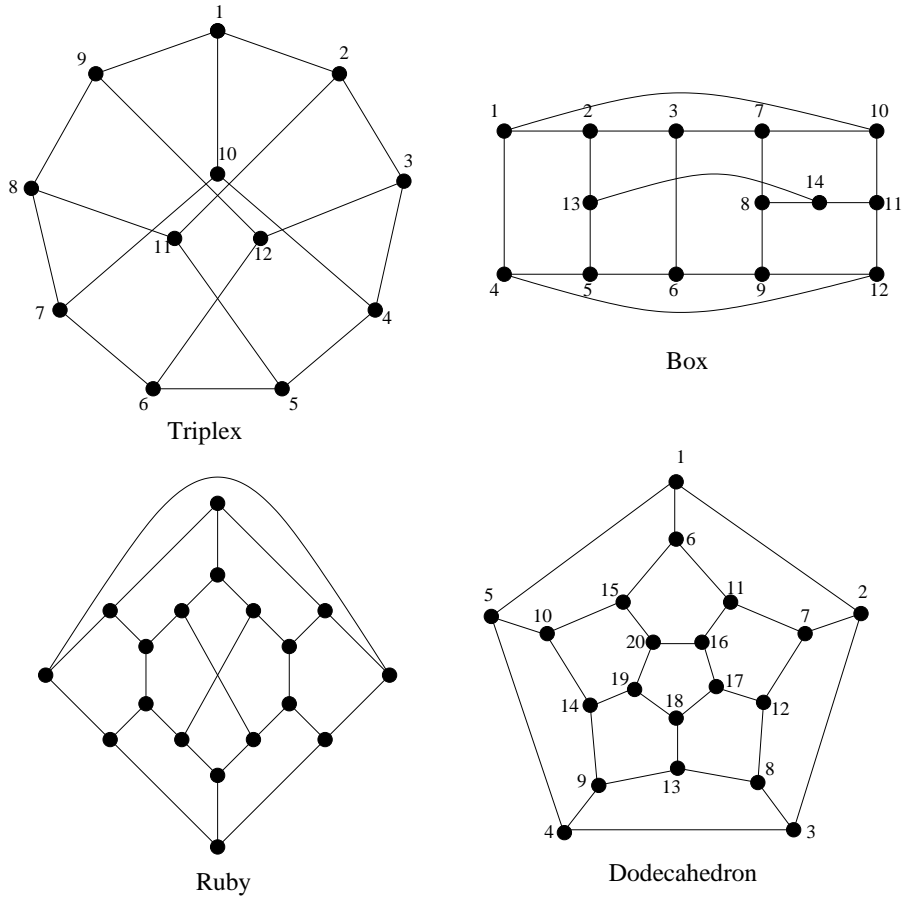}
\caption{Triplex, Box, Ruby and Dodecahedron}
\label{boxfig}
\end{figure}

\ref{starfish0}, proved in section 17, is an easy consequence of a theorem
of a previous paper \cite{RST}, and \ref{main0} is proved in section 18.
The main part of the paper is devoted to proving \ref{none0}.
Our approach is as follows.

A graph $H$ is {\em minimal} with property $P$ if there is no graph
$G$ with property $P$ such that $H$ contains $G$ and $H$
is not isomorphic to $G$.
In Figure 3 we define four more graphs, namely {\em Triplex, Box, Ruby}
and {\em Dodecahedron}.
A theorem of McCuaig~\cite{McCuaig} asserts

\begin{thm}\label{mccuaig}
Petersen, Triplex, Box, Ruby and Dodecahedron are the only graphs
minimal with the property of being cubic and cyclically five-connected.
\end{thm}

We shall prove the following three theorems.

\begin{thm}\label{nonplanar0}
Petersen, Triplex, Box and Ruby are the only graphs minimal
with the property of being cyclically five-connected and non-planar.
\end{thm}

A graph $G$ is {\em dodecahedrally-connected} if it is cubic
and cyclically five-connected, and for every
$X \subseteq V(G)$ with
$|X|, |V(G)  \setminus X| \geq 7$ and $| \delta_{G} (X) | =5$, $G|X$ cannot
be drawn in a disc $\Delta$ such that the five vertices in $X$ with
neighbours in $V(G) \setminus X$ are drawn in $bd( \Delta )$.

\begin{thm}\label{crossingno0}
Petersen, Triplex and Box are the only graphs minimal with the property
of being dodecahedrally-connected and having crossing number at least two.
\end{thm}

We say $G$ is {\em arched} if $G \setminus e$ is planar for some edge $e$.

\begin{thm}\label{arched0}
Petersen and Triplex are the only graphs minimal with the property
of being dodecahedrally-connected and not arched.
\end{thm}

Then we use \ref{arched0} to find all the graphs minimal with the property of
being dodecahedrally-connected and non-apex (there are six).
Let us say $G$ is {\em die-connected} if it is dodecahedrally-connected and
$|\delta_{G} (X) | \geq 6$ for every
$X \subseteq V(G)$ with $|X|, |V(G)  \setminus X| \geq 9$.
We use the last result to find all graphs minimal with the property
of being die-connected and non-apex (there are nine); and
then use that to find the minimal graphs with the property of
being theta-connected and non-apex.
There are three, namely Petersen, Starfish, and Jaws, and from this \ref{main0} follows.

\section{Extensions}

It will be convenient to denote by $ab$ or $ba$ an edge with ends $a$
and $b$ (since we do not permit parallel edges, this is unambiguous).
Let $ab$ and $cd$ be distinct edges of a graph $G$.
They are {\em diverse} if $a,b , c,d$ are all distinct and $a,b$ are not
adjacent to $c$ or $d$.
We denote by $G + (ab, cd)$ the graph obtained from $G$ as follows:
delete $ab$ and $cd$, and add two new vertices $x$ and $y$ and five new
edges $xa, xb, yc, yd, xy$.
We call $x,y$ (in this order) the {\em new vertices} of
$G+ (ab, cd)$.
Multiple applications of this operation are denoted in the natural way;
for instance, if $e,f \in E(G)$ are distinct, and
$G' = G + (e,f)$, and $g,h \in E(G')$ are distinct, we write
$G+ (e,f) + (g,h)$ for $G' + (g,h)$.

Similarly, let $ab, cd, ef$ be distinct edges of $G$, where
$a,b,c,d,e,f$ are all distinct.
We denote by $G+ (ab, cd, ef)$ the graph obtained by deleting
$ab, cd$ and $ef$, and adding four new vertices $x,y,z,w$, and nine new
edges $xa, xb, yc, yd, ze, zf, wx, wy, wz$; and call
$x,y,z,w$ (in this order) the {\em new vertices} of
$G+ (ab, cd, ef)$.

A path has no ``repeated'' vertices or edges.
Its first and last vertices are its {\em ends}, and its first
and last edges are its {\em end-edges}. Its other vertices and edges are
called {\em internal} vertices and edges.
A path with ends $s$ and $t$ is called an $(s,t)$-{\em path}.
If $P$ is a path and $s,t \in V(P)$, the subpath of $P$ with ends $s$
and $t$ is denoted by $P[s,t]$.
Let $\eta$ be a \he of $G$ in $H$.
An $\eta$-{\em path} in $H$ is a path $P$ with distinct ends
both in $V(\eta (G))$, but with no other vertex or edge in $\eta (G)$.
Let $G,H$ both be cubic, and let $\eta$ and $P$ be as above, where
$P$ has ends $s$ and $t$, with $s \in V(\eta (e))$ and
$t \in V(\eta (f))$.
We can sometimes use $P$ to obtain a new \he $\eta'$ of $G$ in $H$, equal to $\eta$ except as follows:
\begin{itemize}
\item
If $e=f$, let $e=uv$, where $\eta (u),s,t,\eta (v)$ lie in
$\eta (e)$ in order.
Define 
$$\eta '(e)
= \eta (e) [\eta (u),s] \cup P \cup \eta (e)[t, \eta (v)].$$
\item
If $e \neq f$ but they have a common end, let $e=uv$ and $f=vw$ say, and let $g$ be the third edge of $G$ incident with $v$.
Define $\eta'$ by:
\begin{eqnarray*}
\eta'(v) &=& t,\\ 
\eta'(e) &=& \eta(e)[\eta(u),s] \cup P,\\ 
\eta'(f) &=& \eta(f)[t,\eta(w)],\\
\eta'(g) &=& \eta(g)\cup \eta(f)[\eta(v),t].
\end{eqnarray*}
\item
If $e,f$ have no common end, but one end of $e$ is adjacent to one end of
$f$, let $e=uv$, $f=wx$ and $g=vw$ say. Let $h,i$ be the third edges at $v,w$ respectively.
Define $\eta'$ by:
\begin{eqnarray*}
\eta'(v)&=&s,\\ 
\eta'(w) &=& t,\\ 
\eta'(e) &=& \eta(e)[\eta(u),s],\\
\eta'(f) &=& \eta(f)[t, \eta(x)],\\
\eta'(g) &=& P,\\
\eta'(h) &=& \eta(h)\cup \eta(e)[s,\eta(v)],\\
\eta'(i) &=& \eta(i)\cup \eta(f)[\eta(w),t].
\end{eqnarray*}
\end{itemize}
In the first two cases we say that $\eta'$ is obtained from $\eta$
by {\em rerouting e along P}, and in the third case by
{\em rerouting g along P}.
If $\eta$ is a \he of $G$ in $H$, an $\eta$-{\em bridge}
is a connected subgraph $B$ of $H$ with
$E(B \cap \eta(G)) = \emptyset$, such that either
\begin{itemize}
\item
$|E(B) | =1$, $E(B) = \{e\}$ say, and both ends of $e$ are in
$V(\eta (G))$, or
\item
for some component $C$ of $H\setminus V(\eta (G))$, $E(B)$
consists of all edges of $H$ with at least one end in $V(C)$.
\end{itemize}
It follows that every edge of $H$ not in $\eta (G)$ belongs to a
unique $\eta$-bridge.
We say that an edge $e$ of $G$ is an $\eta$-{\em attachment} of
an $\eta$-bridge $B$ if $\eta (e) \cap B$ is non-null.

\section{Frameworks}

We shall often have a cubic graph $G$, such that
$G$ (or sometimes, most of $G$)
is drawn
in a surface, possibly with crossings, and also a \he
$\eta$ of $G$ in another cubic graph $H$; and we wish to show
that the drawing of $G$ can be extended to a drawing of $H$ 
without introducing any more crossings.
For this to be true, one necessary condition is that for each
$\eta$-bridge $B$, all its attachments belong to the same
``region'' of $G$.
Each region of the drawing
is bounded either by a circuit (if no crossings involve any edge
incident with the region) or by one or more paths, whose
first and last edges cross others and no internal edges cross others.
For instance, in Figure 2, one region is bounded by the path
$6\d 1\d 2\d 3\d 8$; and another by two paths
$6\d 1\d 13\d 18$ and
$15\d 20\d 8\d 3$.
If we list all these circuits and paths we obtain some set of
subgraphs of $G$, and it is convenient to work with this set rather
than explicitly with regions of a drawing of $G$.

Sometimes, the drawing is just of a subgraph $G'$ of $G$ rather
than of all of $G$, and therefore 
all the circuits and paths in the set are
subgraphs of $G'$.
In this case we shall always be able to arrange that
$\eta (e)$ has only one edge, for every edge $e$ of $G$ not in $G'$.
This motivates the following definition.

We say $(G,F, {\cal C})$ is a {\em framework} if $G$ is cubic,
$F$ is a subgraph of $G$,
and
${\cal C}$ is a set of subgraphs of $G\setminus E(F)$, satisfying
(F1)--(F7) below. We say distinct edges $e,f$ are {\em twinned} if there exist
distinct $C_{1},C_{2} \in {\cal C}$ with $e,f \in E(C_{1} \cap C_{2})$.
\begin{itemize}

\item[{\bf(F1)}]
Each member of ${\cal C}$ is an induced subgraph of $G\setminus E(F)$, with at least
three edges, and is either a path or a circuit.

\item[{\bf(F2)}] Every edge of $G\setminus E(F)$ belongs to some member of ${\cal C}$, and for every
two edges $e,f$ of $G$ with a common end not in $V(F)$, there
exists $C \in {\cal C}$ with $e,f \in E(C)$.

\item[{\bf(F3)}] If $C_{1},C_{2} \in {\cal C}$ are distinct and
$v \in V(C_{1} \cap C_{2})$, then either
$V(C_{1} \cap C_{2}) = \{v\}$, or $v$ is incident with an edge in
$C_{1} \cap C_{2}$, or $v \in V(F)$.

\item[{\bf(F4)}] If $C_{1} \in {\cal C}$ is a path, then every member of ${\cal C}$ containing an end-edge of $C_1$ is a path. Moreover, if also
$C_2\in {\cal C}\setminus \{C_1\}$ is a path, then every component of $C_1\cap C_2$ contains an end of $C_1$, 
and every edge of $C_1\cap C_2$ is an end-edge of $C_1$.

\item[{\bf(F5)}] If $C \in {\cal C}$ is a circuit then
$|V(C \cap F) | \leq 1$, and every vertex in
$C \cap F$ has degree $1$ in $F$; and if $C\in {\cal C}$ is a path then
every vertex in $C \cap F$ is an end of $C$ and has degree $0$ or $2$
in $F$.

\item[{\bf(F6)}] If $e,f$ are twinned and $C \in {\cal C}$ with $e\in E(C)$, then $|V(C)|\le 6$, and either
\begin{itemize}
\item $f\in E(C)$, and $C$ is a circuit, and $e,f$ have a common end in
$V(F)$, and no path in ${\cal C}$ contains any vertex of $e$ or $f$, or
\item $f\in E(C)$, and $C$ is a path with end-edges $e,f$, and $C\cap F$ is null, or
\item $f\notin E(C)$, and $C$ is a path with $|E(C)| = 3$, and $e$ is an end-edge of $C$, 
and no end of $e$ belongs to $V(F)$.
\end{itemize}

\item[{\bf(F7)}] Let $C\in \cal{C}$ be a path of length five, with twinned end-edges $e,f$.
Then $|E(C')| \leq 4$ for every path $C' \in {\cal C}\setminus \{C\}$ containing $e$.
Moreover, let $C$ have vertices $v_0\d v_1\c v_5$ in order; then
there exists $C' \in {\cal C}$ with end-edges $e$ and $f$
and with ends $v_{0}$ and $v_{4}$.

\end{itemize}

We will prove a theorem that says, roughly, that if we have a framework $(G,F, {\cal C})$, and a \he of $G$ in $H$,
where $H$ is appropriately cyclically connected, then either the drawing of $G$ extends to an drawing of the whole of $H$,
or there is some bounded enlargement of $\eta(G)$ in $H$ to which the drawing does not extend, and this enlargement
still has high cyclic connectivity.

These seven axioms are a little hard to digest, and before we go on it may help to see how they will be used.
In all our applications of (F1)--(F7) we have some particular graph $G$ in mind and a drawing of it
that defines the framework. 
We could replace (F1)--(F7) just by the hypothesis that $(G,F, {\cal C})$ arises from one of these particular cases, but there
are nine of these cases, and it seemed clearer to try to abstract the properties that we really use.
Here are three examples that might help.
\begin{itemize}
\item 
The simplest application
is to prove \ref{nonplanar0}; we take $G$ to be Dodecahedron, and $F$ null, and $\cal C$ to be the set of region-bounding circuits in the drawing of $G$ in Figure 3.
Suppose now some $H$ contains $G$; our result will tell us that either the embedding of $G$ extends to an embedding of $H$ (and hence $H$ is planar),
or $H$ contains a non-planar subgraph, a bounded enlargement of $\eta(G)$ with high cyclic connectivity. We enumerate all the possibilities
for this enlargement, and check they all contain one of Petersen, Ruby, Box, Triplex. From this, \ref{nonplanar0} will follow.
\item When we come to try to understand the graphs that contain Jaws and not Petersen, we take $G$ to be Jaws, and $(G,F, {\cal C})$
to be defined by the drawing in Figure 2. Thus, $F$ is null; $\cal {C}$ will contain the seven circuits in Figure 2 that bound regions
and do not include any of the four edges that cross, together with eight paths (four like $6 \d 1 \d 2 \d 3 \d 8$; 
two like $1\d 6 \d 5 \d 4 \d 3 \d8$; and two like $6\d 1\d 13\d 18$.) 
\item A last example, one with $F$ non-null; when we prove \ref{arched0}, we take $G$ to be Box, and $(G,F, {\cal C})$ to be defined by the drawing
in Figure 3, and $E(F) =\{f\}$ where $f$ is the edge 13-14. 
In this case, take the drawing of Box given in Figure 3, and delete the edge
$f$, and we get a drawing of $G\setminus f$ without crossings; let $\cal{C}$ be the set of circuits that bound regions
in this drawing. The only twinned edges are 2-13 with 5-13, and 8-14 with 
11-14.
\end{itemize}

(F1)--(F7) have a number of easy consequences, for instance, the following four results.

\begin{thm}\label{frame1}
Let $(G,F, {\cal C})$ be a framework. 
\begin{itemize}
\item $F$ is an induced subgraph of $G$.
\item Let $e\in E(G)\setminus E(F)$. Then $e$ belongs to at least two members of $\cal{C}$, and to more than two if and only if $e$
is an end-edge of a path in ${\cal C}$ and neither end of $e$ is in
$V(F)$; and in this case $e$ belongs to exactly four members of ${\cal C}$,
all paths, and it is an end-edge of each of them.
\item For every two edges $e,f$ of $G$ with a common end with degree three in $G\setminus E(F)$, there
is at most one $C \in {\cal C}$ with $e,f \in E(C)$.
\end{itemize}
\end{thm}
\Proof
Let $e=uv$ be an edge of $E(G)\setminus E(F)$.
We claim that $|\{u,v\} \cap V(F)| \leq 1$.
For by (F2) there exists $C \in {\cal C}$ with $e \in E(C)$.
If $C$ is a circuit the claim follows from (F5), and if $C$
is a path then one of $u,v$ is internal to $C$, and again it
follows from (F5). Thus the first claim holds.

For the second claim, again let $e=uv$ be an edge of $E(G)\setminus E(F)$.
We may assume that $u\not\in V(F)$.
Let $u$ be incident with $e,e_{1},e_{2}$.
By (F2) there
exist $C_{1},C_{2} \in {\cal C}$ with
$e,e_{i} \in E(C_{i})\; (i=1,2)$.
Hence $C_{1} \neq C_{2}$, so $e$ belongs to at least two members of
${\cal C}$.

No other member of ${\cal C}$ contains $e$ and either $e_{1}$
or $e_{2}$, by (F6), since $u \not\in V(F)$.
Hence every other $C \in {\cal C}$ containing $e$ is a path with one end $u$.
If $e$ is not an end-edge of any path in ${\cal C}$ the second claim
is therefore true, so we assume it is.
Hence by (F4), $C_{1}$ and $C_{2}$ are both paths with end-edge $e$,
and both have one end $v$.
If $v \in V(F)$, there is no path in ${\cal C}$ containing $e$
with one end $u$, by (F5), so we may assume that
$v \not\in V(F)$.
Let $v$ be incident with $e,e_{3},e_{4}$; then by (F2)
there exist $C_{3},C_{4} \in {\cal C}$ with
$e,e_{i} \in E(C_{i})\; (i=3,4)$; and $C_{3},C_{4}$ both
have one end $u$.
Hence $C_{1},\ldots, C_{4}$ are all distinct, and no other member of
${\cal C}$ contains $e$. This proves the second claim.

For the third claim, let $v\in V(G)$ be incident with edges $e,f,g\in E(G)\setminus E(F)$.
Suppose there exist distinct $C,C'\in \cal C$ both containing $e,f$. 
Thus $e,f$ are twinned. If $C$ is a circuit, then by (F6) $v\in V(F)$, and by (F5) $v$ has degree one in $F$, a contradiction. Thus
$C$ is a path. By (F6) both $e,f$ are end-edges of $C$, and hence $C$ has length two, a contradiction. 
This proves the third claim, and hence proves
\ref{frame1}.~\bbox

\begin{thm}\label{frame2}
Let $C_{1},C_{2} \in {\cal C}$ be distinct.
Then $|E(C_{1} \cap C_{2}) | \leq 2$, and if equality holds, then either
\begin{itemize}
\item $C_1,C_2$ are both circuits, and 
$C_{1} \cap C_{2}$ is a $2$-edge path with middle vertex $v$ in
$V(F)$, and $v$ has degree one in $F$, or
\item
$C_1,C_2$ are both paths with the same end-edges $e,f$ say, and 
$C_{1} \cap C_{2}$ consists of the disjoint edges $e,f$ and their ends,
and $C_1,C_2$ are disjoint from $F$.
\end{itemize}
\end{thm}
\Proof
Let $e,f \in E(C_{1} \cap C_{2} )$ be distinct.
If $C_{1}$ is a path then by (F6) and (F4), so is $C_{2}$,
and both $C_{1}$ and $ C_{2}$ have end-edges $e,f$, and no end of $e$
or $f$ is in $V(F)$, and by (F5) $C_1,C_2$ are disjoint from $F$.
But then by (F6) $|E(C_{1} \cap C_{2}) | =2$
(for any third edge in $E(C_{1} \cap C_{2} )$ would also have to be an
end-edge of $C_{1}$, which is impossible); and if
$v \in V(C_{1} \cap C_{2} )$ is not incident with $e$ or $f$, then $v$ is
internal to both paths and hence is incident with an edge of
$C_{1} \cap C_{2}$, a contradiction.
Thus in this case the theorem holds.
We may assume then that $C_{1}$ and $C_{2}$ are both circuits.
By (F6), $e,f$ have  a common end, $v$ say, in $V(F)$.
By (F5) no other vertex of $C_{1}$ or $C_{2}$ is in $V(F)$,
and $v$ has degree one in $F$.
By (F6), $E(C_{1} \cap C_{2}) = \{e,f\}$, and hence the theorem holds.
This proves \ref{frame2}.~\bbox

\begin{thm}\label{frame3}
Let $C_{1},C_{2} \in {\cal C}$ be distinct with
$|E(C_{1} \cap C_{2}) | \geq 2$.
Then $|E(C_{1})| \geq 4$.
\end{thm}
\Proof
Suppose that $C_{1}$ is a circuit.
If $|E(C_{1})| =3$, then since $C_{2}$ is an induced subgraph of $G\setminus E(F)$ and
$|E(C_{1} \cap C_{2}) | \geq 2$ it follows that $C_{1}$
is a subgraph of $C_{2}$ which is impossible.
Hence the result holds if $C_{1}$ is a circuit.
Now let $C_{1}$ be a path.
Let $e,f \in E(C_{1} \cap C_{2})$ be distinct; then by (F6), $e$
and $f$ are end-edges of $C_{1}$, and by (F4) $C_{2}$ is a path with
end-edges $e,f$.
Hence again $C_{1}$ is not a subgraph of $C_{2}$, and so since $C_{2}$
is an induced subgraph of $G\setminus E(F)$ it follows that $|E(C_{1}) | \geq 4$.
This proves \ref{frame3}.~\bbox

\begin{thm}\label{uniquetwin}
Let $(G,F, {\cal C})$ be a framework, and let
$e,f_{1},f_{2} \in E(G)$ be distinct.
If $e,f_{1}$ are twinned then $e,f_{2}$ are not twinned.
\end{thm}
\Proof
Let $C_{1},C_{1}' \in {\cal C}$ be distinct with
$e,f_{1} \in E(C_{1} \cap C_{1}')$, and suppose that there
exist $C_{2}, C_{2}' \in {\cal C}$, distinct, with
$e,f_{2} \in E(C_{2} \cap C_{2}')$.
At least three of $C_1,C'_1,C_2,C'_2$ are distinct, and they all contain $e$, and so by 
\ref{frame1} all of $C_1,C'_1,C_2,C'_2 $ are paths and $e$ is an end-edge of each of them.
By (F6) $C_1$ has end-edges $e$ and $f_{1}$, and
$f_{2} \not\in E(C_{1})$.
Since $e,f_{1} \in E(C_{1})$, by \ref{frame3}
$|E(C_{1})| \geq 4$; but since
$f_{2} \not\in E(C_{1})$, by (F6)
$|E(C_{1}) | \leq 3$, a contradiction.
This proves \ref{uniquetwin}.~\bbox

Let $F,G,H$ be graphs, where $F$ is a subgraph of $G$, and let $\zeta, \eta$ be homeomorphic embeddings of $F,G$ 
into $H$ respectively. We say that $\eta$ {\em extends} $\zeta$ if $\eta(e)= \zeta(e)$ for all $e \in E(F)$ and
$\eta (v) = \zeta (v)$ for all $v \in V(F)$.

Let $(G,F, {\cal C})$ be a framework, let $\eta_F$ be a \he of $F$ into $H$, and
let $J$ be the subgraph of $F$ obtained by deleting all
vertices with degree one in $F$.
Let $G'$ be a cubic graph with $J$ a subgraph of $G'$.
A \he $\eta$ of $G'$ in $H$ is said to {\em respect}
$\eta_F$ if $\eta$ extends the restriction of $\eta_F$ to $J$.

Again, let $(G,F, {\cal C})$ be a framework, and let $\eta_F$ be a \he of $F$ into $H$.
We list a number of conditions on the framework, $H$ and $\eta_F$ that we shall
prove have the following property. Suppose that these conditions are satisfied, and there is
a \he of $G$ in $H$ extending $\eta_F$; then the natural drawing of $G\setminus E(F)$ (where the members of
$\cal C$ define the region-boundaries) can be extended to one of $H\setminus E(\eta_F(F))$.
They are the following seven conditions (E1)--(E7). 

\begin{itemize}
\item[{\bf(E1)}]
$H$ is cubic and cyclically four-connected, and if
$(G,F, {\cal C})$ has any twinned edges, then $H$ is cyclically
five-connected.
Also, 
$\eta_{F} (e)$ has only one edge for every $e \in E(F)$.
\end{itemize}

\begin{itemize}
\item[{\bf(E2)}]
Let $e,f \in E(G)\setminus E(F)$ be distinct.
If there is a \he of $G+(e,f)$ in $H$ respecting $\eta_F$,
then there exists $C \in {\cal C}$ with
$e,f \in E(C)$.
\end{itemize}

If $e,f,g$ are distinct edges of $E(G)$ such that no member of ${\cal C}$
contains all of $e,f,g$, but one contains $e,f$, one contains $e,g$
and one contains $f,g$, we call $\{e,f,g\}$ a {\em trinity}.
A trinity is {\em diverse} if every two edges in it are diverse in $G\setminus E(F)$.
\begin{itemize}

\item[{\bf(E3)}]
For every diverse trinity $\{e,f,g\}$ there is no \he of
$G+(e,f,g)$ in $H$ extending $\eta_F$.

\item[{\bf(E4)}]
Let $v$ have degree one in $F$, incident with
$g \in E(F)$.
Let $C_{1},C_{2}$ be the two members of ${\cal C}$ containing $v$.
For all $e_{1} \in E(C_{1})  \setminus  E(C_{2})$ and
$e_{2} \in E(C_{2} )  \setminus E(C_{1})$ such that $e_{1}$ and $e_{2}$ have
no common end, there is no \he of
$G+(e_{1},g)+(e_{2},vy)$ in $H$ respecting $\eta_F$, where
$G+(e_{1},g)$ has new vertices $x,y$.

\item[{\bf(E5)}]
Let $v$ have degree one in $F$, incident with
$g \in E(F)$.
Let $u$ be a neighbour of $v$ in $G\setminus E(F)$, and let $C_{0}$ be the (unique, by \ref{frame1})
member of ${\cal C}$ that contains $u$ and not $v$.
Let $u$ have neighbours $v,w_{1},w_{2}$.
Let $G'=G+(uw_{1},g)$ with new vertices $x_{1},y_{1}$; and
let $G'' = G'+(uw_{2},vy_{1})$ with new vertices $x_{2},y_{2}$.
Let $i=1$ or $2$, and let $e=ux_{i}$.
Let $f$ be an edge of $C_{0}$ not incident with $w_{1}$ or $w_{2}$,
and with no end adjacent to $w_{i}$.
(This is vacuous unless $|E(C_{0})| \geq 6)$.)
There is no \he of $G''+(e,f)$ in $H$ respecting $\eta_F$.
\end{itemize}

Two edges of $G\setminus E(F)$ are {\em distant} if they are diverse in $G$ and not twinned.
Let $C \in {\cal C}$.
We shall speak of a sequence of vertices and/or edges of $C$ as being
{\em in order} in $C$, with the natural meaning (that is, if $C$
is a path, in order as $C$ is traversed from one end, and if $C$
is a circuit, in order as $C$ is traversed from some starting point).
\begin{itemize}
\item If $e,f,g,h$ are distinct edges of $C$, in order, and $e,g$
are distant and so are $f,h$, we call
$G+(e,g)+(f,h)$ 
a {\em cross extension}
({\em of $G$, over $C$}) {\em of the first kind}.
\item If $e,uv,f$ are distinct edges of $C$, and either $e,u,v,f$ are in order, or $f,e,u,v$ are in order,
and $e,uv$ are distant and so are $uv,f$, we call
$G+(e,uv)+(uy,f)$ a {\em cross extension of the second kind},
where $G+(e,uv)$ has new vertices $x,y$.
\item If $u_{1}v_{1}$ and $u_{2}v_{2}$ are distant edges of $C$ and
$u_{1}, v_{1},u_{2},v_{2}$ are in order, we call
$G+(u_{1}v_{1},u_{2}v_{2})+ (xv_{1},yv_{2})$ a
{\em cross extension of the third kind}, where
$G+(u_{1}v_{1},u_{2}v_{2})$ has new vertices $x,y$.
\end{itemize}

\begin{itemize}
\item [{\bf(E6)}]
For each $C \in {\cal C}$ and every cross extension $G'$ of $G$ over $C$
of the first, second or third kinds, there is no \he of
$G'$ in $H$ extending $\eta_F$.
\item [{\bf(E7)}]
Let $C \in {\cal C}$ be a path with $|E(C)|=5$, with vertices
$v_{0}\c  v_{5}$ in order, and let $v_{0}v_{1}$ and
$v_{4}v_{5}$ be twinned.
Let $G_{1}=G+(v_{0}v_{1},v_{4}v_{5})$ with new vertices
$x_{1},y_{1}$; let $G_{2}=G_{1}+(v_{1}v_{2},y_{1}v_{5})$ with
new vertices $x_{2},y_{2}$; and let
$G_{3}=G_{2}+(v_{0}x_{1},y_{2}v_{5})$.
There is no \he of $G_{3}$ in $H$ extending $\eta_F$.
\end{itemize}

In the proofs to come, when we need to apply (E1)--(E7), it is often cumbersome to 
indicate the full \he involved, and we use some shortcuts. For instance,
when we apply (E2), with $e,f,\eta$ as in (E2), let $g$ be the new edge 
of $G+(e,f)$, and let $H'$ be the graph obtained from $\eta(G+(e,f))$ by deleting the interior of the path
$\eta(g)$; we normally say
``by (E2) applied to $H'$ with edges $e,f$'', and leave the reader to figure out the appropriate 
\he and the path $\eta(g)$.

Whenever we wish to apply our main theorem, we have to verify directly that
(E1)--(E7) hold, and this can be a lot of case-checking. We have therefore
tried to design (E1)--(E7) to be as easily checked as possible consistent 
with implying the main result. Nevertheless, there is still a great deal
of case-checking, and we have omitted almost all the details. We are
making available in~\cite{RSTappendix} both the case-checking and all the 
graphs of the paper in computer-readable form.

\section{Degenerate trinities}

Now (E3) was a statement about diverse trinities; our first
objective is to prove the same statement about non-diverse trinities.

A trinity is a {\em Y-trinity} if some two edges in it (say $e$ and $f$)
have a common end $u$, the third edge in it ($g$ say) is not incident
with $u$, and if $h$ denotes the third edge incident with $u$
then there exist $C_{1},C_{2} \in {\cal C}$ with
$e,g,h\in E(C_{1})$ and $f,g,h \in E(C_{2})$.
(Consequently $g,h$ are twinned.)
It is {\em circuit-type} or {\em path-type} depending
whether $g$ and $h$ have a common end or not.

\begin{thm}\label{pathtrinity}
Let $(G,F, {\cal C})$ be a framework and let 
$H, \eta_{F}$ satisfy $(E1)$--$(E7)$.
For every path-type $Y$-trinity $\{e,f,g\}$ there is no
\he of $G+(e,f,g)$ in $H$ extending $\eta_F$.
\end{thm}
\Proof
Let $u,h,C_{1},C_{2}$ be as above.
Since the twinned edges $g,h$ have no common end, it follows from (F6)
that $C_{1}$ and $C_{2}$ are both paths with end-edges $g,h$, and both are vertex-disjoint from $F$.
Let $e=uw_{1},f=uw_{2}$.
Suppose that $\eta$ is a \he of $G$ into $H$ extending
$\eta_F$, and $e,f,g$ are all $\eta$-attachments of
some $\eta$-bridge $B$.

By \ref{frame3}, $|E(C_{1}) | \geq 4$, and so $g$ is not incident with $w_{1}$,
and similarly not with $w_{2}$.
By (F7), at least one of $C_1,C_2$ has length at most four, and so we may assume that the edges of
$C_1$ in order are $h,e,g_1, g$ say.
Let $\eta '$ be obtained from $\eta$ by
rerouting $g_1$ along an $\eta$-path in $B$ from $\eta (g)$ to $\eta (e)$.
Then $\eta '$ extends $\eta_F$, and $g_1$ and $f$ are both
$\eta'$-attachments of an $\eta'$-bridge.
By (E2) applied to $\eta'(G)$ with edges $g_1,f$, there exists $C \in {\cal C}$ with $g_1,f \in E(C)$,
and hence with $e \in E(C)$ since $C$ is an induced subgraph of $G\setminus E(F)$.
But then $e,g_1 \in E(C \cap C_{1})$, and $C_{1} \neq C$, 
so $e,g_1$ are twinned edges, and yet their common end $w_1$
is not in $V(F)$, contrary to (F6).
There is therefore no such $\eta$.
This proves \ref{pathtrinity}.~\bbox

Let $\{e,f,g\}$ be a circuit-type $Y$-trinity, where
$e=xw_{1}$, $f=xw_{2}$ and $g=vw_{3}$, where $v,w_3\ne x$ and $v,x$ are adjacent in $G$.
Let $h = vx$, and let $w_4$ be the third neighbour of $v$. Since $g,h$ are twinned and share an end, \ref{frame1}
implies that $vw_4\in E(F)$. Hence $w_4\ne w_1, w_2$, since no member of $\cal C$ contains both $v,w_4$.
(See Figure 4.)
We wish to consider three rather similar graphs
$G_{1},G_{2},G_{3}$ called {\em expansions} of the
$Y$-trinity $\{e,f,g\}$.
Let $G'$ be obtained from $G$ by deleting $x$ and the edge $vw_3$,
and adding five new vertices $x_{1},x_{2},x_{3},y_{1},y_{2}$
and nine new edges $x_{1}w_{1},x_{2}w_{2},x_{3}w_{3}$ and
$x_{i}y_{j}$ for $1 \leq i \leq 3$ and $1 \leq j \leq 2$.
Let $G_{1},G_{2},G_{3}$ be obtained from $G'$ by deleting
the edge $y_{2} a$ (where $a$ is $x_{1},x_{2}$ and $x_{3}$
respectively), and adding two new edges
$vy_{2},va$. Let $x_4= v$. (The reason we did not just
replace $v$ by a new vertex $w_4$, is that the edge $vw_4$ belongs to $F$ and we want to preserve it.) 
Thus $F$ is a subgraph of $G_1,G_2$ and $G_3$.
(See Figure 4.)

\begin{figure} [h!]
\centering
\includegraphics{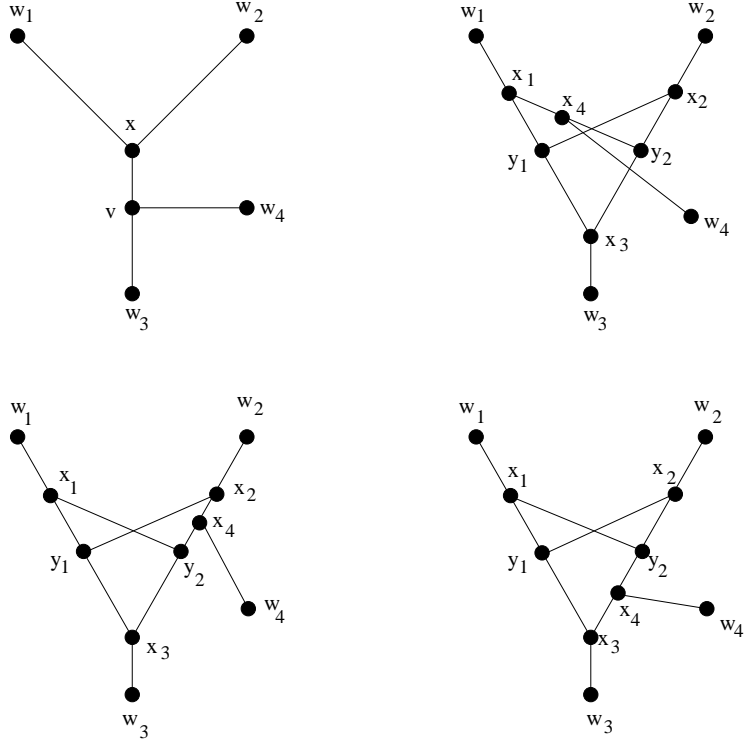}
\caption{A circuit-type $Y$-trinity, and its three expansions.}
\label{trinityfig}
\end{figure}

\begin{thm}\label{trinity2}
Let $(G,F, {\cal C})$ be a framework, and let $H, \eta_{F}$
satisfy $(E1)$--$(E7)$.
Let $\{e,f,g\}$ be a circuit-type $Y$-trinity,
and $G_{1},G_{2},G_{3}$ its three expansions.
Then there is no \he of $G_1,G_2$ or $G_3$ in $H$
extending $\eta_F$.
In particular, there is no \he of $G+\{e,f,g\}$ in $H$
extending $\eta_F$.
\end{thm}
\Proof
Let $v,x,w_{1},\ldots,w_{4}$ be as in Figure 4 and let $G_{1},G_{2},G_{3}$ be
labelled as in Figure 4, where $e=xw_1$, $f=xw_2$, and $g=vw_3$. 

Suppose that there is a \he $\eta$ of some $G_{k}$ in $H$
extending $\eta_F$. Let $A$ be the subgraph of $G_{k}$ induced on
$\{x_{1},x_{2},x_{3},x_{4},y_{1},y_{2} \}$, and
$B$ the subgraph of $G_{k}$ induced on the complementary set of vertices. It follows that
there is a \he $\zeta$ of $G_k$ in $H$ such that:
\begin{itemize}
\item $\zeta$ extends the restriction of $\eta$ to $B$ (and in particular, $\zeta(z) = \eta(z)$ for every vertex or edge $z$
of $F$ different from $x_4,w_4x_4$)
\item $\zeta(w_4x_4)$ is a path with one end $\eta(w_4)$ containing the one-edge path $\eta(w_4x_4)$.
\end{itemize}
(To see this, take $\zeta = \eta$.) Let $Z_i = \zeta(x_{i}w_{i})$ for $i = 1\l 4$.
Let us choose $k$ and $\zeta$ such that
\\
\\
(1) {\em $Z_1\cup Z_2\cup Z_3$ is minimal, and subject to that $Z_4$ is minimal.}

\bigskip

Since $H$ is cyclically five-connected by (E1) since there are twinned edges in $G$, there are five disjoint paths
$P_{1},\ldots,P_{5}$ of $H$ from $\zeta (A)$ to $\zeta (B)=\eta(B)$.
Choose $P_1\l P_5$ to minimize the number of edges of $P_1\cup\cdots\cup P_5$ that do not belong to $Z_1\cup\cdots\cup Z_4$.
It follows that 
each $P_i$ has only its first vertex $a_i$ say in $ V(\zeta(A))$, and only its last
vertex $b_i$ say in $V(\eta (B))$.
Now one of $a_{1},\ldots, a_{5}$ is different from $\zeta(x_{1}), \zeta (x_{2})$,
$\zeta(x_{3}), \zeta(x_{4})$, say $a_{5}$. Let $a_{5} \in V(\zeta(h_1))$, where $h_1 \in E(A)$.
From (1) (or the theory of augmenting paths for network flows) it follows easily that
$\{a_{1},\ldots,a_{4}\} = \{\zeta(x_{1}),\ldots, \zeta (x_{4})\}$,
and we may assume that $a_{i} = \zeta (x_{i})\;(1 \leq i \leq 4)$.

Let $p$ be the first vertex (that is, closest to $a_{5}$) in $P_{5}$
that belongs to
$\eta (B) \cup Z_1\cup Z_2\cup Z_3\cup Z_4$
(this exists since $b_{5} \in V(\eta (B))$), and let $P = P_{5} [a_{5},p]$. 
\\
\\
(2) {\em $p \in V(\eta (B))$.}
\\
\\
\Subproof
Suppose not; then $p \in V(Z_i)$ for some $i$.
If $i=4$, then by replacing
$Z_4 [ \zeta (x_{4}),p]$
by $P$ we obtain a \he of some
$G_{k'}$ (where possibly $k' \neq k$), contradicting (1), since
$Z_4$ is replaced by a proper subpath and $Z_1,Z_2,Z_3$ remain unchanged.
So $1\le i\le 3$.

If $h_1$ is incident with $x_{i}$, then by rerouting $h_1$ along $P$ we obtain a contradiction to (1).
Now suppose that $h_1=ab$ where $a$ is adjacent to $x_{i}$.
By rerouting $ax_{i}$ along $P$, we again obtain
a contradiction to (1).

Thus, neither end of $h_1$ is adjacent to $x_{i}$.
Consequently, $h_1 \neq y_{2}x_{4}$, and $y_1$ is not incident with $h_1$, since $1\le i\le 3$.
The only remaining possibility is that there is a four-vertex
path of $G_k$ with vertices $x_{i},a,b,x_{j}$ in order, for some $j \neq i$, where
$\{a,b\}= \{y_{2},x_{4}\}$, and $h_1=bx_{j}$.
But then there is a \he of some $G_{k'}$ in $H$ mapping $G_{k'}$ to the graph obtained
from $\zeta(G)\cup P$ by deleting the interior of $\zeta (x_{i}a)$, contradicting (1).
This proves (2).

\bigskip
Hence $p \in V(\eta (h_{2}))$ for some $h_{2} \in E(B)$.
Now we examine the possibilities for $h_{1}$ and $h_{2}$.
Since $\eta (h_{2})$ has an interior vertex, it follows from the choice of $\zeta$ that
$h_2\notin E(F)$.
We recall that $v\in V(G)\cap V(F)$.
Let $C_{1},C_{2} \in {\cal C}$ be the two members of ${\cal C}$
that contain $v$, and let $C_{0} \in {\cal C}$ contain $e$ and $f$.
Thus $C_0, C_1, C_2$ are circuits by (F6), and $v$ is the only vertex of $F$ in $V(C_1 \cup C_2)$.
\\
\\
(3) {\em $h_2$ belongs to at most one of $C_0,C_1,C_2$.}
\\
\\
\Subproof
By \ref{frame2}, $E(C_1\cap C_2)$ contains at most two edges, and since it contains both $g,vx$, it follows that
$h_2\notin E(C_1\cap C_2)$. Since $C_1$ is a circuit and $v\in V(F)$, (F5) implies that $x,w_1\notin V(F)$, and so
neither end of $xw_1$ is in $V(F)$. Since $xw_1\in E(C_0\cup C_1)$, \ref{frame2} implies that $|E(C_0\cap C_1)|=1$ and so
$h_2\notin E(C_0\cap C_1)$; and similarly $h_2\notin E(C_0\cap C_2)$. This proves (3).
\\
\\
(4) {\em $k = 1$ or $2$.}
\\
\\
\Subproof
Suppose that $k=3$.
First, suppose that $h_{1}$ is incident with $y_{1}$.
By restricting $\zeta$ to $G_{3} \setminus y_{1}$ we obtain
a \he $\eta'$ of $G$ in $H$ respecting $\eta_F$,
such that $e,f,g$ and $h_{2}$ are all $\eta'$-attachments in $E(G)\setminus E(F)$
of some $\eta'$-bridge.
Since $C_1,C_2$ are the only members of $\mathcal{C}$ containing $g$, it follows
from (E2), applied to $\eta'(G)$ with the edges $g,h_2$, that 
$h_{2} \in E(C_{1} \cup C_{2})$.
Since $C_{1}$ and $C_{0}$ are the only members of ${\cal C}$ containing
$e$ it follows from (E2) (with the edges $e,h_2$) 
that $h_{2} \in E(C_{0} \cup C_{1} )$, and similarly
$h_{2} \in E(C_{0} \cup C_{2} )$.
Thus $h_{2}$ belongs to two of $C_{0},C_{1},C_{2}$, contrary to (3).
This proves that $h_{1}$ is not incident with $y_{1}$.

Suppose next that 
$h_{1}$ is incident with $y_{2}$. By restricting $\eta$ to
$G_{3} \setminus y_{2}$ we obtain a \he $\eta'$ of $G$
in $H$ respecting $\eta_F$ such that $e,f$ and $h_{2}$ are all
$\eta'$-attachments of some $\eta'$-bridge. So 
$h_{2} \in E(C_{0} \cup C_{1})$, by (E2) applied to $\eta'(G)$ with edges $e,h_2$,
and similarly 
$h_{2} \in E(C_{0} \cup C_{2})$. By (3) it follows that $h_2\in E(C_0)$, and $h_2\notin E(C_1\cup C_2)$.
Let $H'$ be the graph obtained from
$\zeta(G_3)$ by deleting the interiors of $\zeta(x_1y_2)$ and $\zeta(x_3y_1)$. 
There is a \he of $G$ in $H$ respecting $\eta_F$, mapping $G$ onto $H'$; and from (E2) applied to
$H'$ with edges $f,h_2$, we deduce that $h_2\in E(C_1\cup C_2)$, a contradiction.
This proves that $h_1$ is not incident with $y_2$.

Thus, $h_{1}=x_{3}x_{4}$.
From (E2) applied to the restriction of $\zeta$ to
$G_{3}\setminus y_{1}$ and the edges $g,h_2$, it follows that $h_{2} \in E(C_{1} \cup C_{2})$;
and from the symmetry between $C_1,C_2$, we may assume that $h_2\in E(C_2)$ without loss of generality.
By \ref{frame2}, $w_1\notin V(C_2)$, and it follows that $h_2,e$ are disjoint edges of $G$.
From (E4)
applied to the restriction of $\zeta$ to
$G_{3} \setminus y_{2}$, we obtain from the paths
$\zeta (x_{1}y_{2}) \cup \zeta (x_{4}y_{2})$ and $P$ that
$h_{2} \not\in E(C_{2})$, a contradiction.
This proves (4).

\bigskip

From (4) and the symmetry between $w_{1}$ and $w_{2}$ (exchanging
$G_{1}$ and $G_{2}$) we may therefore assume that $k=1$. There are three homeomorphic embeddings of $G$ in $H$ respecting $F$ that we need:
\begin{itemize}
\item let $H_1$ be the graph obtained from $\zeta(G_1)$ by deleting the interiors of $\zeta(x_1x_4)$
and $\zeta(x_3y_1)$
\item let $H_2$ be obtained from $\zeta(G_1)$ by deleting the interiors of $\zeta(x_1x_4)$ and $\zeta(x_2y_2)$
\item let $H_3$ be obtained from $\zeta(G_1)$ by deleting the interiors of $\zeta(x_3y_1)$ and $\zeta(x_2y_2)$.
\end{itemize}
For $i = 1,2,3$ there is a \he $\eta_i$ of $G$ in $H_i$ respecting $F$, with $\eta_i(z) = \eta(z)$ for each vertex and edge $z$ of $B$.
\\
\\
(5) {\em $h_2\in E(C_0\cup C_1)$.}
\\
\\
\Subproof Suppose not. 
By (E2) applied to $H_1$ and the edges $e,h_2$, it follows that 
$$h_1\ne x_1y_1,x_3y_1,x_1x_4,x_2y_1,$$
and so $h_1$ is incident with $y_2$. 
By (E2) applied to $H_3$ and the edges $g,h_2$, we deduce that $h_2\in E(C_2)$.
Consequently $e,h_2$ are disjoint, since $w_1\notin V(C_2)$; but then this contradicts
(E4)
applied to $H_2$ and the paths
$\zeta(x_1x_4)$ and $P$ (extended by a subpath of $\zeta(x_2y_2)$ if necessary).
\\
\\
(6) {\em  $h_2\in E(C_0\cup C_2)$.}
\\
\\
\Subproof
Suppose not.
By (E2) applied to $H_3$ and the edges $f,h_2$, it follows that 
$$h_1\ne x_1y_1,x_2y_1,x_3y_1,x_2y_2,$$
and so $h_1$ is one of $x_1x_4,x_4y_2,x_3y_2$. By (5), $h_2\in E(C_1)$, and so $f,h_2$ are disjoint, since $w_2\notin V(C_1)$.
But this contradicts (E4)
applied to $H_2$ and the paths
$\zeta(x_2,y_2)$ and $P$ (extended by a subpath of $\zeta(x_1x_4)$ if necessary).

\bigskip
From (3) and (6), it follows that $h_2\in E(C_0)$, and $h_2\notin E(C_1\cup C_2)$.
By (E2) applied to $H_3$ and the edges $g,h_2$, we deduce that $h_1\ne x_3y_1,x_3y_2,x_2y_2,x_4y_2$;
and by (E2) applied to $H_3$ and the edges $vx,h_2$, we deduce that $h_1\ne x_1x_4$.
Thus $h_1$ is one of $x_1y_1,x_2y_1$. 

We recall that $\eta_2$ is a \he of $G$ in $H_2$.
Suppose that $h_2$ is incident with $w_1$. 
Let $\eta'$ be obtained from $\eta_2$ by rerouting $e$ along $P$; then the paths 
$\zeta(x_2y_2)$ and $\zeta(x_1w_1)\cup \zeta(x_1x_4)$ violate (E4).
Similarly, if $h_2$ is incident with $w_2$, let $\eta'$ be obtained from $\eta_2$
by rerouting $f$ along $P$; then the paths $\zeta(x_1x_4)$ and $\zeta(x_2y_2)\cup \zeta(x_2w_2)$ violate (E4).

Thus $w_1,w_2$ are not incident with $h_2$. 
Next suppose that $h_1 = x_1y_1$ and one end $a$ say of $h_2$ is adjacent to $w_1$.
Let $\eta'$ be obtained from $\eta_2$ by rerouting $aw_1$ along $P$;
then the paths $\zeta(x_1x_4), \zeta(x_2y_2)$ violate (E4). 
Next suppose that $h_1 = x_2y_1$ and 
one end $a$ of $h_2$ is adjacent to $w_2$.
Let $\eta'$ be obtained from $\eta_2$ by rerouting $aw_2$ along $P$; then the paths
$\zeta(x_1x_4), \zeta(x_2y_2)$ violate (E4). 
In summary, then, we have shown that $h_2\in E(C_0)$,
incident with neither of $w_1,w_2$, and for $i = 1,2$, if $h_1 = x_iy_1$ then no end of $h_2$
is adjacent to $w_i$. But this contradicts (E5).

There is therefore no such $\eta$, and the first statement of the
theorem holds.
The second statement of the theorem follows from the first, since
$G+(e,f,g)$ is isomorphic to $G_{3}$ (and the isomorphism fixes
$F$.) This proves \ref{trinity2}.~\bbox

\begin{thm}\label{trinity3}
Let $(G,F,{\cal C})$ be a framework, and let
$H, \eta_{F}$ satisfy $(E1)$--$(E7)$.
Let $\{e_{1},e_{2},e_{3}\}$ be a trinity such that no vertex is incident
with all of $e_{1},e_{2},e_{3}$.
Then there is no \he of $G+(e_{1},e_{2},e_{3})$ in $H$
extending $\eta_F$.
\end{thm}
\Proof
For $i=1,2,3$ there exists $C_{i} \in {\cal C}$ with
$\{e_{1},e_{2},e_{3}\}  \setminus  \{e_{i}\} \subseteq E(C_{i})$
and $e_{i} \not\in E(C_{i})$, since
$\{e_{1},e_{2},e_{3}\}$ is a trinity.
Suppose first that $e_{1},e_{2}$ have a common end $v$ say; and let
$h$ be the third edge incident with $v$.
By hypothesis $h \neq e_{3}$.
If $v \in V(F)$ then since $v$ has degree two in $C_{3}$, $C_{3}$
is a circuit, and hence by (F4), $e_{1}$ is not an end-edge of $C_{2}$;
and if 
$v \not\in V(F)$ then by (F3) either $e_{1}$ is not an end-edge
of $C_{2}$, or $e_{2}$ is not an end-edge of $C_{1}$, and we may assume the
first.
Hence in either case $e_{1}$ is not an end-edge of $C_{2}$.
Since $e_{1} \in E(C_{2})$ and $e_{2} \not\in E(C_{2})$, it follows
that $h \in E(C_{2})$.
By (F3), since $e_{3} \in E(C_{1} \cap C_{2})$, it follows that
$h \in E(C_{1})$, since $v$ not in $V(F)$ by (F5); and so
$\{e_{1},e_{2},e_{3}\}$ is a $Y$-trinity, contrary to \ref{pathtrinity} and \ref{trinity2}.

Thus, no two of $e_{1},e_{2},e_{3}$ have a common end.
Suppose that there is a \he of $G+(e_{1},e_{2},e_{3})$
in $H$ extending $\eta_F$.
Then there is a \he of $G$ in $H$ extending $\eta_F$,
such that $e_{1},e_{2},e_{3}$ are all $\eta$-attachments of the same
$\eta$-bridge $B$ say.
By (E3), $\{e_{1},e_{2},e_{3}\}$ is not diverse in $G\setminus E(F)$, so we may assume that
$e_{1} =a_{1}b_{1}$ and $e_{2}=a_{2}b_{2}$, where
$a_{1},a_{2}$ are adjacent in $G\setminus E(F)$.
Let $a_{1}a_{2}=e_{0}$.

Since $e_{1},e_{2} \in E(C_{3})$ and $C_{3}$ is an induced subgraph of $G\setminus E(F)$, it follows that
$e_{0} \in E(C_{3})$.
Let $a_{1}$ have neighbours $b_{1},a_{2},c_{1}$ and $a_{2}$
have neighbours $a_{1},b_{2},c_{2}$ in $G$.

Since $e_{0}$ is not an end-edge of $C_{3}$, it is not an end-edge
of $C_{1}$ or $C_{2}$, by (F4).
Since $e_{0}$ and $e_{3}$ are disjoint, and
$e_{3} \in E(C_{1} \cap C_{2})$, it follows from (F6)
that $e_{0} \not\in E(C_{1} \cap C_{2})$; we assume that
$e_{0} \not\in E(C_{1})$ without loss of generality.
Suppose that $e_{0} \in E(C_{2}).$
Since $e_{0},e_{1} \in E(C_{2} \cap C_{3})$, it follows from (F6)
that $C_{2} ,C_{3}$ are both circuits,
$a_{1} \in V(F)$ and $a_{1}c_{1} \in E(F)$.
Hence $a_{2}c_{2} \in E(C_{2})$ (since $e_{2} \not\in E(C_{2})$).
Moreover by (F3),  $a_{2}$ is incident with an edge in
$C_{1} \cap C_{2}$, since $E(C_{1} \cap C_{2}) \neq \emptyset$ and
$a_{2} \in V(C_{1} \cap C_{2})$.
Since $e_{0} \not\in E(C_{1})$ and $e_{2} \not\in E(C_{2})$ it follows
that $a_{2}c_{2} \in E(C_{1})$.
Since $E(C_{1} \cap C_{2})$ contains $e_{3}$ and $a_{2}c_{2}$ and
$C_{2}$ is a circuit, it follows from (F6) that
$c_{2} \in V(F)$, and so
$a_{1},c_{2} \in V(C_{2} \cap F)$ contrary to (F5).
This proves that $e_{0} \not\in E(C_{2})$.

If $a_{1} \in V(F)$ then $a_{1}c_{1} \not\in E(C_{2})$ by (F5),
and so $e_{1}$ is an end-edge of $C_{2}$.
By (F4), $C_{3}$ is a path, and $a_{1} $ is an internal vertex of it,
contrary to (F5).
Hence $a_{1} \not\in V(F)$, and similarly $a_{2} \not\in V(F)$.

Now $e_{1},e_{2},e_{3}$ are all $\eta$-attachments of $B$.
Let $P$ be an $\eta$-path in $B$ with ends in $\eta(e_{1})$ and
$\eta (e_{2})$, and let $\eta '$ be obtained by rerouting $e_{0}$ along $P$.
Then $\eta'$ is a \he of $G$ in $H$ extending $\eta_F$.
Since $e_{3}$ is an $\eta$-attachment of $B$, it follows that
$e_{0}$ and $e_{3}$ are $\eta'$-attachments of some $\eta'$-bridge.
By (E2), applied to $\eta'(G)$ with edges $e_0,e_3$, there exists $C_{4} \in {\cal C}$ with
$e_{0},e_{3} \in E(C_{4})$.
Since $a_{1} \not\in V(F)$ it follows from (F6) that 
$e_{1} \not\in E(C_{4})$.
But from (F4) applied to $C_{3} $ and $C_{4}$, $e_{0}$ is not an end-edge
of $C_{4}$.
By (F3) applied to $C_{2}$ and $C_{4}$,
$a_{1}c_{1} \in E(C_{2} \cap C_{4})$.
Since $E(C_{2} \cap C_{4})$ contains both $a_{1}c_{1}$ and $e_{3}$, it follows
that $e_3, a_1c_1$ are twinned, and similarly so are $e_3, a_2c_2$, contrary to \ref{uniquetwin}.
Thus there is no such $\eta$. This proves \ref{trinity3}.~\bbox

Next we need the following lemma.

\begin{thm}\label{gateleg}
Let $\eta$ be a \he of a cubic graph $G$ in a cyclically four-connected
cubic graph $H$.
Let $v \in V(G)$, incident with edges $e_{1},e_{2},e_{3}$,
and suppose that $e_{1},e_{2},e_{3}$ are $\eta$-attachments of some
$\eta$-bridge.
Then there is a \he $\eta'$ of $G$ in $H$, such that
$\eta' (u) = \eta (u)$ for all
$u \in V(G)  \setminus  \{v\}$, and $\eta' (e) = \eta (e)$
for all $e \in E(G)  \setminus  \{e_{1},e_{2},e_{3}\}$, and such that for some edge
$e_{4} \neq e_{1},e_{2},e_{3}$ of $G,e_{1},e_{2},e_{3},e_{4}$ are
$\eta'$-attachments of some $\eta'$-bridge.
\end{thm}
\Proof
For $1 \leq i \leq 3$, let $e_{i}$ have ends $v$ and $v_{i}$.
Let $G' = G+(e_{1},e_{2},e_{3})$, with new vertices
$x_{1},x_{2},x_{3},w$.
By hypothesis, there is a \he $\eta'$ of $G'$ in $H$ such that
$\eta' (u) = \eta (u)$ for all $u \in V(G) \setminus  \{v\}$, and
$\eta' (e) = \eta (e)$ for all $e \in E(G)  \setminus  \{e_{1},e_{2},e_{3}\}$.
Choose $\eta'$ such that
\[
\eta' (v_{1}x_{1}) \cup \eta' (v_{2}x_{2}) \cup \eta' (v_{3}x_{3})
\]
is minimal.
Since $H$ is cyclically four-connected, there is an $\eta'$-path with
one end in
\[
\bigcup (V(\eta' (vx_{i})) \cup V(\eta' (wx_{i})): 1 \leq i \leq 3 )
\]
and the other end, $t$, in
\[
V(\eta (G \setminus v) \cup \eta' (v_{1}x_{1}) \cup
\eta' (v_{2}x_{2}) \cup \eta' (v_{3}x_{3})) .
\]
From the choice of $\eta'$ it follows that $t$ belongs to none of
$\eta' (v_{1}x_{1})$, $\eta' (v_{2}x_{2})$,
$\eta' (v_{3}x_{3})$, and so it belongs to
$\eta' (e_{4}) = \eta (e_{4})$ for some
$e_{4} \in E(G\setminus v)$.
This proves \ref{gateleg}.~\bbox

\begin{thm}\label{trinity5}
Let $(G,F,{\cal C})$ be a framework, and let $H,\eta_{F}$
satisfy $(E1)$--$(E7)$.
Let $\{e_{1},e_{2},e_{3}\}$ be a trinity.
There is no \he of $G+(e_{1},e_{2},e_{3})$ in $H$
extending $\eta_F$.
\end{thm}
\Proof
By \ref{trinity3} we may assume that $v\in V(G)$ is incident with $e_{1},e_{2}$
and $e_{3}$.
Suppose $\eta$ is a \he of $G+(e_{1},e_{2},e_{3})$ in
$H$ extending $\eta_F$.
By \ref{gateleg} there is an edge $e_{4} \neq e_{1},e_{2},e_{3}$ of $G$
such that there are homeomorphic embeddings of each of
$G+(e_{2},e_{3},e_{4})$, $G+(e_{1},e_{3},e_{4})$,
$G+(e_{1},e_{2},e_{4})$ in $H$ extending $\eta_F$.
It follows that $e_{4} \notin E(F)$.
Since no vertex is incident with all of
$e_{2},e_{3},e_{4}$, it follows from \ref{trinity3} that
$\{e_{2},e_{3},e_{4}\}$ is not a trinity; and yet (E2), applied to $\eta(G)$ with edges each pair of $e_2,e_3,e_4$, implies
that every two of $e_2,e_3,e_4$ are contained in a member of $\mathcal{C}$. Consequently there exists
$C_{1} \in {\cal C}$ with $e_{2},e_{3},e_{4} \in E(C_{1})$.
Similarly there exist $C_{2},C_{3} \in {\cal C}$ with
$e_{1},e_{3},e_{4} \in E(C_{2})$ and
$e_{1},e_{2},e_{4} \in E(C_{3})$.
Since $\{e_{1},e_{2},e_{3}\}$ is a trinity,
$e_{i} \not\in E(C_{i}) (1 \leq i \leq 3 )$, and so
$C_{1},C_{2},C_{3}$ are all distinct.

Now if $e_{4}$ is not the end-edge of any path in ${\cal C}$, then since
$C_{2} \cap C_{3}$ contains $e_{1}$ and $e_{4}$ it follows from
(F6) that $e_{1}$ and $e_{4}$ have a common end, and
similarly so do $e_{i}$ and $e_{4}$ for $i=1,2,3$, which is impossible.
Hence $e_{4}$ is an end-edge of some path in ${\cal C}$.
By (F4) $C_{1},C_{2}$ and $C_{3}$ are all paths.
Since $e_{3},e_{4} \in E(C_{1} \cap C_{2})$, $C_{1}$ has end-edges
$e_{3}$ and $e_{4}$; and since
$e_{2},e_{4} \in E(C_{1} \cap C_{3})$, $C_{1}$ has end-edges
$e_{2}$ and $e_{4}$,
a contradiction.
This proves \ref{trinity5}.~\bbox

\begin{thm}\label{trinity6}
Let $(G,F, {\cal C})$ be a framework, and let $H, \eta_{F}$
satisfy $(E1)$--$(E7)$.
Let $\eta$ be a \he of $G$ in $H$ extending $\eta_F$.
For every $\eta$-bridge $B$ there exists $C \in {\cal C}$ such
that $e\in E(C)$ for every $\eta$-attachment $e$ of $B$.
\end{thm}
\Proof
Since $\eta$ extends $\eta_F$, it follows that 
$Z \subseteq E(G)\setminus E(F)$, where $Z$ is the
set of all $\eta$-attachments of $B$.
Suppose, for a contradiction, that there is no $C \in {\cal C}$ with
$Z \subseteq E(C)$, and choose $X \subseteq Z$ minimal such that there
is no $C \in {\cal C}$ with $X \subseteq E(C)$.
By (F2), $|X| \geq 2$; by (E2), applied to $\eta(G)$ with edges the members of $X$, 
$|X| \neq 2$; and by \ref{trinity5}, $|X| \neq 3$.
Hence $|X| \geq 4$.
Let $X= \{e_{1},\ldots, e_{k}\}$ say, where $k \geq 4$.
For each $i\in \{1\l k\}$, there exists $C\in \cal C$ including $X\setminus \{e_i\}$, from the minimality
of $X$. All these members of $\cal C$ are different, and so every two members of $X$
are twinned, contrary to \ref{uniquetwin}.
This proves \ref{trinity6}.~\bbox

\section{Crossings on a region}

Let $\eta$ extend $\eta_F$, and let $B$ be an $\eta$-bridge. Since $\eta$ extends $\eta_F$, it follows
that no $\eta$-attachment of $B$ is in $E(F)$, and so by \ref{trinity6}, there exists $C\in \cal C$
such that every $\eta$-attachment
of $B$ belongs to $C$. If $C$ is unique, 
we say that $B $ {\em sits on C}.

Our objective in this section is to show that if $\eta$ extends $\eta_F$,
then for every $C \in {\cal C}$ all the bridges that sit on $C$
can be simultaneously drawn within the ``region'' that $C$ bounds.
There may be some bridges that sit on no member of ${\cal C}$, but we
worry about them later.

Let $C$ be a path or circuit in a graph $J$. We say paths $P,Q$ of $J$ {\em cross}
with respect to $C$, if $P,Q$ are disjoint, and $P$ has distinct ends $p_1,p_2\in V(C)$, and $Q$ has distinct ends $q_1,q_2\in V(C)$,
and no other vertex of $P$ or $Q$ belongs to $C$, and these ends can be numbered such that
either $p_1,q_1,p_2,q_2$ are in order in $C$, or $q_1,p_1,q_2,p_2$ are in order in $C$.
We say that $J$ is $C$-{\em planar} if $J$ can be drawn in a closed disc
$\Delta$ such that every vertex and edge of $C$ is drawn in the boundary
of $\Delta$.
We shall prove:

\begin{thm}\label{extend}
Let $(G,F, {\cal C})$ be a framework, and let $H$, $ \eta_{F}$
satisfy $(E1)$--$(E7)$.
Let $\eta$ be a \he of $G$ in $H$ that extends $\eta_F$,
let $C \in {\cal C}$, and let ${\cal A}$ be a set of $\eta$-bridges
that sit on $C$.
Let $J = \eta (C) \cup \bigcup (B\;:\;B \in {\cal A})$.
Then $J$ is $\eta (C)$-planar.
\end{thm}

\ref{extend} is a consequence of the following.

\begin{thm}\label{extendbranch}
Let $(G,F, {\cal C})$ be a framework, and let
$H, \eta_{F}$ satisfy $(E1)$--$(E7)$.
Let $\eta$ be a \he of $G$ in $H$ that extends $\eta_F$, and let
$C \in {\cal C}$.
Let $P$, $Q$ be $\eta$-paths that cross with respect to $\eta(C)$. 
Then for one of $P$, $Q$, the $\eta$-bridge that contains it does not sit on $C$.
\end{thm}

\noindent{\bf Proof of \ref{extend}, assuming \ref{extendbranch}.}

Suppose that $X,Y \subseteq V(J)$ with $X \cup Y=V(J)$ and $V(C) \subseteq Y$,
such that $|X \setminus Y| \geq 2$ and no edge of $J$ has one end in $X \setminus Y$ and the other
in $Y \setminus X$.
We claim that $|X \cap Y | \geq 4$.
For let $Y' = Y \cup (V(H) \setminus X)$;
then no edge of $H$ has
one end in $X \setminus Y'$ and the other in $Y' \setminus X$, and $X \cup Y'=V(H)$,
and $|X \setminus Y'| \geq 2$, and so $X$ and $Y'$ both includes the vertex set
of a circuit of $H$.
Since $H$ is cyclically four-connected, it follows that
$|X \cap Y' | \geq 4$, and so $|X \cap Y | \geq 4$ as claimed.

From this and theorems 2.3 and 2.4 of \cite{RS9}, it follows, assuming
for a contradiction that $J$ is not $\eta (C)$-planar, that there are
$\eta$-paths $P,Q$ in $J$ that cross with respect to $\eta(C)$.
By \ref{extendbranch} the $\eta$-bridge containing one of $P,Q$ does not sit on $C$
and hence does not belong to ${\cal A}$, a contradiction.
This proves \ref{extend}.~\bbox

\bigskip

\noindent{\bf Proof of \ref{extendbranch}.}

We remark, first, that  
\\
\\
(1) {\em If $B$ is an $\eta$-bridge that sits on $C$, and $e\in E(C)$ is an $\eta$-attachment of $B$, then
there is an $\eta$-attachment $g\in E(C)$ of $B$ such that $g\ne e$ and $g$ is not twinned with $e$.}
\\
\\
\Subproof By \ref{frame1} it follows that $B$ has at least two $\eta$-attachments. 
Suppose that every $\eta$-attachment different from $e$ is twinned with $e$;
then by \ref{uniquetwin} there is only one other, say $f$, and $e,f$ are twinned, 
and therefore there exists $C'\ne C$ in ${\cal C}$ containing all $\eta$-attachments of $B$,
contradicting that $B$ sits on $C$. This proves (1).

\bigskip
For $e,f\in E(C)$, let
\[
\epsilon (e,f) = \left\{
\begin{array}{ll}
3 & {\rm if}\, e=f, \\
2 & {\rm if}\, e \neq f,\, {\rm and}\, e,f \, {\rm are\,twinned} \\
0 & {\rm if}\, e \neq f,\, {\rm and}\, e,f\, {\rm are\,not\,twinned.}
\end{array}
\right.
\]
Let $P$ have ends $p_1,p_2$, and let $Q$ have ends $q_1,q_2$; and let $B_1, B_2$ be the $\eta$-bridges containing $P,Q$ respectively.
Let $p_{i} \in V(\eta (e_{i}))$ and $q_{i} \in V(\eta (f_{i}))$ for $i = 1,2$, and let 
$N = \epsilon (e_{1},e_{2}) + \epsilon (f_{1},f_{2})$. 
We prove by induction on
$N$ that one of $B_1,B_2$ does not sit on $C$. We assume they both sit on $C$, for a contradiction.
\\
\\
(2) {\em Either $e_1,e_2$ are different and not twinned, or $f_1,f_2$ are different and not twinned.}
\\
\\
\Subproof Suppose that $e_{1}$ and $e_{2}$ are equal or twinned, and so are
$f_{1},f_{2}$. We claim that 
$$|\{e_1,e_2,f_1,f_2\}|\le 2,$$ 
and if this set has two members then they are twinned.
For suppose that $e_1 = e_2$. Since $P,Q$ cross, it follows that one of $f_1,f_2$ equals $e_1$,
say $f_1 = e_1$; and since either $f_2 = f_1$ or $f_2$ is twinned with $f_1$, the claim follows. So we may assume that
$e_1,e_2$ are twinned, and similarly so are $f_1,f_2$. But by (F5) and (F6), 
only one pair of edges of $C$ is twinned, and so again the claim holds.

Since $B_1$ sits on $C$, by (1) it has an $\eta$-attachment
$g\ne e_1$ that is not twinned with $e_1$; and so $g\ne e_1,e_2,f_1,f_2$.
Take a minimal path $R$ in $B_1$ between
$V(P \cup Q)$ and $V(\eta (g))$, and let its end $r$ in $P\cup Q$ be a vertex of $S$, say, where
$\{S,T\} = \{P,Q\}$. 
Let $S'$ be a path consisting of the union
of $R$ and a subpath of $S$ from $r$ to an
appropriate end of $S$, chosen such that $S',T$ cross. This contradicts the inductive 
hypothesis on $N$, and so proves (2).
\\
\\
(3) {\em $e_1\ne e_2$ and $f_1\ne f_2$.}
\\
\\
\Subproof Suppose that $e_{1}=e_{2}$, say.
Since $P,Q$ cross, one of
$f_{1},f_{2}$ equals $e_{1}$, say $f_{1}=e_{1}= e_{2}$; and by (2),
$f_{2} \neq f_{1}$, and $f_{1},f_{2}$ are not twinned.
By (1), $B_1$ has an
$\eta$-attachment $g\in E(C)$ not twinned with $e_1$.
Hence there is a minimal path $R$ of $B_1$ from $V(P)$ to
$V(Q) \cup \eta (g)$.
If it meets $\eta (g)$, we contradict the inductive hypothesis as before, 
so we assume $R$ has one end in $V(P)$ and the other in $V(Q)$.

Let $f_1 = uv$, and let $G' = G+(f_{1},f_{2})$ with new vertices $x,y$.
By adding $Q$ to $\eta (G)$ we see that there is a \he
$\eta ''$ of $G'$ in $H$ extending $\eta_F$ such that $ux,vx$ and $xy$
are all $\eta''$-attachments of some $\eta''$-bridge
(including $P \cup R$).
From \ref{gateleg}, we may choose $\eta''$ extending $\eta_F$ such that
$ux,vx,xy$ and some fourth edge $g$ are all $\eta''$-attachments
of some $\eta''$-bridge.
In other words, we may choose a \he $\eta'$ of $G$ in $H$
extending $\eta_F$ such that there exist 
\begin{itemize}
\item an $\eta'$-path $P'$ with ends
$p_{1}',p_{2}'$ in $V(\eta'(f_{1}))$, 
\item an $\eta'$-path $Q'$ with ends
$q_{1}',q_{2}'$ disjoint from $P'$, where $q_{1}'$ lies in
$\eta' (f_{1})$ between $p_{1}'$ and $p_{2}'$, and
$q_{2}' \in V(\eta' (f_{2}))$, 
\item a path $R'$
with one end in $P'$, the other end in $Q'$, and with no other
vertex or edge in $\eta' (G) \cup P' \cup Q'$, and
\item a path $S'$ with one end in $P' \cup R'$, the other end in
$\eta' (g)$ where $g \neq f_{1}$, and with no other vertex
or edge in $\eta' (G) \cup P' \cup Q' \cup R'$.
\end{itemize}

Let $B'$ be the $\eta'$-bridge containing
$P' \cup Q' \cup R' \cup S'$.
By \ref{trinity6}, there exists $C' \in {\cal C}$ such that all $\eta'$-attachments
of $B'$ are in $E(C')$.
Now $f_{1} \neq f_{2}$ and they are not
twinned, so $C'=C$, and hence $B'$ sits on $C$.
Let $T$ be an $\eta'$-path in $P'\cup R'\cup S'$ 
with one end in $\eta'(f_1)$ and the other in $\eta'(g)$, chosen such that 
$Q',T$ cross with respect to $\eta'(C)$. Then both $Q',T$ are contained in $B'$, and yet
$B'$ sits on $C$, and $\epsilon(f_1,g)<\epsilon(f_1,f_1)$, contrary to the inductive hypothesis.
This proves (3).
\\
\\
(4) {\em $e_1,e_2$ are not twinned, and $f_1,f_2$ are not twinned.}
\\
\\
\Subproof Suppose that $f_{1},f_{2}$ are twinned, say.
Let  $f_{1} = v_{1}x_{1}$ and $f_{2}=v_{2}x_{2}$ where either $C$ is a circuit and 
$v_{1}=v_{2} \in V(F)$, or $C$ is
a path with ends $v_{1},v_{2}$.
By (1), there is an $\eta$-attachment of $B_2$ different from $f_1,f_2$; and so there is a minimal $\eta$-path $R$ in $B_2$
from $V(Q)$ to $V(P) \cup V(\eta(C \setminus \{f_{1},f_{2}\}))$.
From the inductive hypothesis, $R$ does not meet
$\eta (C \setminus \{f_{1},f_{2}\})$,
and so it meets $P$.
Let $R$ have ends $r_{1} \in V(P)$ and $r_{2} \in V(Q)$, and for $i = 1,2$, let
$P_{i} = P[p_{i},r_{1}]$ and 
$Q_{i} = Q[q_{i},r_{2}]$.

Now for $i = 1,2$, $x_i\notin V(F)$ by (F5) (since if $C$ is a circuit then $v_1\in V(F)$ by (F6)).
For $i = 1,2$, let $g_i$ be the edge of $G$ not in $C_i$ incident with $x_i$, and let
$h_i$ be the edge of $C$ different from $f_i$ that is incident with $x_i$.

Now since either $C$ is a path and $f_{1},f_{2}$ are end-edges
of $C$, or $C$ is a circuit and $f_{1},f_{2}$ have a common end, and
since $P,Q$ cross, we may assume that 
$e_{1}=f_{1}$, and $p_1$ lies in $\eta(f_1)$ between $q_1$ and $\eta(v_1)$.
It follows that $e_{2} \neq f_{1},f_{2}$ by (2).

Suppose first that either $e_{2}=h_{1}$ or $x_{1}$ is adjacent to an end of
$e_{2}$.
By rerouting $h_{1}$ along $P$, we obtain a \he $\eta'$
of $G$ in $H$ extending $\eta_F$, such that $g_{1},h_{1}$ and
$f_{2}$ are all $\eta'$-attachments of some $\eta'$-bridge
(containing $Q\cup R$).
Since no member of ${\cal C}$ contains all of $g_{1},h_{1}$
and $f_{2}$, this contradicts \ref{trinity6}.
Hence $e_{2} \neq h_{1}$ and $x_{1}$ is not adjacent to any end of $e_{2}$.

By (F6), $|V(C)| \leq 6$, so either $e_{2}=h_{2}$, or $C$ is a circuit and $x_{2}$
is adjacent to an end of $e_{2}$. Suppose first that $C$ is a path; so $e_2 = h_2$.
By rerouting $h_{2}$ along $P_{2} \cup R \cup Q_{2}$ and adding
$P_{1}$ and $Q_{1}$, we obtain a \he (in $H$, extending $\eta_F$) of a cross extension of $G$ 
over $C$ of the third kind, contrary to (E6).
Thus $C$ is a circuit, and so $v_1\in V(F)$, and therefore $\{f_1,g_2,h_2\}$ is a circuit-type $Y$-trinity.
But then by rerouting $h_{2}$ along $P_{2} \cup R \cup Q_{2}$ and adding
$P_{1}$ and $Q_{1}$  we obtain a \he (in $H$, extending $\eta_F$) of an expansion of this 
$Y$-trinity of the first or second type, contrary to \ref{trinity2}.
This proves (4).
\\
\\
(5) {\em $e_1,e_2$ have no common end, and $f_1, f_2$ have no common end.}
\\
\\
\Subproof Suppose that $e_{1},e_{2}$ have a common end, $v$ say.
Since $e_1,e_2$ are not twinned by (4), it follows from \ref{frame1} that $v$ has degree three in $G\setminus E(F)$;
and so by (F5), $v\notin V(F)$. Since $P,Q$ cross, 
we may assume that $f_{1} = e_1$  and $p_1,q_1,\eta(v)$ are in order in $\eta(e_1)$.
Let $f,e_{1},e_{2}$ be the three edges of $G$ incident with $v$.
Let $\eta'$ be obtained from $\eta$ by
rerouting $e_{1}$ along $P$.
Then $\eta'$ is a
\he of $G$ in $H$ extending $\eta_F$, and $f,f_{1}$
and $f_{2}$ are $\eta'$-attachments in $E(G)\setminus E(F)$ of the $\eta'$-bridge
containing $Q$.
From \ref{trinity6}, there exists $C' \in {\cal C}$ with
$f,f_{1},f_{2} \in E(C')$.
Since $f \not\in E(C)$ it follows that $C' \neq C$, and so
$f_{1},f_{2}$ are twinned, contrary to (4). This proves (5).

\bigskip

Thus $e_{1},e_{2}$ have no common  end, and nor do $f_{1},f_{2}$.
By (E6), we may assume that one end of $e_{1}$ is adjacent to one
end of $e_{2}$. Since $P,Q$ cross, 
we may therefore assume that for some edge
$g=uv$ of $C$, $u$ is an end of $e_{1}$, $v$ is an end of
$e_{2}$, $f_{1} \in \{e_{1},g\}$, and if $f_1 = e_1$ then $p_1,q_1,\eta(u)$ are in order in $\eta(e_1)$.
Let $u$ be incident with $g,e_{1},g_{1}$ and $v$ with
$g,e_{2},g_{2}$.

Suppose that $u\notin V(F)$. Let $\eta'$ be the \he obtained from $\eta$ by rerouting $g$ along $P$. By 
By (E2), applied to $\eta'(G)$ with edges $f_2,g_1$, it follows
that there exists $C_1\in \mathcal{C}$ with $f_2,g_1\in E(C_1)$. By (F3), $C_1$ contains one of $e_1,g$,
say $h$. Hence $h,f_2$ are twinned; and since $f_1,f_2$ are not twinned, it follows that $\{e_1,g\} = \{f_1,h\}$.
If $h = g$, then $f_1 = e_1$; and since $g$ is not an
end-edge of $C$, (F6) implies that $f_2 = e_2$ and $g_2\in E(F)$.
But then $\{e_1,g_1,e_2\}$ is a $Y$-type trinity, and adding $P$ and $Q$ provides
a \he (in $H$, extending $\eta_F$) of an expansion of this trinity, contrary
to \ref{trinity2}. Thus $h\ne g$, and so $h = e_1$ and $f_1 = g$,
and $f_2\ne e_1,g,e_2$.
If $C$ is a circuit, (F6) implies that the
end of $e_1$ different from $u$ belongs to $V(F)$; but then by (F5), $v\notin V(F)$, and the symmetry between 
$u$ and $v$ implies that the end of $e_2$ different from $v$ belongs to $V(F)$, contrary to (F5). If $C$
is a path, then (F6) implies that $e_1$ is an end-edge of $C$, and $v\notin V(F)$;
but then the symmetry between $u,v$ implies that $e_2$ is also an end-edge of $C$, a contradiction.

This proves that $u\in V(F)$. Consequently $v\notin V(F)$, and it follows 
(by exchanging $P,Q$, and exchanging $e_1,e_2$) that $f_2\ne e_2$.
Since $C$ contains $e_1,g$, it follows that $u$ is not an end of $C$,
and so by (F5), $g_1\in E(F)$. 
By (F2) there exists $C_2\in \cal C$ containing $g,g_2$, since $v\notin V(F)$. 
Since $g_1\in E(F)$, we deduce that $e_1\in E(C_2)$. Since $f_1,f_2$
are not twinned, it follows that $f_2\notin E(C_2)$. Thus $g_2\in E(C_2)\setminus E(C)$, 
and $f_2\in E(C)\setminus E(C_2)$, and $f_2,g_2$ have no common end,
since $f_2\ne e_2$.
But rerouting $g$ along $P$ gives a \he of $G$ in $H$ extending $\eta_F$, and adding $\eta(g)$
and $Q$ to it contradicts (E4).
This proves \ref{extendbranch}.~\bbox

\section{The bridges between twins}

To apply these results about frameworks, we have to choose a \he
$\eta$ of $G$ in $H$, and there is some freedom in how we do so.
If we choose it carefully we can make several problems disappear
simultaneously.
The most important consideration is to ensure that each
$\eta$-bridge has at least two $\eta$-attachments, but that
is rather easy.
With more care, we can also discourage $\eta$-bridges from having
$\eta$-attachments in certain difficult places.
To do so, we proceed as follows.

Let $(G,F, \cal{C})$ be a framework, and let 
$\eta$ be a \he of $G$ in $H$ extending $\eta_F$, as usual.
An edge $e$ of $G$ is a {\em twin} if there exists $f$ such that $e,f$ are twinned. 
(Thus, stating that ``$e,f$ are twins'' does not imply that
they are twinned with each other.)
An edge $e\in E(G)\setminus E(F)$ is 
\begin{itemize}
\item {\em central} if it does not belong to any path in ${\cal C}$ and is not a twin; 
\item {\em peripheral} if $e$ is an internal edge of some path in ${\cal C}$ 
\item {\em critical} if either $e$ is a twin or $e$ is an end-edge of some path in ${\cal C}$.
\end{itemize}
By (F4) and (F6), no edge is both peripheral and critical, so every edge of $E(G)\setminus E(F)$ is of exactly one of these three kinds.

An edge $f\in E(H)$ is said to {\em $\eta$-attach}
to $e\in E(G)$ if there is a path $P$ of $H$ with no internal vertex in $V(\eta(G))$ with $f\in E(P)$ and with one end a
vertex of $\eta(e)$. (Thus $f$ $\eta$-attaches to $e$ if and only if either $f\in E(\eta(e))$ or
$f$ belongs to an $\eta$-bridge for which $e$ is an $\eta$-attachment.)
Let 
\begin{itemize}
\item $L_1(\eta)$ be the set of edges in $E(H)$ that $\eta$-attach to some central edge of $G$;
\item $L_2(\eta)$ be the set of edges in $E(H)$ that $\eta$-attach to an edge of $G$ which is either peripheral or central;
\item $L_3(\eta)$ be the set of edges in $E(H)$ that attach to two edges of $G$ that are not twinned; and
\item $L_4(\eta)$ be the set of edges in $E(H)$ that attach to two edges of $G$.
\end{itemize}
We say that $\eta$ is {\em optimal} if it is chosen (among all homeometric embeddings of $G$ in $H$ extending $\eta_F$)
with the four-tuple of cardinalities of these sets lexicographically maximum; that is, for every \he $\eta'$ extending $\eta_F$, 
there exists $j\in \{1\l 5\}$ such that $|L_i(\eta)| = |L_i(\eta')|$ 
for $1\le i<j$, and $|L_j(\eta)| > |L_j(\eta')|$ if $j\le 4$.
In this section we study the properties of optimal embeddings.

\begin{thm}\label{twoatt}
Let $\eta$ be an optimal \he of $G$ in $H$ extending $\eta_F$.
Then every $\eta$-bridge has at least two $\eta$-attachments.
\end{thm}
\Proof
Let $e \in E(G)  \setminus  E(F)$. Let us say an $\eta$-bridge is {\em singular} if $e$ is its only $\eta$-attachment, and {\em nonsingular}
otherwise. Suppose that there is a singular
$\eta$-bridge.
Let $e = uv$, let $p_1\l p_r$ be the set of vertices of $\eta(e)$ that belong to
nonsingular $\eta$-bridges, and let $p_0 = \eta(u)$ and $p_{r+1} = \eta(v)$, numbered such that
$p_0,p_1\l p_{r+1}$ are in order in $\eta(e)$. For $0\le i\le r$ let $P_i =\eta(e) [p_i, p_{i+1}]$.
Choose $j$ with $0\le j\le t$ such that some singular $\eta$-bridge
contains a vertex of $P_j$. Since $H$ is three-connected, there is an $\eta$-bridge $B$ containing a vertex $b$ of the interior of $P_j$ and containing
a vertex $a$ of $\eta(G)$ not in $P_j$. From the definition of $p_1\l p_r$, it follows that $B$ is singular.
Hence there exists
$i\ne j$ with $0\le i\le r$ such that $a$ belongs to $P_{i}$, and from the symmetry we may assume that $i<j$.
Let $P$ be an $\eta$-path
in $B$ between $a,b$. Let $\eta'$ be obtained from $\eta$ by rerouting $e$ along $P$.
For every edge $f$ of $E(H)$, every $\eta$-attachment of $f$ is also an $\eta'$-attachment.
Consequently $L_i(\eta)\subseteq  L_i(\eta')$ for $1\le i\le 4$.
But the edge of $P_j$ incident with $p_j$ belongs to $L_4(\eta')\setminus L_4(\eta)$, 
contrary to the optimality of $\eta$. This proves \ref{twoatt}.~\bbox

\begin{thm}\label{peripheral}
Let $\eta$ be an optimal \he of $G$ in $H$ extending $\eta_F$.
Let $C\in {\cal C}$ be a path, and suppose that $B$ is an $\eta$-bridge and all its $\eta$-attachments are edges of $C$.
Then its $\eta$-attachments are pairwise diverse in $C$.
\end{thm}
\Proof We claim first
\\
\\
(1) {\em If $e,f$ are edges of $C$ with a common end $v$, and $g$ is the third edge of $G$ incident with $v$, 
then $v\notin V(F)$, and either 
$g$ is central, or $g$ is peripheral and one of $e,f$ is an end-edge of $C$.}
\\
\\
\Subproof Certainly $v\notin V(F)$ by (F5), since $C$ is a path.
If $g$ does not belong to any path of ${\cal C}$ then it is not a twin by (F6), and so it is central. Thus we may assume
that there is a path $C'\in {\cal C}$ containing $g$. By (F4), $C'$ contains 
one of $e,f$, say $e$, and $e$ is an end-edge of both $C,C'$. 
Now (F1) implies that $g$ is not an end-edge of $C'$,
and so by (F6), $g$ is not a twin, and by (F4) $g$ is not an end-edge of any path in $\mathcal{C}$, that is, 
$g$ is peripheral. This proves (1).
\\
\\
(2) {\em No two $\eta$-attachments of $B$ in $C$ have a common end.}
\\
\\
\Subproof Suppose that $e,f$ are $\eta$-attachments of $B$, and they have a common end $v$. Let $g$ be the third edge of $G$ incident with $v$.
Choose a path $P$ in $B$ from a vertex $a$ of $\eta(e)$
to a vertex $b$ of $\eta(f)$. 
Let $\eta'$ be obtained from $\eta$ by rerouting $f$ along $P$.
Then $\eta'$ is a \he of $G$ in $H$ extending $\eta_F$ (note that $g\notin E(F)$ since $v\notin V(F)$
by (1)). Moreover, since no $\eta$-attachment of $B$ is central, it follows that $L_1(\eta)\subseteq L_1(\eta')$, and therefore
equality holds. In particular, the edge of $\eta(e)$ incident with $\eta(v)$ therefore does not belong to $L_1(\eta')$, and so $g$ is not central.
We deduce from (1) that $g$ is peripheral and one of $e,f$ is an end-edge of $C$, 
and from the symmetry we may assume that $e$ is an end-edge of $C$. 
Thus $f$ is peripheral, and it follows that  $L_2(\eta)\subseteq L_2(\eta')$, and therefore
equality holds. But the edge of $\eta(e)$ incident with $\eta(v)$ belongs to $L_2(\eta')$, and does not belong to
$L_2(\eta)$ since $e$ is an end-edge of $C$, a contradiction. This proves (2).

\bigskip
To complete the proof, suppose that some two $\eta$-attachments $e,f$ of $B$ in $C$ are not diverse in $C$.
Then by (2), there are consecutive vertices $u,v,w,x$ of $C$, such that $e = uv$ and $f = wx$.
Let the third edge of $G$ at $v$ be $g$ and at $w$ be $h$. Choose a path $P$ in $B$ from a vertex $a$ of $\eta(e)$
to a vertex $b$ of $\eta(f)$. Let $\eta'$ be obtained from $\eta$ by rerouting $vw$ along $P$.
Then $\eta'$ 
is a \he of $G$ in $H$ extending $\eta_F$. Since no $\eta$-attachment of $B$ is central, it follows that $L_1(\eta)\subseteq L_1(\eta')$, and therefore
equality holds. In particular,  the edge of $\eta(e)$ incident with $\eta(v)$ does not belong to $L_1(\eta')$, and so $g$ is not central. 
From (1), it follows that $g$ is peripheral and $e$ is an end-edge of $C$. Similarly $h$ is peripheral and $f$ is an end-edge of $C$. Hence
$L_2(\eta)\subseteq L_2(\eta')$, and therefore equality holds. But the edge of $\eta(e)$ incident with $\eta(v)$ 
belongs to $L_2(\eta')$ and not to $L_2(\eta)$ since $e$ is an end-edge of $C$, a contradiction. 
This proves \ref{peripheral}.~\bbox

\bigskip

If $C \in {\cal C}$, we denote by ${\cal A}(C)$ the set of all $\eta$-bridges
that sit on $C$.
If $e,f$ are twinned edges of $G$, 
we denote by ${\cal A}(e,f)$ the set of all $\eta$-bridges
that have no attachments different from $e,f$. 
Thus, if $\eta$ is optimal, then by \ref{twoatt} every bridge belongs to ${\cal A}(C)$ 
for some $C$ or to ${\cal A}(e,f)$ for some $e,f$, and to only one such set 
(except that ${\cal A}(e,f) = {\cal A}(f,e)$).
The next four theorems are all about a pair of twinned edges $e,f$, and it is convenient first to set up some notation.
Thus, let $e,f$ be twinned edges of $G$. Let there be $r$ vertices $p_{1},\ldots,p_{r}$ of $\eta (e)$
that belong to an $\eta$-bridge with an $\eta$-attachment different
from $e$ and $f$, and let $\eta (e)$ have ends $p_{0}$ and
$p_{r+1}$, numbered such that $p_{0},\ldots,p_{r+1}$ are in order
in $\eta (e)$.
For $0 \leq i \leq r$, let $P_{i}$
$= \eta (e) [p_{i},p_{i+1}]$.
Let $q_{0},\ldots,q_{s+1} \in V(\eta (f))$ and
$Q_{0},\ldots,Q_{s}$ be defined similarly.

\begin{thm}\label{twinbridge1}
Let $\eta$ be an optimal \he of $G$ in $H$ extending $\eta_F$, and let $e,f$ be twinned edges of $G$. 
With notation as above, for every $B \in {\cal A}(e,f)$ there exist $i$ and $j$ with
$0 \leq i \leq r$ and $0 \leq j \leq s$ such that
$B \cap \eta (e) \subseteq P_{i}$ and $B \cap \eta (f) \subseteq Q_{j}$.
\end{thm}
\Proof
Suppose that some member $B$ of ${\cal A}(e,f)$ meets both $P_i$ and $P_j$,
where $0 \leq i < j \leq r$.
Let $P$ be an $\eta$-path
in $B$ between some $a\in V(P_i)$ and some $b\in V(P_j)$.
Let $\eta'$ be obtained from $\eta$ by rerouting $e$ along $P$.
Since no $\eta$-attachment of $B$ is central or peripheral, and no edge of $B$ is in $L_3(\eta)$,
it follows that $L_i(\eta)\subseteq L_i(\eta')$ for $i = 1,2,3$, and so equality holds in all three.
Let $B'$ be an $\eta$-bridge
containing $p_i$; then $B'$ has an $\eta$-attachment different from $e,f$, say $g$. Consequently $e,g$ are not twinned, and in particular,
the edge of $P_j$ incident with $p_j$ is in $ L_3(\eta')$, a contradiction. This proves \ref{twinbridge1}.~\bbox

\begin{thm}\label{twinbridge2}
Let $\eta$ be an optimal \he of $G$ in $H$ extending $\eta_F$, and let $e,f$ be twinned edges of $G$.
Suppose that $e,f$ have a common end $v$, and let $e = uv$ and $f = vw$. Then ${\cal A}(e,f)$ can be numbered as $\{B_1\l B_k\}$, 
such that 
\begin{itemize}
\item $B_i$ has only one edge $c_id_i$ for $1\le i\le k$; 
\item $\eta(u),c_1\l c_k, \eta(v)$ are in order in $\eta(e)$, and $\eta(w),d_1\l d_k,\eta(v)$ are in order in $\eta(f)$; and
\item for $1\le i<k$, one of $\eta(e)[c_i, c_{i+1}]$, $\eta(f)[d_i,d_{i+1}]$ contains a vertex of some $\eta$-bridge not in  ${\cal A} (e,f)$.
\end{itemize}
\end{thm}
\Proof
Using the notation established earlier, we may assume that $\eta(v) = p_0 = q_0$.
\\
\\
(1) {\em Suppose that $M,N$ are disjoint $\eta$-paths, from $m$ to $m'$ and from $n$ to $n'$ respectively, such that
\begin{itemize}
\item $\eta(u),m,n,\eta(v),m',n',\eta(w)$ are in order in the path $\eta(e)\cup \eta(f)$; and
\item no edge of $M\cup N$ belongs to $L_2(\eta)$.
\end{itemize}
Then there exist $i,j$ with $0 \leq i \leq r$ and $0 \leq j \leq s$ 
such that $m,n$ belong to $P_i$ and $m',n'$ belong to $P_j$.}
\\
\\
\Subproof Let $m$ be in $P_i$ and $n$ be in $P_{h}$ where $0\le h<i\le r$.
Let
$$\eta'(e) = \eta(e)[\eta(u),m]\cup M\cup \eta(f)[m',\eta(v)].$$
 and
$$\eta'(f) = \eta(e)[\eta(v),n]\cup N\cup \eta(f)[n',\eta(w)].$$
Then $\eta'$ is a \he of $G$ in $H$ extending $\eta_F$. Since no edge of 
the $\eta$-bridges containing $M$ or $N$ belongs to $L_1(\eta)$ or to $L_2(\eta)$, and $e,f$
are critical, it follows that $L_i(\eta)\subseteq L_i(\eta')$ for $i = 1,2$, and so equality holds in both. Let $B$ be the $\eta$-bridge containing
$p_i$. Then there is an $\eta$-attachment $g\ne e,f$ of $B$. Choose $C\in {\cal C}$ containing $e,g$ (this is possible by (E2) 
applied to $\eta(G)$ with edges $e,g$). From
(F6), $C$ is a circuit, and so $g$ is not critical from (F5). Hence $g$ is either central or peripheral, and so the edges of $\eta(e)$
incident with $p_i$ belongs to $L_2(\eta')$, a contradiction. This proves (1).

\bigskip
To complete the proof,
for $0 \leq i \leq r$ and $0 \leq j \leq s$ let ${\cal A}_{ij}$
be the set of all $B \in {\cal A} (e,f)$ with $B \cap \eta(e)\subseteq P_{i}$ and
$B \cap \eta(f)\subseteq Q_{j}$.
From (1), ${\cal A}(e,f) = \bigcup {\cal A}_{ij}$.
For each $i,j$ let $J_{ij}$ be the union of all members of ${\cal A}_{ij}$.
Suppose that some $|E(J_{ij}) | \geq 2$.
Since $H$ is cyclically five-connected by (E1), we may assume (by
exchanging $e$ and $f$ if necessary) that there are
$b_{1},b',b_{2}$ in $P_{i}$, in order, such that $b_{1}$ and $b_{2}$
both belong to $J_{ij}$, and $b'$ belongs to some
$\eta$-bridge $B' \not\in {\cal A}_{ij}$.
Since $b' \neq p_{1},\ldots,p_{r}$ it follows that
$B' \in {\cal A} (e,f)$, and so $B' \in {\cal A}_{ij'}$,
for some $j'\ne j$. In particular, $J_{ij}$ and $J_{ij'}$ are disjoint. By \ref{twoatt} it
follows that there is a path $M$ in $J_{ij}$ and a path $N$ in $J_{ij'}$ violating (1) (possibly with $M,N$ exchanged).
This proves that each $J_{ij}$ has at most one edge, and in particular from \ref{twinbridge1}, each $\eta$-bridge in ${\cal A} (e,f)$ has only one edge.
The result follows from (2) applied to the paths of length one formed by these $\eta$-bridges. This proves \ref{twinbridge2}.~\bbox

\begin{thm}\label{twinbridge3}
Let $\eta$ be an optimal \he of $G$ in $H$ extending $\eta_F$, and let $e,f$ be twinned edges of $G$. 
Suppose that $e,f$ are disjoint, and there is no path $C\in {\cal C}$ of length five with end-edges $e,f$.
Then 
\begin{itemize}
\item there is at most one $\eta$-bridge in ${\cal A}(e,f)$, and any such $\eta$-bridge has only one edge; 
\item no other $\eta$-bridge
contains any vertex of $\eta(e)\cup \eta(f)$; and
\item ${\cal A}(C) = \emptyset$ for every member of ${\cal C}$ containing $e$ or $f$.
\end{itemize}
\end{thm}
\Proof
Now there is a path in ${\cal C}$ with end-edges $e,f$, and so every member $C$ of ${\cal C}$ containing $e$ or $f$ is a path, by (F4).
Moreover, if $e,f\in E(C)$ then $C$ has length at most four by hypothesis and (F6), and $C$ has end-edges $e,f$, and therefore every member 
of ${\cal A}(C)$ has an $\eta$-attachment some edge of $C$ different from $e,f$. By \ref{peripheral}, this implies that ${\cal A}(C) = \emptyset$.
On the other hand, if $C\in {\cal C}$ contains just one of $e,f$ then $C$ has length three by (F6), and again ${\cal A}(C) = \emptyset$ 
by \ref{peripheral}. This proves the third assertion. Consequently, $r = s = 0$ (in our previous notation). Since $H$ is cyclically five-connected
by (E1), it follows that
the union of all $\eta$-bridges in  ${\cal A} (e,f)$ and the paths $\eta(e), \eta(f)$ contains no circuit; and so there is at most one
$\eta$-bridge in  ${\cal A} (e,f)$ and any such $\eta$-bridge has only one edge. This proves \ref{twinbridge3}.~\bbox

\begin{thm}\label{twinbridge4}
Let $\eta$ be an optimal \he of $G$ in $H$ extending $\eta_F$, and let $e,f$ be twinned edges of $G$.
Suppose that $e,f$ are disjoint, and there exists $C\in {\cal C}$ of length five with end-edges $e,f$.
Then:
\begin{itemize}
\item ${\cal A}(C')$ is empty for every $C'\ne C$ in ${\cal C}$ containing $e$ or $f$;
\item the vertices of $C$ can be numbered in order as $v_0\d v_1\c v_5$, such that for each $B\in {\cal A}(C)$, 
its only $\eta$-attachments are $v_1v_2$ and $v_4v_5$ (and we may assume that $e = v_0v_1$ and $f = v_4v_5$,
possibly after exchanging $e,f$);
\item ${\cal A}(e,f)$ can be numbered as $\{B_1\l B_k\}$ such that $B_i$ has exactly one edge $c_id_i$ for $1\le i\le k$, where $c_i\in V(\eta(e))$ and $d_i\in V(\eta(f))$; and
\item $\eta(v_0),c_1\l c_k, \eta(v_1)$ are in order in $\eta(e)$, and $\eta(v_4),d_1\l d_k,\eta(v_5)$ are in order in $\eta(f)$.
\end{itemize}
\end{thm}
\Proof Let $C\in {\cal C}$ of length five with end-edges $e,f$. 
\\
\\
(1) {\em The first assertion of the theorem is true.}
\\
\\
\Subproof By (F7), every other path in ${\cal C}$ containing $e$ or $f$
has length at most four. If $C' \in {\cal C}$ contains both $e,f$, then ${\cal A}(C') = \emptyset$ by \ref{peripheral}, since each
member of ${\cal A}(C')$ has an $\eta$-attachment in $C$ different from $e,f$; and if $C' \in {\cal C}$ contains just one of $e,f$,
then it has length three by (F6), and again ${\cal A}(C') = \emptyset$ by \ref{peripheral}. This proves (1).
\\
\\
(2) {\em The second assertion is true.}
\\
\\
\Subproof
Let $C$ have vertices $v_0\d v_1\c v_5$ in order, where $e = v_0v_1$ and $f = v_4v_5$. Let $B\in {\cal A}(C)$. By \ref{peripheral},
one of $e,f$ is an $\eta$-attachment of $B$, say $f$; and since $B$ has two $\eta$-attachments in $C$ and they are diverse in $C$ by \ref{peripheral},
and $e,f$ are twinned, it follows that the only other $\eta$-attachment of $B$ is $v_1v_2$.
Let $B' \in {\cal A} (C)$ with $B' \neq B$; we claim that
$v_{1}v_{2}$ and $v_{4}v_{5}$ are the $\eta$-attachments of $B'$.
For if not, then by the previous argument $v_{0}v_{1}$ and $v_{3}v_{4}$ are
$\eta$-attachments of $B'$, contrary to (E6). This proves (2).

\bigskip

In our earlier notation, we may assume that $p_0 = \eta(v_0)$ and $q_0 = \eta(v_4)$. 
Suppose that $B$ is an $\eta$-bridge not in ${\cal A}(e,f)$ that meets $\eta(e)$. Then from \ref{twoatt} and \ref{trinity6},
$B\in {\cal A}(C')$ for some $C'\in {\cal C}$ containing $e$, and hence $B\in {\cal A}(C)$ from (1); but this contradicts (2).
Consequently $r = 0$.
\\
\\
(3) {\em Suppose that $M,N$ are disjoint $\eta$-paths, from $m$ to $m'$ and from $n$ to $n'$ respectively, where
$\eta(v_0),m,n,\eta(v_1)$ are in order in $\eta(e)$, and $\eta(v_4) ,n',m',\eta(v_5)$ are in order in $\eta(f)$.
Then there exists $j$ with $0 \leq j \leq s$ such that $m',n'$ belong to $Q_j$.}
\\
\\
\Subproof
Suppose not; then there exist distinct $j,j'$ with $m'\in V(Q_j)$ and $n'\in V(Q_{j'})$, and consequently $j<j'$. 
Let $B$ be the $\eta$-bridge containing $q_{j'}$; then
$B\notin {\cal A}(e,f)$ from the definition of $q_1\l q_s$, and so $B$ has an $\eta$-attachment $g\ne e,f$. From \ref{trinity6},
and (1) it follows that $B\in {\cal A}(C)$, and $g = v_1v_2$. In particular, $B$ is disjoint from $M,N$. Choose an $\eta$-path $P$
in $B$ from $q_{j'}$ to $V(\eta(v_1v_2))$; then $M,N,P$ contradict (E7). This proves (3).

\bigskip

For $0 \leq j \leq s$ let ${\cal A}_{j}$
be the set of all $B \in {\cal A} (e,f)$ with 
$B \cap \eta(f)\subseteq Q_{j}$.
From (1), ${\cal A}(e,f) = \bigcup {\cal A}_{j}$.
For each $j$ let $J_{j}$ be the union of all members of ${\cal A}_{j}$.
Suppose that some $|E(J_{j}) | \geq 2$.
Since $H$ is cyclically five-connected by (E1), 
there are distinct
$b_{1},b',b_{2}$ in $\eta(e)$, in order, such that $b_{1}$ and $b_{2}$
both belong to $J_{j}$, and $b'$ belongs to some
$\eta$-bridge $B' \not\in {\cal A}_{j}$.
Since $b' \neq p_{1},\ldots,p_{r}$ it follows that
$B' \in {\cal A} (e,f)$, and so $B' \in {\cal A}_{j'}$,
for some $j'\ne j$. In particular, $J_{j}$ and $J_{j'}$ are disjoint. By \ref{twoatt} it
follows that there is a path $M$ in $J_{j}$ and a path $N$ in $J_{j'}$ violating (1) (possibly with $M,N$ exchanged).
This proves that $|E(J_{j}) | \le 1$ for $0\le j\le s$. Thus every $\eta$-bridge in $\mathcal{C}(e,f)$ has only one edge, and
no two of them have ends in the same $Q_j$.
The result follows from (3) applied to the paths of length one formed by these $\eta$-bridges. This proves \ref{twinbridge4}.~\bbox


\section{Flattenable graphs}

Let $(G,F, {\cal C})$ be a framework and let $H, \eta_{F}$
satisfy (E1).
We say that $H$ is {\em flattenable onto $(G,F, {\cal C})$ via $\eta_{F}$} if there is
\begin{itemize}
\item
a \he $\eta$ of $G$ in $H$ extending $\eta_F$
\item
a set of $\eta$-bridges ${\cal B}(C)$, for each $C \in {\cal C}$, and
\item
an edge $N(e)$ of $\eta (e)$, for each edge $e$ of $G\setminus E(F)$ such that for some
edge $f \neq e$, $e$ and $f$ are twinned and have no common end
\end{itemize}
with the following properties.
For each $C \in {\cal C}$, if $C$ is a circuit let $P(C)$ be $\eta (C)$,
and if $C$ is a path let $P(C)$ be the maximal subpath of $\eta (C)$
that contains $\eta (g)$ for every $g \in E(C)$ that is not an
end-edge of $C$, and does not contain any edge $N(e)$.
Then we require:
\begin{itemize}
\item
every $\eta$-bridge belongs to exactly one set ${\cal B} (C)$
\item
if $B \in {\cal B}(C)$ then $B \cap \eta (G) \subseteq P(C)$
\item
for $C \in {\cal C}, P(C) \cup \bigcup (B:B\in {\cal B}(C))$ is
$P(C)$-planar.
\end{itemize}

The main result, that everything so far has been directed towards,
and of which all the other results in the paper will be a consequence,
is the following.

\begin{thm}\label{mainthm}
Let $(G,F, {\cal C})$ be a framework, and let
$H, \eta_{F}$ satisfy
$(E1)$--$(E7)$. Suppose that there is a \he of $G$ in $H$ extending $\eta_F$.
Then $H$ is flattenable onto $(G,F, {\cal C})$ via $\eta_{F}$.
\end{thm}
\Proof
Since there is a \he of $G$ in $H$ extending $\eta_F$, there is an optimal one, say $\eta$.
We will prove that $\eta$ provides the required flattening. We begin with
\\
\\
(1) {\em If $e,f \in E(G)$ are twinned and have a common end, there exists
$C \in {\cal C}$ containing $e,f$ such that
\[
\eta (C) \cup \bigcup (B:\;B\in {\cal A} (C) \cup {\cal A} (e,f))
\]
is $\eta (C)$-planar.}
\\
\\
\Subproof
Let the two members of ${\cal C}$ that contain $v$ be
$C_{1},C_{2}$, where $v$ is the common end of $e$ and $f$. Let $e = uv$ and $f = vw$, and let
$c_{1}d_{1},\ldots,c_{k}d_{k}$ be the edges of $H$ with
one end in $\eta (e)$ and the other in $\eta (f)$ (these are
the edges of the bridges in ${\cal A}(e,f)$) numbered as in \ref{twinbridge2}).
By \ref{extend} we may assume that $k \geq 1$.
Now
\[
\eta (C_{i}) \cup \bigcup (B:\;B \in {\cal A} (e,f))
\]
is $\eta (C_{i})$-planar for $i=1,2$.
We claim that
for either $i=1$ or $i=2$, no member of ${\cal A}(C_{i})$ meets
$\eta (e) \cup \eta (f)$ between $c_{1}$ and $d_{1}$.
For if not, there are disjoint $\eta$-paths $R_{1},R_{2}$ such that
for $i=1,2$, $R_{i}$ has one end $r_{i}$ in
$\eta (e) \cup \eta (f)$ between $c_{1}$ and $d_{1}$, and its other end
$s_{i}$ is in $\eta (C_{i})$ and not in $\eta (e) \cup \eta (f)$.
Let $s_{i} \in V(\eta (g_{i}))\;(i=1,2)$.
If $g_{1},g_{2}$ have no common end, this contradicts (E4),
and if they have a common end, this contradicts \ref{trinity2}.
(To see this, in each case delete an appropriate end-edge of the subpath of $\eta(e)\cup \eta(f)$ between $c_1,d_1$.)
We may therefore assume that no member of ${\cal A}(C_{1})$ meets
$\eta (e) \cup \eta (f)$ between $c_{1}$ and $d_{1}$.
But then by \ref{extend}, the claim holds.
This proves (1).

\bigskip
For edges $e,f$ as in (1), let $D(e,f)$ be some $C\in \mathcal{C}$ satisfying (1).
\\
\\
(2) {\em Let $e,f$ be twinned, with no common end.
Then there are edges $N(e)$ of $\eta(e)$ and $N(f)$ of $\eta (f)$,
and distinct paths $C_1,C_2\in \mathcal{C}$, both with end-edges $e,f$ and with the following property, 
where for $i=1,2$, $P(C_i)$
denotes the component of $\eta (C_i) \setminus \{N(e),N(f)\}$
containing $\eta(g)$ for each internal edge $g$ of $C$.
\begin{itemize}
\item $\mathcal{A}(C)=\emptyset$ for all $C\in \mathcal{C}$ containing either $e$ or $f$ and different from $C_1$
\item $B \cap \eta (G) \subseteq P(C_{1})$ for all $B \in \mathcal{A}(C_1)$
\item $B \cap \eta (G) \subseteq P(C_{2})$ for all $B \in {\cal A} (e,f)$, and
\[
P(C_{2}) \cup \bigcup (B:\;B \in {\cal A} (e,f))
\]
is $P(C_{2})$-planar.
\end{itemize}}
\Subproof
By \ref{frame2} there are at least two paths in ${\cal C}$ with end-edges $e,f$, and by (F6) every such path has length at most five.
If there is no path in $\mathcal{C}$ with end-edges $e,f$ and with length exactly five, 
the claim follows from \ref{twinbridge3}, 
so we assume that some such path has length five, say $C_1$.
By \ref{twinbridge4}, ${\cal A}(C)$ is empty for every $C\ne C_1$ in ${\cal C}$ containing $e$ or $f$, so the first assertion
of the claim holds.
Moreover, also by \ref{twinbridge4}, 
\begin{itemize}
\item the vertices of $C_1$ can be numbered in order as $v_0\d v_1\c v_5$, such that for each $B\in {\cal A}(C)$,
its only $\eta$-attachments are $v_1v_2$ and $v_4v_5$ (and we may assume that $e = v_0v_1$ and $f = v_4v_5$,
possibly after exchanging $e,f$);
\item ${\cal A}(e,f)$ can be numbered as $\{B_1\l B_k\}$ such that $B_i$ has exactly one edge $c_id_i$ for $1\le i\le k$, where $c_i\in V(\eta(e))$ and $d_i\in V(\eta(f))$; and
\item $\eta(v_0),c_1\l c_k, \eta(v_1)$ are in order in $\eta(e)$, and $\eta(v_4),d_1\l d_k,\eta(v_5)$ are in order in $\eta(f)$.
\end{itemize}

Let $N(e)$ be the
edge of $\eta (e)$ incident with $\eta (v_{1})$, and $N(f)$
be the edge of $\eta (f)$ incident with $\eta (v_{5})$.
Then $B \cap \eta (G) \subseteq P(C)$ for all
$B \in {\cal A} (C)$,
so the second assertion holds.

By (F7) there exists $C_{2} \in {\cal C}$ with end-edges $e$ and $f$
and with ends $v_{1}$ and $v_{5}$. It follows that $N(e)$ and $N(f)$ are the end-edges of $C_2$, and so
$B \cap \eta (G) \subseteq P(C_{2})$ for all
$B \in {\cal A} (e,f)$. From the second and third bullets above,
\[
P(C_{2} ) \cup \bigcup (B:B \in {\cal A} (e,f))
\]
is $P(C_{2})$-planar. So the third assertion holds.
This proves (2).

\bigskip

For $e,f$ as in (2), choose $C_1,C_2$ as in (2), and define $D(e,f) = C_2$.
For each edge $e$ that is twinned with an edge disjoint from $e$, choose $N(e)$ as in (2).
Since no edge of $e$ is twinned with more than one other edge,
by \ref{uniquetwin}, this is well-defined. For each $C\in \mathcal{C}$, if $C$ is a circuit let $P(C) = C$,
and if $C$ is a path let $P(C)$
be the maximal subpath of $\eta (C)$
that contains $\eta (g)$ for every $g \in E(C)$ that is not an
end-edge of $C$, and does not contain any edge $N(e)$.
\\
\\
(3) {\em For every path $C\in \mathcal{C}$, $B \cap \eta (G) \subseteq P(C)$ for each $B \in {\cal A} (C)$.}
\\
\\
For let $C\in \mathcal{C}$ be a path. If $P(C)=C$ the claim is true, so we may assume that
some edge $e$ of $C$ is twinned with some other edge $f$ disjoint from $e$, and so $N(e)$
is defined. Choose $C_1,C_2$ satisfying (2), where $C_2 = D(e,f)$. If $C\ne C_1$ then 
${\cal A} (C)=\emptyset$ and the claim is trivial, by the first assertion of (2); while if $C=C_1$
then the claim holds by the second assertion of (2). This proves (3).

\bigskip
For each $C\in \mathcal{C}$, let $\mathcal{B}(C)$ be the following set of $\eta$-bridges:
\begin{itemize}
\item if $C = D(e,f)$ for some pair $e,f$ of twinned edges with a common end, let
$\mathcal{B}(C) = {\cal A} (C) \cup {\cal A} (e,f)$
\item if $C= D(e,f)$ for some pair $e,f$ of twinned edges with no common end, let
$\mathcal{B}(C) = {\cal A} (e,f)$
\item otherwise, let $\mathcal{B}(C)=\mathcal{A}(C)$.
\end{itemize}

Now let $B$ be an $\eta$-bridge. We claim that $B$ belongs to exactly one set ${\cal B} (C)$. For if $B$ sits on some
$C'\in \mathcal{C}$, then for $C\in \mathcal{C}$, $B\in \mathcal{C}$ if and only if $C=C'$; and otherwise,
$B$ belongs to $\mathcal{A}(e,f)$ for a unique pair $e,f$ of twinned edges, and then for $C\in \mathcal{C}$,
$B\in \mathcal{B}(C)$
if and only if $C=D(e,f)$.

Also, we claim that if $B \in {\cal B}(C)$ then $B \cap \eta (G) \subseteq P(C)$; for this is trivial
if $C$ is a circuit, so we assume that $C$ is a path. By (3) the claim holds if $B \in {\cal A} (C)$,
so we may assume that $B=D(e,f)$ for some pair $e,f$ of disjoint twinned edges, and $B\in \mathcal{A}(e,f)$.
But then the claim holds by the third assertion of (2).

Finally, we claim that $P(C) \cup \bigcup (B:B\in {\cal B}(C))$ is
$P(C)$-planar for each $C \in {\cal C}$. If $C=D(e,f)$ for some pair $e,f$ 
with a common end, the claim follows from (1) and the definition 
of $D(e,f)$. If $C=D(e,f)$ for some pair of disjoint twinned edges, the claim follows from 
the third assertion of (2) and the definition of $D(e,f)$, since $\mathcal{A}(C)=\emptyset$ from the first assertion of (2).
And otherwise, the claim follows from \ref{extend}. This proves that 
$\eta$ provides a flattening satisfying the theorem, and so proves \ref{mainthm}.~\bbox

\section{Augmenting paths}

We need three more techniques for the second half of the paper, all
developed in \cite{RST}, and in this section we describe the first.
If $F$ is a subgraph of $G$ and of $H$, and $\eta$ is a
\he of $G$ in $H$, we say it {\em fixes} $F$
if $\eta (e) =e$ for all $e \in E(F)$ and
$\eta (v) =v$ for all $v \in V(F)$.

Let $G$ be cubic, and let $F$ be a subgraph of $G$ with
minimum degree $\geq 2$ (possibly null).
Let $X \subseteq V(G)$, such that $\delta_{G} (X) \cap E(F) = \emptyset$.
Let $n \geq 1$, let $G_{0}  =G$, and inductively for
$1 \leq i \leq n$ let $G_{i} = G_{i-1} +(e_{i},f_{i})$
with new vertices $u_{i},v_{i}$, where $e_{i},f_{i}$ are
edges of $G_{i-1}$ not in $E(F)$.
Let $\eta_{F}$ be the identity \he of $G_{0}$
to itself; and for $1 \leq i \leq n$, let $\eta_{i}$ be obtained
from $\eta_{i-1}$ by replacing
$e_{i}$ and $f_{i}$
by the corresponding two-edge paths of $G_{i}$.
Thus $\eta_{i}$ is a \he of $G$ in $G_{i}$; it
fixes $F$, and
$\eta (v) =v$ for all $v \in V(G)$, and $\eta (e)=e$ for all $e \in E(G)$
except edges of $G$ in $\{e_{1},f_{1} ,\ldots,e_{i},f_{i} \}$.

Let $\delta_{G} (X) = \{x_{1}y_{1},\ldots, x_{k}y_{k} \}$, where
$x_{1},\ldots, x_{k} \in X$
are all distinct, and $y_{1},\ldots,y_{k} \in V(G)  \setminus X$ are all
distinct.
Suppose in addition:
\begin{itemize}
\item
$e_{1} \in E(G)$ has both ends in $X$, and $f_{n} \in E(G)$ with
both ends in $V(G) \setminus X$
\item
for $1 \leq i <n$ there exists $j \in \{1,\ldots,k\}$ such that
$f_{i}$ is the edge of $\eta_{i-1}(x_{j}y_{j})$ incident with
$y_{j}$, and $e_{i+1}$ is the edge of $\eta_{i}(x_{j}y_{j})$
incident with $v_{i}$ and not with $y_{j}$
\item
if $f_{1} \in E( \eta (x_{j}y_{j}))$ (that is, $f_{1}=x_{j}y_{j}$)
where $1 \leq j \leq k$, then $e_{1}$ is not incident with $x_{j}$
in $G$, and no end of $e_{1}$ is adjacent in $G\setminus E(F)$ to $x_{j}$;
similarly, if $e_{n} \in E(\eta (x_{j}y_{j}))$ then $e_{n}$
is not incident with $y_{j}$ in $G$, and no end of $e_{n}$ 
is adjacent in $G\setminus E(F)$ to $y_{j}$
\item
for $2 \leq i \leq n-1$, let $e_{i} \in E(\eta_{i-1} (x_{j}y_{j}))$
and $f_{i} \in E(\eta_{i-1} (x_{j'}y_{j'}))$; then
$j' \neq j$, and $x_{j}$ is not adjacent to $x_{j'}$ in $G\setminus E(F)$,
and $y_{j}$ is not adjacent to $y_{j'}$ in $G\setminus E(F)$.
\end{itemize}
(See Figure 5.)

\begin{figure} [h!]
\centering
\includegraphics{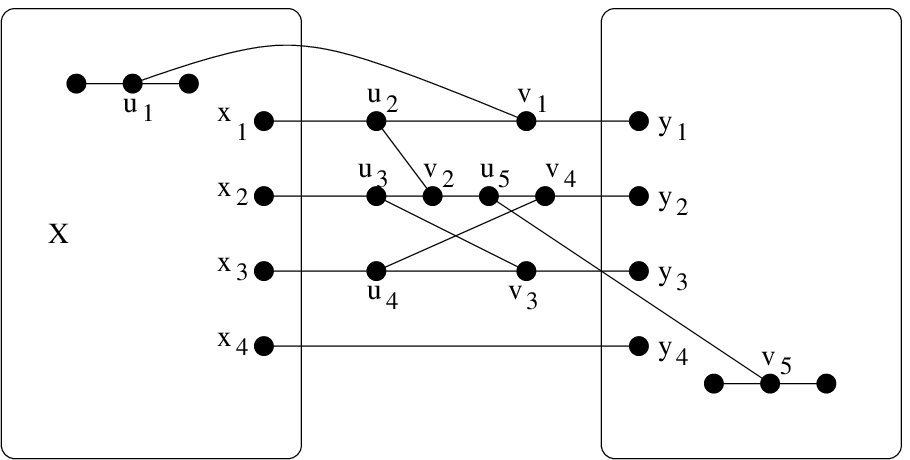}
\caption{An $X$-augmentation of a graph, with $k=4$ and $n=5$.}
\label{augmentationfig}
\end{figure}

In these circumstances we call $G_{n}$ an
$X$-{\em augmentation} of $G$ (modulo $F$), and
$(e_{1},f_{1}),\ldots,(e_{n},f_{n})$ an $X$-{\em augmenting sequence}
of $G$ (modulo $F$).
Note that we permit $n=1$.
The following is proved in lemma 3.4 of~\cite{RST}, applied to $F$, $H\setminus E(F)$ and $X$.

\begin{thm}\label{augment}
Let $G$ be cubic and let $F$ be a subgraph of $G$ with minimum degree
at least two.
Let $X \subseteq V(G)$ with $\delta_{G} (X) \cap E(F) = \emptyset$,
such that the edges in $\delta_{G} (X)$ pairwise have no common end.
Let $H$ be cubic such that $F$ is a subgraph of $H$, and let $\eta$
be a \he of $G$ in $H$ fixing $F$.
Then either
\begin{itemize}
\item
there exists $X' \subseteq V(H)$ with
$|\delta_{H} (X') | = | \delta_{G} (X) |$, such that for
$v \in V(G)$, $v \in X$ if and only if $\eta(v)  \in X'$, or
\item
there is an $X$-augmentation $G'$ of $G$ modulo $F$, and a
\he of $G'$ in $H$ fixing $F$.
\end{itemize}
\end{thm}

\section{Jumps on a dodecahedron}

Now we begin the second part of the paper.
First we prove the following variant of \ref{mccuaig} (equivalent to \ref{nonplanar0}).

\begin{thm}\label{nonplanar}
Let $H$ be cyclically five-connected and cubic.
Then $H$ is non-planar if and only if $H$ contains one of Petersen,
Triplex, Box and Ruby.
\end{thm}
\Proof
``If'' is clear.
For ``only if'', let $H$ be cyclically five-connected and cubic,
and contain none of the four graphs.
By \ref{mccuaig} it follows that $H$ contains Dodecahedron.
Let $G$= Dodecahedron, let $F$ and $\eta_F$ be null, and let ${\cal C}$ be the
set of circuits of $G$ that bound regions in the drawing in Figure 3;
then $(G,F, {\cal C})$ is a framework.
We claim that (E1)--(E7) are satisfied.
Most are trivial, because $F$ is null, and there are no twinned edges, and no paths in ${\cal C}$.
Also, (E6) is vacuously true because no member of ${\cal C}$ has
length $\geq 6$; so the only axiom that needs work is (E2).

Let $e,f \in E(G)$ such that no member of ${\cal C}$ contains both
$e$ and $f$; we claim that $G+(e,f)$ contains
one of Petersen, Triplex, Box, Ruby.
Up to isomorphism of $G$ there are five possibilities for $e,f$, namely
(setting $e=ab$ and $f=cd$)
$(a,b,c,d)=(1,2,6,15)$, $(1,2,10,15), (1,2,15,20)$,
$(1,2,18,19), (1,2,19,20)$.
In the first three cases $G+(e,f)$ contains Ruby, and in the last
two it contains Box.

Thus, (E2) holds; and so $H$ is planar, by \ref{mainthm}.
This proves \ref{nonplanar}.~\bbox

Next, a small repair job.
The definition of ``dodecahedrally-connected'' in \cite{RST} differs from
the one in this paper, and our objective of the remainder of this section is to prove
them equivalent. To do so, we essentially have to repeat the proof of \ref{nonplanar} with slightly different hypotheses.

In this section we fix a graph $F$, and we need to look at several graphs such that $F$ is a subgraph of all of them.
If $G,H$ are cubic, and $F$ is a subgraph of them both, and there is a \he of $G$ in $H$ fixing $F$, we say that $H$ {\em $F$-contains} $G$.

Let $G$ be cubic, and let $F$ be a subgraph of $G$, such that
every vertex in $F$ has degree $\geq 2$ in $F$.
Let $C$ be a circuit of $G$ of length four, with vertices
$a_{1},a_{2},a_{3},a_{4}$ in order, none of them in
$V(F)$.
Let $a_{i}$ be adjacent to $b_{i} \not\in V(C)$ for
$1 \leq i \leq 4$, where $b_{1},\ldots,b_{4}$ are all distinct, and not in
$V(F)$, and are pairwise non-adjacent.
A {\em $C$-leap} of $G$ means a graph $G+(e,f)$, where $e\in E(C)$ and $f\in E(G)\setminus E(F)$, with no vertex in $V(C)$.

\begin{thm}\label{Cleap}
Let $G$  be cubic and cyclically four-connected, with
$|V(G)| \geq 8$.
Let $F$ be a subgraph of $G$ such that every vertex in $F$
has degree $\geq 2$ in $F$.
Let $C$ be a circuit of $G$ of length $4$, disjoint from $F$.
Let ${\cal L}$ be a set of cubic graphs such that $F$
is a subgraph of each of them.
Suppose that every $C$-leap of $G$ $F$-contains a member of ${\cal L}$.
Let $H$ be a cyclically five-connected cubic graph containing
$F$ as a subgraph,
that does not $F$-contain any member of ${\cal L}$.
Then $H$ does not $F$-contain $G$.
\end{thm}
\Proof
Let $X=V(C)$.
Then $\delta_{G} (X) \cap E(F) = \emptyset$ since
$X \cap V(F) = \emptyset$.
Since $G$ is cyclically four-connected and $|V(G)| \geq 8$ it follows
that no two members of $\delta_{G} (X)$ have a common end.

Suppose that $H$ $F$-contains $G$.
Let us apply \ref{augment}.
Since $H$ is cyclically five-connected, \ref{augment}(i) does not hold, and so
\ref{augment}(ii) holds.
Let $(e_{1},f_{1}),\ldots , (e_{n},f_{n} )$ be an $X$-augmenting
sequence of $G$, such that there is a \he of the corresponding
$X$-augmentation $G'$ in $H$ fixing $F$.
From condition (iii) in the definition of
``$X$-augmenting sequence'', it follows that $n=1$,
and so $G'= G+(e_{1},f_{1})$.
Thus $G'$ is a $C$-leap of $G$, and therefore $F$-contains a member of ${\cal L}$.
But $H$ $F$-contains $G'$, and so $H$
$F$-contains a member of  ${\cal L}$, a contradiction.
This proves \ref{Cleap}.~\bbox

It is convenient from now on to make the following convention.
When we speak of a graph $G+(e,f)$ and the vertices of $G$
are numbered $1,\ldots, n$, the new vertices of
$G+(e,f)$ will be assumed to be numbered $n+1$ and
$n+2$ (in order), unless we specify otherwise.

Let $G$ be Dodecahedron, and let $F$ be a circuit of $G$ of length five. If $e,f\in E(G)\setminus E(F)$, and at most one of $e,f$ has an end
in $V(F)$, and $e,f$ are not incident with the same region of $G$, we call $G+(e,f)$ a {\em hop extension} of $(G,F)$; and if in addition $e,f$ are
diverse, we call $G+(e,f)$ a {\em jump extension} of $(G,F)$. 
We begin with the following lemma.

\begin{thm}\label{opposite} Let $G$ be Dodecahedron, and let $F$ be a circuit of $G$ of length five. Let $H$ be a cyclically five-connected cubic graph,
such that $F$ is a subgraph of $H$. Suppose that
\begin{itemize}
\item $H$ $F$-contains no jump extension of $(G,F)$ 
\item for every $X\subseteq V(H)\setminus V(F)$ with $|\delta_H(X)| = 5$ and $X\ne V(H)\setminus V(F)$, there is no \he $\eta$ of $G$ in $H$ fixing $F$ such that
$\eta(v)\in X$ for all $v\in V(G)\setminus V(F)$.
\end{itemize}
Let $e,f$ be diverse edges of $G$ not in $E(F)$; then $H$ does not $F$-contain $G+(e,f)$.
\end{thm}
\Proof Suppose it does. Hence $G+(e,f)$ is not a jump extension of $(G,F)$, and so both $e,f$ have ends in $V(F)$.
Let us number the vertices of Dodecahedron as in Figure 3, and from the symmetry we may assume that $e$ is 2-7
and $f$ is 5-10.
Let $G' = G+(e,f)$ with new vertices 21, 22 say.
Let $X = \{6,7,\ldots,20\}$. From the second bullet and \ref{augment}, there is an $X$-augmenting sequence of $G'$ modulo
$F$, say $(e_{1},f_{1}),\ldots,(e_{n},f_{n})$,
and a \he $\eta''$ of the corresponding $X$-augmentation
$G''$ in $H$ fixing $F$.
Now $e_{1}$ ($= a_{1}b_{1}$ say) has both ends in $X$, but $f_{1}$
does not, so $f_{1}$ is one of 1-6, 2-21, 7-21, 3-8, 4-9, 5-22,
10-22, 21-22; and from the symmetry we may assume that $f_{1}$ is one of
1-6, 2-21, 7-21, 3-8, 21-22.

Suppose that $f_1$ is one of 1-6, 3-8. Then $e_1, f_1$ are diverse, from the third condition in the definition
of $X$-augmenting sequence; but then $G+(e_1,f_1)$ is a jump extension of $(G,F)$ $F$-contained in $G'+(e_1,f_1)$ and hence in $H$, a contradiction. 
Similarly if $f_1$ is 7-21 then $G+(e_1,2\d 7)$ is a jump extension $F$-contained in $H$.
Thus $f_1$ is one of 21-22, 2-21, and in particular $n = 1$. Assume $f_1$ is 21-22. Then we may assume that $e_1$, 2-7 are not diverse in $G$
(for otherwise $G+ (e_1, 2\d 7)$ is a jump extension $F$-contained in $H$), and similarly $e_1$, 5-10 are not diverse in $G$. But this is impossible.
Finally, assume that $f_1$ is 2-21. we may assume that $e_1, 2\d 7$ are not diverse in $G$, and so $e_1$ is one of 
$$7\d 11,7\d 12,6\d 11,11\d 16,8\d 12,12\d 17.$$
If $e_1$ is one of 7-12, 8-12, 12-17, rerouting 7-12 along 21-22 gives a jump extension of $(G,F)$ $F$-contained in $H$; and if $e_1$ is one
of 7-11, 6-11, 11-16, rerouting 7-11 along 21-22 gives a jump extension of $(G,F)$ $F$-contained in $H$, again a contradiction.
This proves \ref{opposite}.~\bbox

\begin{thm}\label{hoptojump}
Let $G,F,H$ be as in \ref{opposite}.
Then $H$ $F$-contains no hop extension of $(G,F)$.
\end{thm}
\Proof
Let ${\cal L}$ be the set of all graphs $G+(e,f)$ where $e,f$ are diverse edges of $G$ not in $E(F)$. By \ref{opposite},
$H$ $F$-contains no member of ${\cal L}$.
Let $G$ be labelled as in Figure 3. (We do not specify the circuit $F$ at this stage; it is better to preserve the symmetry.)
Let $G_1 = G+(a,b)$ be a hop extension of $G$, and suppose that $H$ $F$-contains $G_1$. Thus $G_1\notin {\cal L}$.
From the symmetry of $G$, we may therefore assume that $a$ is 15-20 and $b$ is 16-17. Thus the edges 16-17 and 15-20 are not in $E(F)$. 
Since $F$ is a circuit of length five, it follows that 16-20 is not in $E(F)$, and hence 16,20 are not in $V(F)$.
Let $C$
be the circuit 16-20-21-22-16 of $G_1$. Then no vertex of $C$ is in $V(F)$, and $H$ is cyclically five-connected, so we can apply \ref{Cleap}. 
We deduce that $H$ $F$-contains some $C$-leap $G_2 = G_1+ (e,f)$. 

Now $e$ is one of 16-20, 20-21, 21-22, 16-22.
Since $F$ is not yet specified, there is a symmetry of $G_1$ exchanging the edges 16-20 and 21-22; and one exchanging 20-11 and 16-22.
Thus we may assume that $e$ is one of 21-22, 20-21.

Now $f$ is an edge of $G$ not incident with either of 16,20. Since $e$ is one of 21-22, 20-21, and $f\notin E(F)$, $H$ $F$-contains $G+(15\d 20,f)$ in $G_2$, 
and so $G+(15\d 20,f)\notin {\cal L}$. Consequently $f$, 15-20 are not diverse, 
so $f$ is one of
$$6\d 15, 10\d 15, 1\d 6, 6\d 11, 5\d 10, 10\d 14, 14\d 19, 18\d 19.$$

Suppose first that $e$ is 21-22. Then by the same argument, $f$ and 16-17 are not diverse in $G$, and so
$f$ is one of 6-11, 18-19. If $f$ is 6-11, rerouting 6-15 along 24-23-21 gives a member of ${\cal L}$ $F$-contained in $H$ (in future we just say ``works'')
and if
$f$ is 18-19, rerouting 17-18 along 22-23-24 works. Thus the claim holds if $e$ is 21-22.

Now we assume that $e$ is 20-21. If $f$ is one of 1-6,6-11,6-15 then rerouting 6-15 along 23-24 works;
if $f$ is one of 10-15,5-10,10-14, rerouting 10-15 along 23-24 works; and if $f$ is 14-19 or 18-19 then
rerouting 19-20 along 23-24 works. Thus us each case we have a contradiction. This proves \ref{hoptojump}.~\bbox

Next we need another similar lemma. Let $G$ be Dodecahedron, labelled as in Figure 3, and let $G_1$ be $G+(1\d 6,2\d 7)$. Let $G_2 = G_1 + (6\d 21,2\d 22)$.
Thus the edge 1-6 of $G$ has been subdivided to become a path 1-21-23-6 of $G_2$, and 2-7 has become 2-24-22-7.

\begin{thm}\label{crossext}
Let $G,F,H$ be as in opposite. Then $H$ does not $F$-contain $G_2$.
\end{thm}
\Proof Let $X = \{6,7\l 20\}$. By the second bullet hypothesis about $H$, and \ref{augment},
there is an $X$-augmenting sequence of $G_2$ modulo
$F$, say $(e_{1},f_{1}),\ldots,(e_{n},f_{n})$,
and a \he $\eta'$ of the corresponding $X$-augmentation
$G'$ in $H$ fixing $F$.
Now $e_{1}$ ($= a_{1}b_{1}$ say) has both ends in $X$, but $f_{1}$
does not, so $f_{1}$ is one of 
$$1\d 21, 21\d 23, 6\d 23, 2\d 24, 22\d 24, 7\d 22, 3\d 8, 4\d 9, 5\d 10, 21\d 22, 23\d 24,$$
and from the symmetry we may assume that $f_1$ is one of 
$$1\d 21, 21\d 23, 6\d 23, 5\d 10, 4\d 9, 21\d 22.$$
If $f_1$ is one of 5-10, 4-9 then by the third condition in the definition of  $X$-augmenting sequence,
it follows that $e_1,f_1$ are diverse in $G$, and $H$ contains the jump extension $G+(e_1,f_1)$, a contradiction.
Similarly if $f_1$ is 6-23 then $e_1, 1\d 6$ are diverse in $G$, again a contradiction.
Thus $f_1$ is one of 1-21, 21-23, 21-22.
Hence $H$ $F$-contains $G+(1\d 6, e_1)$, and so by \ref{hoptojump}, $G+(1\d 6, e_1)$ is not a hop extension of $(G,F)$.
Consequently $f_1$ is one of 10-15, 6-15, 6-11, 7-11.
If $f_1$ is one of 6-11, 7-11, then rerouting 1-6 along 25-26 gives a jump extension of $(G,F)$ $F$-contained in $H$; while
if $f_1$ is one of 6-15, 10-15, rerouting 6-15 along 25-26, and then rerouting 7-11 along 23-24, give the desired jump extension.
This proves \ref{crossext}.~\bbox

From these lemmas we deduce a kind of variant of \ref{nonplanar}:

\begin{thm}\label{newstrange}
Let $G$ be Dodecahedron, and let $F$ be a circuit of $G$ of length five. Let $H$ be a cyclically five-connected cubic graph,
such that $F$ is a subgraph of $H$. Suppose that
\begin{itemize}
\item $H$ $F$-contains no jump extension of $(G,F)$
\item for every $X\subseteq V(H)\setminus V(F)$ with $|\delta_H(X)| = 5$ and $X\ne V(H)\setminus V(F)$, there is no \he $\eta$ of $G$ in $H$ fixing $F$ such that
$\eta(v)\in X$ for all $v\in V(G)\setminus V(F)$.
\end{itemize}
Then $H$ is planar, and can be drawn in the plane such that $F$ bounds the infinite region.
\end{thm}
\Proof
Let ${\cal C}$ be the set of the following eleven subgraphs
of $G =$ Dodecahedron; the six circuits that bound
regions (in the drawing in Figure 3) that contain no edge incident
with the infinite region,
and for each $e \in E(F)$, the path $C\setminus e$ where
$C \neq F$ is the boundary of a region incident with $e$.
Let $\eta_F$ be the identity \he on $F$.
By hypothesis there is a \he of $G$ in $H$ extending $\eta_F$.
We apply \ref{mainthm} to $(G,F, {\cal C})$ and $H, \eta_{F}$.
There are no twinned edges and all members of ${\cal C}$ have
at most five edges; so we have to check only (E2) and (E6). (Note that in this case,
the paths in ${\cal C}$ are not induced subgraphs of $G$; this is the only one of our applications when this is so.)
But the truth of (E2) and (E6) follows from the three lemmas above~\ref{opposite}, \ref{hoptojump}, \ref{crossext}; and so by \ref{mainthm}, the result follows.
This proves \ref{newstrange}.~\bbox

As we said earlier, we need this to prove the equivalence of the definitions of dodecahedrally-connected given in this paper and in~\cite{RST},
and now we turn to that. 
Let $G$ be Dodecahedron, and let $F$ be a circuit of $G$ of length five.
Let $H$ be a cubic graph, and let $X \subseteq V(H)$.
We say that $H$ is {\em placid} on $X$ if 
\begin{itemize}
\item $|V(H)  \setminus X | \geq 7$, and $\delta_{H} (X)$ is a matching of cardinality five
\item $\{x_{i}y_{i}: 1 \leq i \leq 5 \}$  is an enumeration of $\delta_{H} (X)$, with $x_{i} \in X\;(1 \leq i \leq 5 )$
\item there is a \he of $G$ in $H'$ mapping $F$ to the circuit $y_1\d y_2\d y_3\d y_4\d y_5\d y_1$, and
\item
there is no \he of any jump extension of $(G,F)$ in $H'$ mapping $F$ to $y_1\d y_2\d y_3\d y_4\d y_5\d y_1$,
\end{itemize}
where $H'$ is obtained from $H|(X\cup \{y_1,y_2,y_3,y_4,y_5\})$ by deleting all edges with both ends in $\{y_1,y_2,y_3,y_4,y_5\}$,  
and adding 
new edges $y_1y_2, y_2y_3, y_3y_4,y_4y_5,y_1y_5$.

We say that a graph $H$ is {\em strangely connected} if $H$ is cubic and
cyclically five-connected, and there is no $X\subseteq V(H)$ such that $H$ is placid on $X$.
(This is the definition of ``dodecahedrally-connected'' in~\cite{RST}.)

\begin{thm}\label{strange}
A graph $H$ is dodecahedrally-connected if and only if it is strangely connected.
\end{thm}
\Proof
We may assume that $H$ is cubic and cyclically five-connected.
Suppose first that it is not  dodecahedrally-connected.
Let $X \subseteq V(H)$ with $|X|, |V(H) \setminus X| \geq 7$ and
$|\delta_{H} (X) | =5$, $\delta_{H} (X) = \{x_{1}y_{1},\ldots,x_{5}y_{5}\}$
say where $x_{1},\ldots,x_{5} \in V(H)$, such that $H|X$ can be drawn
in a disc with $x_{1},\ldots,x_{5}$ on the boundary in order.
Let us choose such $X$ with $|X|$ minimum.
Since $H$ is cyclically five-connected it follows that
$x_{1},\ldots,x_{5}$ are all distinct and so are
$y_{1},\ldots,y_{5}$.
Also, from the planarity of $H|X$ it follows that $|X| \geq 9$,
and so from the minimality of $X$, no two of $x_{1},\ldots,x_{5}$
are adjacent.
Let $H'$ be obtained from $H$ as in the definition of ``placid'', and
let $F'$ be the circuit made by the five new edges.
It follows easily that $H'$ is cyclically five-connected, and hence
from \ref{nonplanar0} contains $G= $ Dodecahedron.
Take a planar drawing of $H'$, and choose a \he $\eta$ of $G$
in $H'$ such that the region of $\eta (G)$ including $r$ is minimal,
where $r$ is the region of $H'$ bounded by $F'$.
It follows easily that $F' \subseteq \eta (G)$, and so from
the symmetry of $G$ we may choose $\eta$ mapping $F$ to $F'$.
Hence $H$ is placid on $X$ 
(the final condition in the definition of ``placid'' holds because of the planarity of $H'$) and so $H$ is not
strangely connected, as required.

For the converse, suppose that $H$ is not strangely connected, and let
$X$, $x_{i}y_{i}\; (1 \leq i \leq 5 )$, $F$ and $H'$
be as in the definition of ``strangely connected'', such that $H$
is placid on $X$ via $x_{1}y_{1},\ldots,x_{5}y_{5}$.
Choose $X$ minimal.
By \ref{newstrange}, $H|X$ can be drawn
in a disc with $x_{1},\ldots,x_{5}$ on the boundary in order;
and so $H$ is not dodecahedrally-connected.
This proves \ref{strange}.~\bbox

\section{Adding jumps to repair connectivity}

Now that we have reconciled the two definitions of ``dodecahedrally- connected'', 
we can apply results of \cite{RST} about this kind of connectivity.

The idea behind \ref{Cleap} is that cyclic five-connectivity is better than cyclic four-connectivity,
and we begin with a graph $G$ that is cyclically five-connected, except for the circuit $C$. We use
the cyclic five-connectivity of $H$ to prove that if $H$ contains $G$ then $H$ also contains a slightly larger graph
where the circuit $C$ has been expanded to a circuit of length five by adding an edge to $G$. This can be useful,
as we saw in the previous section. However, it has the defect that the edge we add to $G$ to expand the circuit $C$ might create
a new circuit of length four, with its own problems. We can apply \ref{Cleap} again to this new circuit, but the process
can go on forever. In fact, there is a stronger theorem; one can expand the circuit $C$ to a longer circuit, without adding
any new circuits of length four, just by adding a bounded number of edges. That is essentially
the content of the next result, proved in \cite{RST}. (We also weaken the hypothesis on $G$, allowing it to have
more than one circuit of length four.)
But first we need some definitions.

Let ${\cal L}$ be a set of cubic graphs.
We say that a graph $H$ is {\em killed by}
${\cal L}$ if there is a \he of some $G' \in {\cal L}$ in $H$.
Let $G$ be cubic, and let
$C$ be a circuit of $G$ of length four, with vertices
$a_{1},a_{2},a_{3},a_{4}$ in order.
Let $a_{i}$ be adjacent to $b_{i} \not\in V(C)$ for
$1 \leq i \leq 4$, where $b_{1},\ldots,b_{4}$ are all distinct
and pairwise non-adjacent.
We denote by ${\cal P} (C, {\cal L})$ the set of all pairs $(e,f)$ such that
$f \in E(G)$ is incident with one of $b_{1},\ldots,b_{4}$, say $b_{i}$,
$f \neq a_{i}b_{i}, e \in E(C)$ is incident with $a_{i}$, and
$G+(e,f)$ is not killed by ${\cal L}$.

Let $e=uv$ and $f=wx$ be edges of a cubic graph $G$.
If $u,v \neq w,x$, and $u$ is adjacent to $w$, and no other edge has
one end in $\{u,v\}$ and the other in $ \{w,x\}$, we denote by
$(e,f)^{\ast}$ the pair of edges $(e',f')$, where
$e' \,(\neq e,uw)$ is incident with $u$ and
$f' \, (\neq f, uw)$ is incident with $w$.

We shall frequently have to list the members of some set
${\cal P}(C, {\cal L})$ explicitly, and we can save some writing as follows.
Clearly $(e,f) \in {\cal P} (C, {\cal L})$ if and only if
$(e,f)^{\ast} \in {\cal P} (C, {\cal L})$, and so we really need only
to list half the members of ${\cal P} (C, {\cal L})$.
If $X$ is a set of pairs of edges for which $(e,f)^{\ast}$
is defined for each $(e,f) \in X$, we denote by $X^{\ast}$
the set $X \cup \{(e,f)^{\ast} :(e,f) \in X \}$.

If $e \in E(C)$ and $e,f$ are diverse in $G$, we call
$G+(e,f)$ an {\em A-extension} of $G$.
Now let $e \in E(C)$ and $f \in E(G) \setminus E(C)$
such that $e,f$ are not diverse in $G$ but have no common end.
Let $G' = G+(e,f)$ with new vertices $x_{1},y_{1}$.
Label the vertices of $C$ as $a_{1},\ldots,a_{4}$ in order, and their
neighbours not in $V(C)$ as $b_{1},\ldots,b_{4}$ respectively, as before,
such that $e=a_{1}a_{2}$ and $f$ is incident with
$b_{1},f=b_{1}c_{1}$ say.
If $g \in E(G)$, not incident in $G$ with
$a_{1},b_{1},c_{1},d_{1}$ (where $b_{1}$ is adjacent in $G$ to
$a_{1},c_{1},d_{1}$) we call $G' + (b_{1}y_{1},g)$ a
$B$-{\em extension} ({\em of G}) {\em via} $(e,f)$.
If $g \in E(G)$ incident with $b_{2}$ and not with $c_{1}$
or $a_{2}$, we call $G' + (x_{1}y_{1},g)$ a $C$-{\em extension via}
$(e,f)$ {\em onto g}.
We call $G' + (a_{1}x_{1},a_{3}b_{3})$ a $D$-{\em extension via}
$(e,f)$.
Finally, we say $(e,f)$ and $(e', f')$ are $C$-{\em opposite} if
$e,e' \in E(C)$ and the labelling can be chosen as before with
$e=a_{1}a_{2}$, $f=b_{1}c_{1}$, $e'=a_{3}a_{4}$, and $f'=b_{3}c_{3}$.
Let $(e,f), (e',f')$ be $C$-opposite, with labels as above.
Let $G'' = G' + (e', f')$ with new vertices $x_{2},y_{2}$;
then we call $G'' + (a_{1}x_{1},a_{3}x_{2})$ an
$E$-{\em extension via} $(e,f), (e',f')$.

We say a graph $G$ is {\em quad-connected} if
\begin{itemize}
\item $G$ is cubic and cyclically four-connected
\item $|V(G)|\ge 10$, and if $G$ has more than one circuit of length four then
$|V(G)|\ge 12$
\item for all $X\subseteq V(G)$ with $|\delta_G(X)|\le 4$, one of $|X|,|V(G)\setminus X|\le 4$.
\end{itemize}

The following is a restatement of \ref{Cleap} in this language 
(with $F$ removed, because we no longer need it.)
\begin{thm}\label{Aextn}
Let $G$  be cubic and cyclically four-connected, with
$|V(G)| \geq 8$.
Let $C$ be a circuit of $G$ of length $4$, and let ${\cal L}$ be a set of cubic graphs.
Suppose that  every A-extension of $G$ is killed by ${\cal L}$,
and ${\cal P} (C, {\cal L}) = \emptyset$.
Let $H$ be a cyclically five-connected cubic graph
that is not killed by ${\cal L}$.
Then there is no \he of $G$ in $H$.
\end{thm}

Here is the strengthening, proved in~\cite{RST}.
\begin{thm}\label{extn}
Let $G$ be quad-connected, 
and let $C$ be a circuit of $G$ of length four.
Let ${\cal L}$ be a set of cubic graphs, such that
\begin{itemize}
\item
every A-extension of $G$ is killed by ${\cal L}$
\item
for every $(e,f) \in {\cal P}(C, {\cal L})$, every B-extension via
$(e,f)$ is killed by ${\cal L}$, and so is the D-extension via $(e,f)$
\item
for all $(e,f_{1}),(e,f_{2}) \in {\cal P} (C, {\cal L})$ such that
$f_{1},f_{2}$ have no common end, the C-extension via
$(e,f_{1})$ onto $f_{2}$ is killed by ${\cal L}$, and
\item
for all $C$-opposite $(e_{1},f_{1}),(e_{2},f_{2}) \in {\cal P} (C, {\cal L})$,
the E-extension via
$(e_{1},f_{1}), (e_{2},f_{2})$ is killed by ${\cal L}$.
\end{itemize}
Let $H$ be a dodecahedrally-connected cubic graph such that $H$ is not killed by ${\cal L}$.
Then there is no \he of $G$ in $H$.
\end{thm}

The other result of~\cite{RST} that we need is the following.
Let $n \geq 5$ be an integer, with $n \geq 10$ if $n$ is even.
The $n$-{\em biladder} is the graph with vertex set
$\{a_{1},\ldots,a_{n},b_{1},\ldots,b_{n}\}$, where for
$1 \leq i \leq n, a_{i}$ is adjacent to $a_{i+1}$ and to $b_{i}$,
and $b_{i}$ is adjacent to $b_{i+2}$ (where
$a_{n+1}, b_{n+1}, b_{n+2}$ mean $a_{1},b_{1},b_{2}$).
Thus, Petersen is isomorphic to the 5-biladder, and Dodedahedron to the
10-biladder. The following follows from theorem 1.4 of~\cite{RST}.

\begin{thm}\label{splitter}
Let $G$ be cubic and cyclically five-connected.
Let there be a \he of $G$ in $H$, where $H$ is dodecahedrally-connected.
Then either
\begin{itemize}
\item
there exist $e,f \in E(G)$, diverse in $G$, such that there is a
\he of $G+(e,f)$ in $H$, or
\item
$G$ is isomorphic to an $n$-biladder for some $n$, and there is a \he of the
$(n+2)$-biladder in $H$, or
\item
$G$ is isomorphic to $H$.
\end{itemize}
\end{thm}

\section{Graphs with crossing number at least two}

At the end of the proof of \ref{nonplanar}, there were five statements left
to the reader to verify, that five particular graphs contain either
Ruby or Box.
In the remainder of the paper there will be many more similar statements
left to the reader; unfortunately, we see no way of avoiding this,
since there are simply too many of them to include full details of each.
But perhaps 95\% of them are of the form that
``Graph $G$ contains Petersen'', where $G$ is cubic and cyclically
five-connected; and here is a quick method for checking such a statement.
Choose a circuit $C$ of $G$ with $|E(C) | =5$, arbitrarily
(there always is one, in this paper).
Let $C$ have vertices $v_{1},\ldots,v_{5}$ in order.
Let $u_{1},\ldots,u_{5}$ be vertices of a 5-circuit of Petersen, in order.
Check if there is a \he $\eta$ of Petersen in $G$ with
$\eta (u_{i}) =v_{i} \, (1 \leq i \leq 5)$.
(This is easy to do by hand.)
It is proved in \cite{Truemper} that such a \he exists if
and only if $G$ contains Petersen.

This makes checking for containment of Petersen much easier.
But even so, there are too many cases to reasonably do them all
by hand, and we found it very helpful to write a simple computer
programme to check containment for us. 
We suggest that the reader who wants to check these cases
should do the same thing. There is a computer file available
online with all the details of the case-checking~\cite{RSTappendix}.

In this section, we prove \ref{crossingno0}, which we restate as:

\begin{thm}\label{crossingno}
Let $H$ be dodecahedrally-connected.
Then $H$ has crossing number $\geq 2$ if and only if it contains one of
Petersen, Triplex or Box.
\end{thm}

Dodecahedral connectivity cannot be replaced by cyclic
5-connectivity, because the graph of Figure 6 is a counterexample.
\begin{figure} [h!]
\centering
\includegraphics{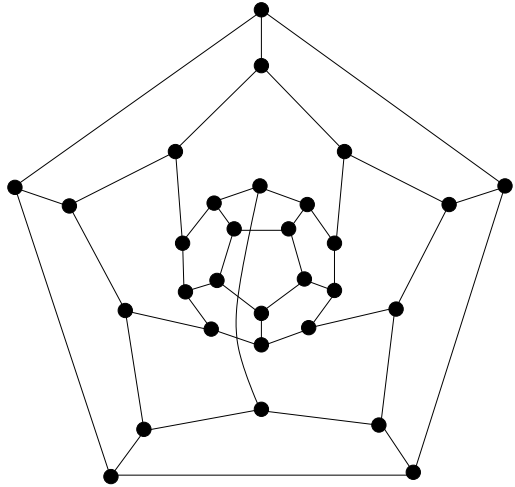}
\caption{A counterexample to a strengthening of \ref{crossingno}.}
\label{counterexamplefig}
\end{figure}
The graphs Window, Antibox, and Drape are defined in Figure 7.
\begin{figure} [h!]
\centering
\includegraphics{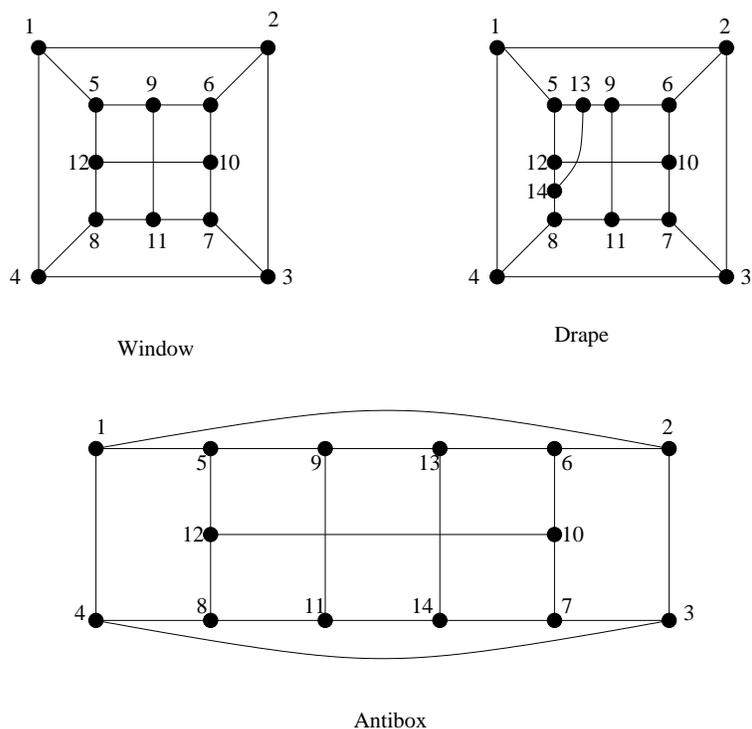}
\caption{Window, Drape and Antibox.}
\label{windowfig}
\end{figure}
We prove \ref{crossingno} in three steps, as follows.

\begin{thm}\label{antibox}
Let $H$ be a dodecahedrally-connected graph containing  Antibox;
then $H$ contains Petersen, Triplex or Box.
\end{thm}

\begin{thm}\label{drape}
Let $H$ be a cyclically five-connected cubic graph containing Drape;
then $H$ contains Petersen, Triplex, Box or Antibox.
\end{thm}

\begin{thm}\label{window}
Let $H$ be a cyclically five-connected cubic graph containing
Window, but not Petersen, Triplex, Box, Antibox or Drape.
Then $H$ has crossing number $\leq 1$.
\end{thm}

\noindent {\bf Proof of \ref{crossingno}, assuming \ref{antibox}, \ref{drape}, \ref{window}.}
``If'' is clear and we omit it.
For ``only if'', let $H$ be dodecahedrally-connected, and contain
none of Petersen, Triplex or Box.
By \ref{antibox} it does not contain Antibox, and by \ref{drape} it does not
contain Drape.
We may assume from \ref{nonplanar} that it contains Ruby (in fact it must,
for no dodecahedrally-connected graph is planar), and
hence Window, since Ruby contains Window.
From \ref{window}, this proves \ref{crossingno}.~\bbox

\bigskip

\noindent{\bf Proof of \ref{antibox}.}
We shall apply \ref{extn}, with $G$ = Antibox, $C$ the
quadrangle of $G$, and ${\cal L} =$ $\{$Petersen, Triplex, Box$\}$.
Thus, $V(C) = \{1,2,3,4\}$.
We find that every $A$-expansion is killed by ${\cal L}$.
In detail, let $G'$ be $G+(ab, cd)$, where $(a,b,c,d)$
is as follows; in each case $G'$ contains  the specified member of
${\cal L}$.
\\
\\
Petersen:
(1, 2, 7, 10), (1, 2, 7, 14), (1, 2, 8, 11), (1, 2, 8, 12),
(1, 2, 9, 11), (1, 2, 11,14), (1, 2, 13, 14), (1, 4, 6, 10),
(1, 4, 6, 13), (1, 4, 7, 10), (1, 4, 7, 14), (1, 4, 9, 13),
(1, 4, 11, 14), (1, 4, 13, 14).
\\
\\
Triplex:
(1, 2, 5, 12), (1, 2, 6, 10), (1, 2, 10, 12), (1, 4, 5, 9),
(1, 4, 8, 11), (1, 4, 9, 11).
\\
\\
Box: (1, 2, 9, 13), (1, 4, 10, 12).
\\
\\
In future we shall omit this kind of detail (because in the future
it will get worse). The full details are in \cite{RSTappendix}.

We find that ${\cal P} (C, {\cal L}) =$
$\{(1\d 2,5\d 9), (1\d 2, 6\d 13),$
$(3\d 4, 8\d 11), (3\d 4,7\d 14)\}^{\ast}$.
Then we verify the hypotheses (ii)-(iv) of \ref{extn}.
This proves \ref{antibox}.~\bbox

\bigskip
\noindent{\bf Proof of \ref{drape}.}
We apply \ref{Aextn}, with $G$ = Drape, $C$ the quadrangle
of $G$ with vertex set $\{5, 12, 13, 14\}$, and
${\cal L} =$ $\{$Petersen, Triplex, Box, Antibox$\}$.
We find that every $A$-extension of $G$ is killed by ${\cal L}$, and
${\cal P} (C, {\cal L}) = \emptyset$, so from \ref{Aextn}, this proves \ref{drape}.~\bbox

\bigskip
\noindent{\bf Proof of \ref{window}.}
Let $G$ be Window, let $F$ and $\eta_F$ be null, and let
${\cal C}$ be the subgraphs of $G$ induced on the following nine sets:
\begin{eqnarray*}
&1, 2, 3, 4;&\\
 &1, 2, 5, 6, 9;&\\
 &2, 3, 6, 7, 10;&\\
&3, 4, 7, 8, 11;& \\
&1, 4, 5, 8, 12;&\\
& 5, 9, 10, 11, 12;&\\
&6, 9, 10, 11, 12;&\\
&7, 9, 10, 11, 12;&\\
&8, 9, 10, 11, 12.&
\end{eqnarray*}
Then $(G, F, {\cal C})$ is a framework.
We claim that (E1)--(E7) hold.
The only twinned edges are 9-11 and 10-12, and
again the only axiom that needs work is (E2).
But if $e,f \in E(G)$ are not both in some member of ${\cal C}$, then
$G+(e,f)$ contains one of Petersen, Triplex, Box, Antibox, Drape,
and so (E2) holds.
From \ref{mainthm}, this proves \ref{window}.~\bbox

\section{Non-projective-planar graphs}

Now we digress, to prove a result that we shall not need;
but it is pretty, and follows easily from the machinery we
have already set up.
The graph {\em Twinplex} is defined in Figure 8.
We shall show the following.

\begin{thm}\label{projp}
Let $H$ be dodecahedrally-connected.
Then $H$ cannot be drawn in the projective plane if and only if
$H$ contains one of Triplex, Twinplex, Box.
\end{thm}
\Proof
``If'' is easy and we omit it.
For ``only if'', suppose that $H$ contains none of Triplex,
Twinplex, Box; we shall show that it can be drawn in the
projective plane.
If $H$ has crossing number $\leq 1$ this is true, so by \ref{crossingno}
we may assume that $H$ contains Petersen.

\begin{figure} [h!]
\centering
\includegraphics{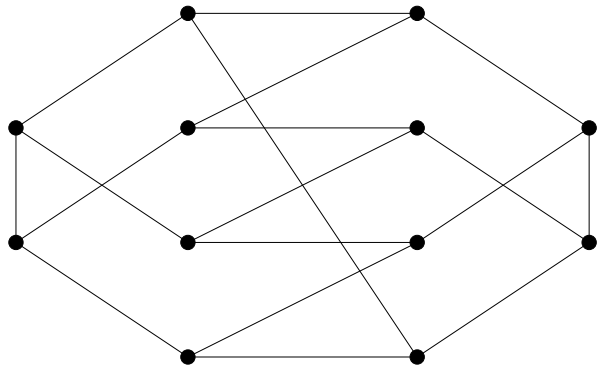}
\caption{Twinplex.}
\label{twinplexfig}
\end{figure}

Let $G_{0}$ = Petersen.
We may assume that $H$ is not isomorphic to $G_{0}$, so by \ref{splitter}
either there are edges $ab, cd$ of $G_{0}$ diverse in $G_{0}$
and a \he of $G_{0} +(ab,cd)$ in $H$, or $H$ contains the
$7$-biladder.
The former is impossible, because from the symmetry of $G_{0}$
we may assume that $(a,b,c,d) = (4,5,6,8)$, and then
$G_{0} + (ab,cd)$ is isomorphic to Twinplex, a contradiction.
Hence there is a \he of $G$ in $H$, where $G$
is the $7$-biladder.
Let $V(G) = \{a_{1},\ldots,a_{7},b_{1},\ldots,b_{7}\}$, as
in the definition of ``biladder''.
Let ${\cal C}$ be the subgraphs of $G$ induced on the following
vertex sets:
\begin{eqnarray*}
&b_{1},b_{2},\ldots,b_{7};&\\&
a_{1},a_{2},a_{3},b_{3},b_{1};&\\&
a_{2},a_{3},a_{4},b_{4},b_{2};&\\&
a_{3},a_{4},a_{5},b_{5},b_{3};&\\&
a_{4},a_{5},a_{6},b_{6},b_{4};&\\&
a_{5},a_{6},a_{7},b_{7},b_{5};&\\&
a_{6},a_{7},a_{1},b_{1},b_{6};&\\&
a_{7},a_{1},a_{2},b_{2},b_{7}.&
\end{eqnarray*}
(These are the face-boundaries of an embedding of $G$ in the
projective plane.)
Let $F$ and $\eta_F$ be null; then
$(G,F, {\cal C})$ is a framework, and we claim that
(E1)--(E7) hold.
All except (E2), (E3) and (E6) are obvious.
To check (E2), let $G'=G+(ab,cd)$ where
$ab,cd \in E(G)$ are not both in any member of ${\cal C}$.
There are twelve possibilities for $(a,b,c,d)$ up to
isomorphism of $G$; in one case $G'$ contain Box, in three
others it contains Twinplex, and in the other eight it contains
Triplex. (As usual, we omit the details; they are also not in the 
appendix~\cite{RSTappendix}, because we don't really need the result.)
Thus, (E2) holds.
For (E3), the only diverse trinity (up to isomorphism of $G$)
is $\{a_{1}a_{2},b_{1}b_{3}, b_{2}b_{7}\}$, and
$G+(a_{1}a_{2},b_{1}b_{3}, b_{2}b_{7})$ contains Twinplex.
Hence (E3) holds.
For (E6), we need only check cross extensions over the circuit
with vertex set $\{b_{1},\ldots,b_{7}\}$, since all other
members of ${\cal C}$ have only five edges.
There are four possibilities (up to isomorphism of $G$).
Let $G' = G+(b_{1}b_{3},b_{2}b_{4})$ with new vertices $x,y$;
then the possibilities are $G'+(ab,cd)$ where
$(a,b,c,d)$ is  $(b_{1},x,b_{2},y)$,
$(b_{1},x,b_{2},b_{7})$, $(b_{1},b_{6},b_{2},b_{7})$,
$(b_{1},b_{6},b_{5},b_{7})$.
The first contains Box, and the other three contain Triplex.
Hence (E6) holds, and from \ref{mainthm}, this proves \ref{projp}.~\bbox

\section{Arched graphs}

We say a graph $H$ is {\em arched} if $H\setminus e$
is planar for some edge $e$.
In this section we prove \ref{arched0}, which we restate as:

\begin{thm}\label{arched}
Let $H$ be dodecahedrally-connected.
Then $H$ is arched if and only if it does not contain Petersen or Triplex.
\end{thm}

We start with the following lemma.

\begin{thm}\label{box}
Let $G$ be Box, let $G'$ be obtained by deleting the edge 13-14, and let ${\cal C}$ be the set of circuits of
$G'$ that bound regions in the drawing in Figure 3.
Let $e,f \in E(G)$, with no common end, and not both in any
member of ${\cal C}$.
Then either $G+(e,f)$ has a Petersen or Triplex minor, or
(up to exchanging $e$ and $f$, and automorphisms of $G$)
$e$ is 13-14 and $f$ is 1-2 or 1-4.
\end{thm}

We leave the proof to the reader (the details are in the Appendix~\cite{RSTappendix}).

\begin{thm}\label{box3}
Let $G$ be Box, and let $H$ be cyclically five-connected, and not contain
Petersen or Triplex.
Let $\eta$ be a \he of $G$ in $H$ such that $\eta (13\d 14)$
has only one edge, $g$ say.
Then $H \setminus g$ is planar, and so $H$ is arched.
\end{thm}
\Proof
We apply \ref{mainthm}, taking $F$
to be the subgraph of $G$ consisting of 13-14 and its ends,
and $\eta_F$ the restriction of $\eta$ to $F$.
Let ${\cal C}$ be as in \ref{box}.
Then $(G,F, {\cal C})$ is a framework, and we claim that
(E1)--(E7) hold.
(E2) follows from \ref{box}, and (E5) and (E6) are vacuously true, because
all members of ${\cal C}$ have five edges. Also, (E3) and (E7) are vacuously true.
For (E4), it suffices from symmetry to check
\begin{eqnarray*}
&G+(1\d 2,13\d 14) + (3\d 6,13\d 16)&\\&
G+(1\d 2,13\d 14) + (3\d 6,14\d 16)&\\&
G+(1\d 2,13\d 14) + (5\d 6,13\d 16)&\\&
G+(1\d 4,13\d 14) + (3\d 6,13\d 16),&
\end{eqnarray*}
but all four contain Triplex.
Hence (E4) holds, so from \ref{mainthm}, this proves \ref{box3}.~\bbox

The graph Superbox is defined in Figure 9.
(It is isomorphic to Box + $(1\d 4,13\d 14)$.)

\begin{figure} [h!]
\centering
\includegraphics{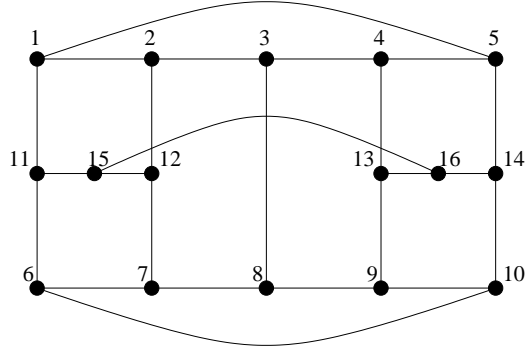}
\caption{Superbox.}
\label{superboxfig}
\end{figure}

\begin{thm}\label{superbox}
Let $G$ be Superbox, let $G'$ be obtained by deleting the edge 15-16, and let ${\cal C}$ be the set of circuits of
$G'$ that bound regions in the drawing in Figure 9.
Let $e,f \in E(G)$ with no common end, and not both in any member
of ${\cal C}$.
Then either $G+(e,f)$ has a Petersen or Triplex minor, or
(up to exchanging $e, f$ and automorphisms of $G$) $e$ is
15-16 and $f$ is 1-2 or 1-11.
\end{thm}

We leave the proof to the reader.
(Actually, it follows quite easily from \ref{box}.)

\begin{thm}\label{superbox5}
Let $G$ be Superbox, and let $H$ be cyclically five-connected, and
not contain Petersen or Triplex.
Let $\eta$ be a \he of $G$ in $H$ such that
$\eta (15\d 16)$ has only one edge, $g$ say.
Then $H \setminus g$ is planar, and so $H$ is arched.
\end{thm}
\Proof
We apply \ref{mainthm} to $(G,F, {\cal C})$, where
$F$ consists of
15-16 and its ends, and $\eta_F$ is the restriction of $\eta$ to $F$, and ${\cal C}$ is as in
\ref{superbox}.
Because of \ref{superbox}, it remains to verify (E4), (E5) and (E6), because (E3), (E7) are vacuous.
Checking (E4) is exactly like in \ref{box3} (indeed, by deleting
14-16 from $G$ we obtain Box, so actually we could
deduce that (E4) holds now from the fact that it held in the proof
of \ref{box3}).
For (E5), we must check
\[
G+(1\d 11,15\d 16) + (6\d 11,15\d 18 )+
(ab,cd)
\]
where $(ab,cd)$ is either $(11\d 17,10\d 14)$ or
$(11\d 19,5\d 14)$; and both contain Triplex.
Thus (E5) holds.
For (E6), we need only check cross extensions over the circuit
bounding the infinite region, since all other members of
${\cal C}$ have length five; and from symmetry, it suffices to check
\begin{eqnarray*}
&G+(1\d 11,10\d 14) + (1\d 17,10\d 18)&\\&
G+(1\d 11,10\d 14) + (6\d 11,5\d 14)&\\&
G+(1\d 11,10\d 14) + (1\d 5,6\d 10)&\\&
G+(1\d 5,6\d 10) + (1\d 17,10\d 18).&
\end{eqnarray*}
All four contain Petersen.
Hence (E6) holds, and from \ref{mainthm}, this proves \ref{superbox5}.~\bbox

\bigskip

\noindent{\bf Proof of \ref{arched}.}
``Only if'' is easy and we omit it.
For ``if'', let $H$ be dodecahedrally-connected, and not
contain Petersen or Triplex.
Since graphs of crossing number $\leq 1$ are arched, we may assume
from \ref{crossingno} that $G$ contains Box.
Choose a \he of $G$ in $H$, where $G$ is either Box or Superbox,
such that $|E(S)|$ is minimum, where
$S = \eta (15\d 16)$ if $G$ is Box, and $S= \eta (17\d 18)$
if $G$ is Superbox.
We claim that $|E(S)| =1$.
For suppose not.
Since $H$ is three-connected, there is an $\eta$-path $P$ with one end in $V(S)$
and the other, $t$, in $V(\eta (G)) \setminus V(S)$.
Let $t \in \eta (f)$ say, and let
$e=15\d 16$ if $G$ is Box, and $e=17\d 18$ if $G$ is Superbox.
If $e,f$ have a common end in $G$, then by rerouting $f$ along $P$
we contradict the minimality of $|E(S)|$.
If some edge $g$ of $G$ joins an end of $e$ to an end of $f$,
then by rerouting $g$ along $P$ we contradict the minimality of $|E(S)|$.
Hence $e,f$ are diverse in $G$.
By the symmetry we may therefore assume, by \ref{box} and \ref{superbox}, that
either $G$ is Box and $f=1\d 4$, or $G$ is Superbox and
$f=1\d 2$.
In the first case, by adding $P$ to $\eta (G)$ we obtain a \he
of Superbox contradicting the minimality of $|E(S)|$.
In the second case, by adding $P$ to
$\eta (G \setminus \{3\d 8,6\d 7 \})$ we obtain a
\he of Box contradicting the minimality of $|E(S)|$.

This proves our claim that $|E(S)| =1$.
From \ref{box3} and \ref{superbox5}, $H$ is arched.
This proves \ref{arched}.~\bbox

\section{The children of Drum}

The graph {\em Drum} is defined in Figure 10.

\begin{figure} [h!]
\centering
\includegraphics{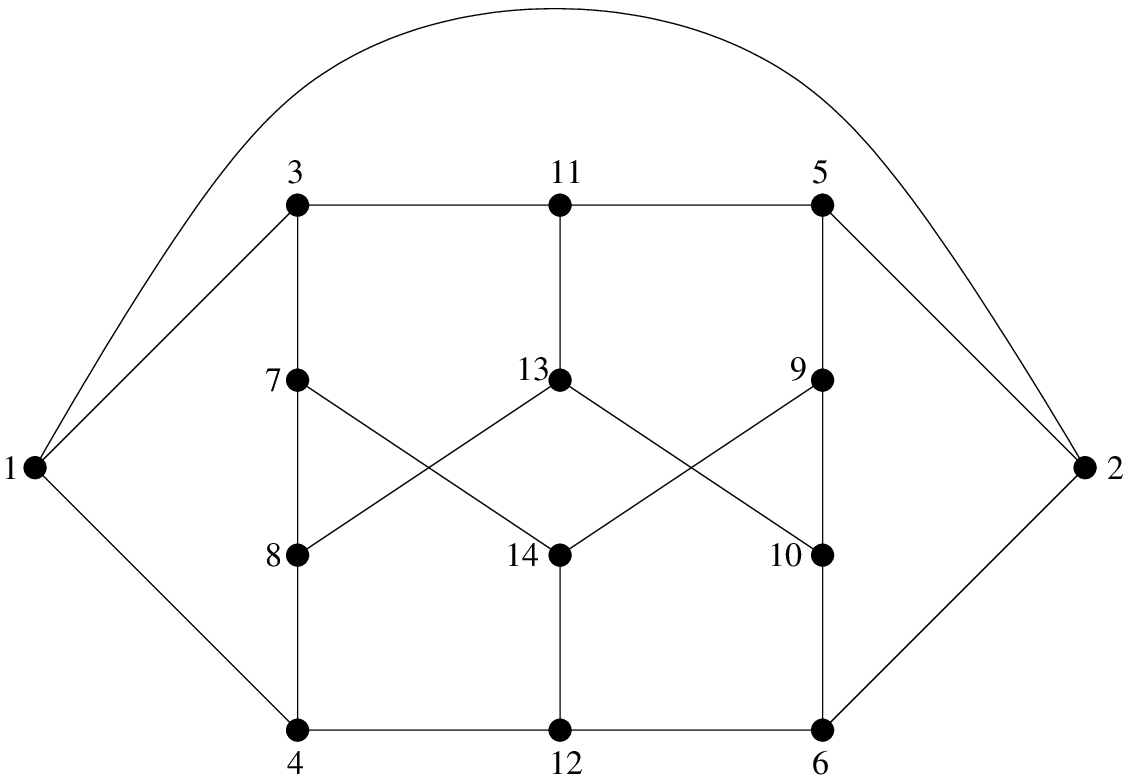}
\caption{Drum.}
\label{drumfig}
\end{figure}

\begin{thm}\label{drum}
Let $H$ be dodecahedrally-connected, and not isomorphic to Triplex.
Then $H$ is arched if and only it contains none of Petersen, Drum.
\end{thm}
\Proof
Since Drum contains Triplex (delete 9-10) ``only if'' follows
from \ref{arched}.
For ``if'', let $H$ be dodecahedrally-connected, not isomorphic to Triplex,
and not arched, and suppose that $H$ does not contain Petersen.
We must show that $H$ contains Drum.
By \ref{arched},
$H$ contains Triplex; and so by \ref{splitter}, since Triplex is not a
biladder, it follows that $H$ contains Triplex + $(e,f)$, where
$e,f$ are diverse edges of Triplex.
But for all such choices of $e,f$, Triplex + $(e,f)$ either contains
Petersen or is isomorphic to Drum.
This proves \ref{drum}.~\bbox

\begin{figure} [h!]
\centering
\includegraphics{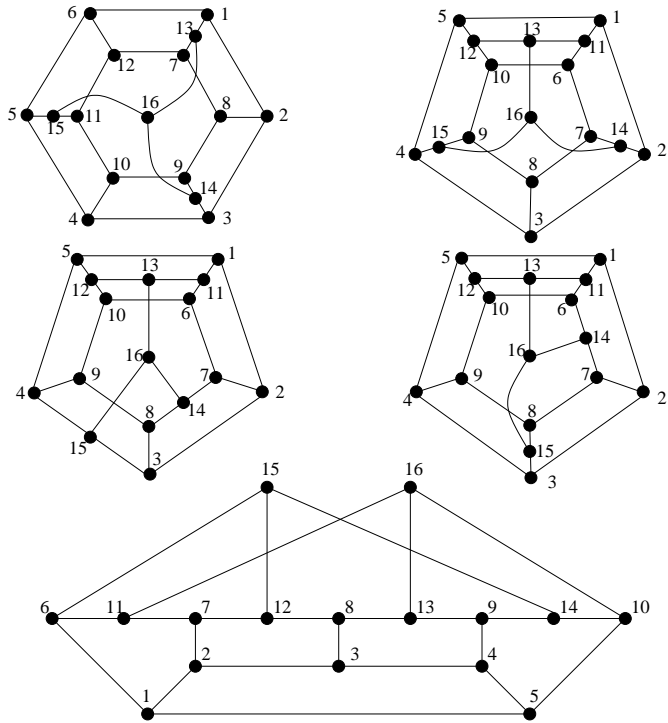}
\caption{Firstapex, Secondapex, Thirdapex, Fourthapex and Sailboat.}
\label{sailboatfig}
\end{figure}

In Figure 11 we define the graphs {\em Firstapex, Secondapex},
{\em Thirdapex, Fourthapex}, and {\em Sailboat}.
They all contain Drum.
We call the first four of them {\em Apex-selectors}.

\begin{thm}\label{selector}
Let $H$ be dodecahedrally-connected, and not isomorphic to
Triplex or Drum.
Then $H$ is arched if and only if it contains none of Petersen, an
Apex-selector, or Sailboat.
\end{thm}
\Proof
As in \ref{drum}, ``only if'' is easy, and for ``if'' we may assume that
$H$ contains Drum, by \ref{drum}.
By \ref{splitter} $H$ contains Drum + $(e,f)$ where $e,f$ are diverse edges
of Drum.
There are (up to isomorphism of Drum) 26 possibilities for
$\{e,f\}$; let $e=ab,f=cd$, and $G' =$ Drum + $(ab,cd)$.
If $(a,b,c,d)$ is one of 
$$(1,2,11,13), (1, 3, 8, 13), (3, 7, 5, 9), (3, 11, 9, 14), (7, 14, 11, 13),$$
$G$ is
isomorphic to Firstapex, Secondapex, Thirdapex, Fourthapex
and Sailboat respectively, and in all other cases $G$ contains
Petersen.
This proves \ref{selector}.~\bbox

Let us say $H$ is {\em doubly-apex} if it has two vertices
$u,v$ such that the graph obtained from $H$ by identifying $u$ and $v$
is planar.
Sailboat is doubly-apex (identify 15 and 16) but the Apex-selectors
are not, and Petersen is not.
The main result of this section is the following.

\begin{thm}\label{apexs}
Let $H$ be dodecahedrally-connected.
Then $H$ is either arched or doubly-apex if and only if it does not
contain Petersen or an Apex-selector.
\end{thm}

\ref{apexs} follows from the following.

\begin{thm}\label{sailboat}
Let $H$ be dodecahedrally-connected, and contain Sailboat but
not Petersen or any Apex-selector.
Then $H$ is doubly-apex.
\end{thm}

\noindent{\bf Proof of \ref{apexs} assuming \ref{sailboat}.}

``If'' is easy, and we omit it.
For ``only if'', let $H$ not contain Petersen or an Apex-selector.
If $H$ is isomorphic to Triplex or Drum it is doubly-apex as required.
Otherwise, by \ref{selector} either it is arched or it contains Sailboat;
and in the latter case by \ref{sailboat} it is doubly-apex.
This proves \ref{apexs}.~\bbox

It remains to prove \ref{sailboat}.
That will require several lemmas.
Let ${\cal C}$ be the set of the subgraphs of Sailboat induced on
the following vertex sets (which bound the regions when Sailboat is drawn 
in the plane with 15 and 16 identified):
\begin{eqnarray*}
&1, 2, 3, 4, 5;&\\&
 1, 2, 7, 11, 6;&\\&
2, 3, 8, 12, 7;&\\&
3, 4, 9, 13, 8;&\\&
4, 5, 10, 14, 9;&\\&
 15, 6, 1, 5, 10, 16; &\\&
15, 6, 11, 16;&\\&
16,11, 7, 12, 15; &\\&
15, 12, 8, 13, 16; &\\&
16, 13, 9, 14, 15; &\\&
15, 14, 10, 16.&
\end{eqnarray*}

Let Boat(1),\ldots,Boat(7) be Sailboat + $(ab,cd)$ where
respectively $(a,b,c,d)$ is 
$$(2, 7, 12, 15), (7, 12, 6, 15),
(1, 6, 11, 16), (2, 7, 11, 16), (6, 11, 12, 15),
(9, 14, 12, 15), (6, 15, 12, 15).$$

\begin{thm}\label{yards}
Let $G$ be Sailboat, and let $ab$ and $cd$ be edges of $G$
such that no member of ${\cal C}$ contains them both.
Then $G+(ab,cd)$ contains Petersen or an Apex-selector or
one of Boat(1),\ldots,Boat(7).
\end{thm}
\Proof
If $a=c$ then since no member of ${\cal C}$ contains $ab$ and $cd$
it follows that $a=15$ or $16$, and then $G+(ab,cd)$ is isomorphic to
Boat(7).
We assume therefore that $a,b \neq c,d$.

Up to the symmetry of Sailboat and exchanging $ab$ with $cd$,
there are 88 cases to be checked.
Let $G'=G+(ab,cd)$.
If $(a,b,c,d)$ is (1, 6, 11, 16) or (6, 15, 7, 11),
$G'$ is (isomorphic to) Boat(3).
If $(a,b,c,d)$ is (7,12,6,11) or (2,7,11,16),
$G'$ is Boat(4).
If $(a,b,c,d)$ is (2, 7, 12, 15) or (7, 11, 8, 12),
$G'$ is Boat(1).
If $(a,b,c,d)$ is (1, 6, 14, 15) or (6, 11, 12, 15), $G'$
is  Boat(5).
If $(a,b,c,d)$ is (7, 12, 6, 15) or (8, 12, 14, 15), $G'$
is Boat(2).
If $(a,b,c,d)$ is (9, 14, 12, 15) or (10, 14, 6, 15), $G'$
is Boat(6).
If $(a,b,c,d)$ = (2, 3, 12, 15), $G'$ contains Firstapex;
if $(a,b,c,d)$ = (1, 6, 10, 14), (3, 8, 12, 15) or (7, 11, 8, 13) it
contains Secondapex; if $(a,b,c,d)$ is one of 
$$(1, 5, 6, 11),
(1, 5, 14, 15), (1, 2, 6, 15), (1, 2, 11, 16), (2, 7, 6, 15),
(8, 12, 11, 16)$$
$G'$ contains Thirdapex;
and in the remaining 66 cases, $G'$ contains Petersen.
This proves \ref{yards}.~\bbox

\begin{thm}\label{boat}
Let $H$ be dodecahedrally-connected, and not contain Petersen or an
Apex-selector.
Then $H$ contains none of Boat(1),\ldots,Boat(7).
\end{thm}
\Proof
\\
\\
(1) {\em H does not contain Boat(1).}
\\
\\
\Subproof
Let ${\cal L}_{1}$ consist of Petersen and the four Apex-Selectors,
and let $C$ be the quadrangle of Boat(1).
Then every $A$-extension of Boat(1) is killed by
${\cal L}_{1}$, and ${\cal P}(C, {\cal L}_{1}) = \emptyset$,
so the claim follows from \ref{Aextn}.
This proves (1).
\\
\\
(2) {\em H does not contain Boat(2).}
\\
\\
\Subproof
Let $C$ be the quadrangle of Boat(2).
Then every $A$-extension of Boat(2) is killed by ${\cal L}_{1}$, and
\[
{\cal P}(C,{\cal L}_{1}) = \{(17\d 18,6\d 11),
(17\d 18,7\d 11) \}^{\ast}.
\]
The result follows from \ref{extn}.
This proves (2).
\\
\\
(3) {\em $H$ does not contain Boat(3) or Boat(4).}
\\
\\
\Subproof
Let $G$ be Boat(3) or Boat(4), and
${\cal L}_{3} = {\cal L}_{1} \cup \{Boat(2)\}$.
Let $C$ be the quadrangle of $G$.
Then every $A$-extension of $G$ is killed by ${\cal L}_{3}$, and
${\cal P}(C,{\cal L}_{3})= \emptyset$, so the result follows
from (2) and \ref{Aextn}. This proves (3).
\\
\\
(4) {\em $H$ does not contain Boat(5) or Boat(6).}
\\
\\
\Subproof
Let $G$ be Boat(5) or Boat(6), and let
\[
{\cal L}_{4} = {\cal L}_{3} \cup \{Boat(3), Boat(4)\} .
\]
Let $C$ be the quadrangle of $G$.
Then every $A$-extension of $G$ is killed by ${\cal L}_{4}$,
and ${\cal P}(C,{\cal L}_{4}) = \emptyset$, so the result
follows from (2), (3) and \ref{Aextn}. This proves (4).
\\
\\
(5) {\em H does not contain Boat(7).}
\\
\\
\Subproof
Let $G$ be Boat(7), and let $C$ be its circuit of length 3.
Let $X=V(C)$.
Suppose that there is a \he of $G$ in $H$; then by \ref{augment},
there is a $X$-augmenting sequence
$(e_{1},f_{1}),\ldots,(e_{n},f_{n})$ of $G$ such that $H$ contains
$G+(e_{1},f_{1})+\ldots+ (e_{n},f_{n})$.
From the definition of ``$X$-augmentation'' it follows that
$n=1$ since $|E(C)| =3$; and so $H$ contains
$G(e_{1},f_{1})$ for some $e_{1} \in E(C)$ and
$f_{1} \in E(G \setminus X)$.
But for all such $e_{1},f_{1}, G+(e_{1},f_{1})$ contains a member of
${\cal L}_{1}$ or one of Boat(2), Boat(5), Boat(6), a contradiction
by (2) and (4). This proves (5).

\bigskip

From (1)--(5), this proves \ref{boat}.~\bbox

\bigskip
\noindent{\bf Proof of \ref{sailboat}.}

Let $H$ be dodecahedrally-connected and not contain Petersen or an
Apex-selector.
Let $\eta$ be a \he of $G$ in $H$, where $G$ is Sailboat.
Let $V(F) = \{15, 16\}$ and
$E(F) = \emptyset$; and let $\eta_F$ be the restriction of $\eta$ to $F$. 
Let ${\cal C}$ be as before.
Then $(G,F, {\cal C})$ is a framework, and we claim that
(E1)--(E7) hold.
By \ref{boat} $H$ contains none  of Boat(1),..., Boat(7), so by \ref{yards}
(E2) holds.
All the others are clear except for (E6), and for (E6) we need only
consider cross-extensions of $G$ on some of the paths in ${\cal C}$,
namely the ones with vertex sets 
$$\{15, 6, 1, 5, 10, 16\},\{16, 11, 7, 12, 15\}, \{15, 12, 8, 13, 16\}$$
(and two more, that
from symmetry we need not consider).
We need to examine
\begin{eqnarray*}
&G+(6\d 15, 10\d 16) + (1\d 6,16\d 18)&\\&
G+(6\d15, 10\d 16) + (6\d 17,16\d 18) &\\&
G+(6\d 15,5\d 10) + (6\d 17,10\d 18) &\\&
G+(6\d 15,5\d 10 ) + (1\d 6,10\d 16) &\\&
G+(11\d 16,12\d 15) + (11\d 17,15\d 18 ) &\\&
G+(12\d 15,13\d 16 ) + (12\d 17,16\d 18) ;&
\end{eqnarray*}
they contain Thirdapex, Boat(3), Boat(3), Petersen, Boat(3) and Boat(3)
respectively.
Hence (E6) holds, and from \ref{mainthm}, this proves \ref{sailboat}.~\bbox

\section{Dodecahedrally connected non-apex graphs}

The graphs {\em Diamond, Bigdrum} and {\em Concertina} are
defined in Figure 12.

\begin{figure} [h!]
\centering
\includegraphics{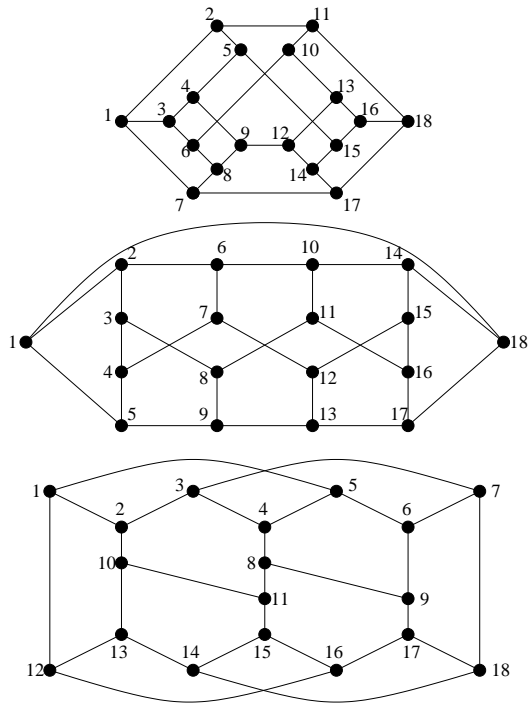}
\caption{Diamond, Bigdrum and Concertina.}
\label{diamondfig}
\end{figure}
In this section we prove the following.

\begin{thm}\label{apex}
Let $H$ be dodecahedrally-connected.
Then $H$ is apex if and only if it contains none of Petersen, Jaws,
Starfish, Diamond, Concertina, Bigdrum.
\end{thm}

Let Square(1) be Secondapex + $(14\d 16,11\d 13)$.
Let Square(2),..., Square(5) be Fourthapex + $(ab, cd)$ where
$(a,b,c,d)$ is 
$$(1, 5, 10, 12), (1, 11, 6, 10), (6, 14, 13, 16), (12, 13, 15, 16)$$
respectively.
Let Square(6) and Square(7) be Thirdapex + $(ab,cd)$ where
$(a,b,c,d)$ is $(3, 15, 14, 16)$ and $(2, 3, 8, 9)$ respectively.

\begin{thm}\label{square}
Let $H$ be dodecahedrally-connected, and not contain any of Petersen,
Jaws, Starfish, Diamond, Concertina, Bigdrum.
Then it contains none of Square(1),..., Square(7).
\end{thm}
\Proof
\\
\\
(1) {\em H does not contain Square(1).}
\\
\\
\Subproof
Let $G$ be Square(1), let $C$ be the quadrangle of $G$, and let
\[
{\cal L}_{1} = {\rm \{Petersen, Jaws, Starfish, Diamond, Concertina, Bigdrum}\}.
\]
Every $A$-extension of $G$ is killed by ${\cal L}_{1}$
(indeed, by $\{$Petersen, Jaws, Starfish$\}$), and
\[
{\cal P}(C,{\cal L}_{1} ) = \{ (13\d 18,5\d 12),
(13\d 18,10\d 12 ), (13\d 18,1\d 11),
(13\d 18, 6\d 11 ) \}^{\ast}.
\]
(Note that $G+(13\d 18,1\d 5)$ is isomorphic to Jaws, and
$G+(16\d 17,3\d 8)$ to Starfish.)
Then we verify the hypotheses of \ref{extn}; and find that all the
various extensions listed in \ref{extn} contain Petersen, except for the
$B$-extensions
\begin{eqnarray*}
&G+(13\d 18,12\d 5) + (12\d 20,4\d 15) &\\&
G+(13\d 18,12\d 5 ) + (12\d 20,11\d 18)&\\&
G+(13\d 18,12\d 5 ) + (12\d 20,15\d 16).&
\end{eqnarray*}
(which contain Jaws, Diamond, and Concertina respectively) and the $C$-extension
\[
G+(13\d 18,12\d 5 ) +(19\d 20,1\d 11) 
\]
(which contains Jaws), and isomorphic extensions.
Hence, from \ref{extn}, this proves (1).

\bigskip

Now let
\begin{center}
${\cal L}_{2}$ = $\{$Petersen, Square(1), Diamond, Concertina, Bigdrum$\}$
\end{center}
(Jaws and Starfish are no longer necessary, since they both contain
Square(1).)
\\
\\
(2) {\em $H$ does not contain Square(2).}
\\
\\
\Subproof
We apply \ref{Aextn} to the quadrangle $C$ of Square(2), with
${\cal L} = {\cal L}_{2}$.
All $A$-extensions are killed by ${\cal L}_{2}$, and
${\cal P}(C,{\cal L}_{2} ) = \emptyset$, so the result follows from \ref{Aextn}.
This proves (2).
\\
\\
(3) {\em $H$ does not contain Square(3).}
\\
\\
\Subproof
Let $C$ be the quadrangle of $G$ = Square(3); we apply \ref{extn},
with ${\cal L}= {\cal L}_{2}$.
All $A$-extensions are killed by ${\cal L}_{2}$, and
\[
{\cal P} (C, {\cal L}_{2}) = \{(6\d 11,13\d 16),
(6\d 11,14\d 16)\}^{\ast}.
\]
We verify the hypotheses of \ref{extn}. This proves (3).
\\
\\
(4) {\em $H$ does not contain Square(4).}
\\
\\
\Subproof
Now let ${\cal L}_{4} = {\cal L}_{2} \cup$
\{Square(2), Square(3)\}.
The result follows from \ref{Aextn}, applied to the quadrangle of Square(4)
and ${\cal L}_{4}$, using (2) and (3). This proves (4).
\\
\\
(5) {\em $H$ does not contain Square(5).}
\\
\\
\Subproof
Let ${\cal L}_{5} = {\cal L}_{4} \cup$ \{Square(4)\}, and
$C$ the quadrangle of $G$ = Square(5).
Then all $A$-extensions are killed by ${\cal L}_{5}$, and
\[
{\cal P}(C, {\cal L}_{5}) = \{(13\d 17,6\d 11)\}^{\ast};
\]
and we verify the hypotheses of \ref{extn} to prove (5).
\\
\\
(6) {\em $H$ does not contain Square(6).}
\\
\\
\Subproof
Let ${\cal L}_{6} = {\cal L}_{5} \cup$ \{Square(5)\}, and $C,G$ as usual.
All $A$-extensions are killed by ${\cal L}_{6}$, and
\[
{\cal P}(C,{\cal L}_{6}) = \{(17\d 18,3\d 8),
(17\d 18,8\d 14)\}^{\ast};
\]
and again the result follows from \ref{extn}. This proves (6).
\\
\\
(7) {\em $H$ does not contain Square(7).}
\\
\\
\Subproof
Let ${\cal L}_{7} = {\cal L}_{6} \cup$ \{Square(6)\}, and $C,G$ as usual.
Then all $A$-extensions are killed by ${\cal L}_{7}$, and
${\cal P}(C, {\cal L}_{7}) = \emptyset$, so (7) follows from \ref{Aextn}.

\bigskip
From (1)--(7), this proves \ref{square}.~\bbox

The graph {\em Extrapex} is defined in Figure 13.
\begin{figure} [h!]
\centering
\includegraphics{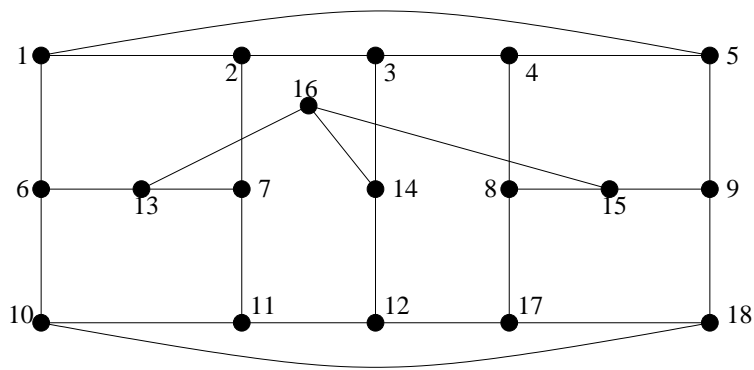}
\caption{Extrapex.}
\label{extrapexfig}
\end{figure}
We say that $G$ is an {\em Apex-forcer} if either it is an Apex-selector
or it is Extrapex.
By the {\em Non-apex family} we mean
\begin{center}
\{Petersen, Diamond, Concertina, Bigdrum, Square(1),..., Square(7)\}.
\end{center}

\begin{thm}\label{forcer}
Let $G$ be an Apex-forcer.
Let ${\cal C}$ be the set of circuits that bound regions in the planar
drawing of $G \setminus 16$.
If $ab$ and $cd$ are edges of $G$ with $a,b \neq c,d$, and no member of
${\cal C}$ contains them both, then either $G+(ab,cd)$
contains a member of the Non-apex family, or one of $a,b,c,d$ is $16$ and the other three belong to some member of ${\cal C}$.
\end{thm}

We leave the proof to the reader (the details are in the Appendix~\cite{RSTappendix}).
If $G$ is an Apex-forcer, and $\eta$ is a \he of $G$ in $H$,
we define the {\em spine} of $\eta$ to be
$\eta (13\d 16) \cup \eta (14\d 16) \cup \eta (15\d 16)$.

\begin{thm}\label{spine}
Let $H$ be cubic and cyclically four-connected, and contain no member
of the Non-apex family.
Let $H$ contain some Apex-forcer.
Then there is a \he $\eta$ of some Apex-forcer in $H$ such that
its spine has only three edges.
\end{thm}
\Proof
Choose an Apex-forcer $G$ and a \he $\eta$ of $G$ in $H$, such that
its spine is minimal.
Suppose its spine has more than three edges; then since $H$ is
cyclically four-connected, there is an $\eta$-path $P$ with one end in
$\eta (e)$ and the other in $\eta (f)$, where $f$ is one of
$13\d 16,14\d 16,15\d 16$ and $e$ is not incident
with $16$.
If $e$ and $f$ have a common end then by rerouting $e$ along $P$
we obtain a new \he with smaller spine, a contradiction.
Similarly, it follows that no edge of $G\setminus 16$
joins an end of $e$ to an end of $f$.
Let $\mathcal{C}$ be as in \ref{forcer}.
By \ref{forcer} there exists $C \in {\cal C}$ such that
$e \in E(C)$ and $f$ has an end in $V(C)$.
Let $e = ab$ and let $f$ be incident with $c,16$.
Now we must examine cases.

If $G$ is Firstapex, we may assume that $(a,b,c)=(2,8,13)$ from the
symmetry.
Then $\eta(G \setminus 6\d 12) \cup P$ yields a \he
of Secondapex with smaller spine, a contradiction. (We apologize for this awkward notation;
by $G \setminus 6\d 12$ we mean the graph obtained from $G$ by deleting the edge 6-12. We use
the same notation below.)

If $G$ is Secondapex, there are three possibilities for
$(a,b,c): (1,5,13)$ (when $\eta (G \setminus 6\d 10) \cup P$
yields a \he of Firstapex), (1, 11, 14) (when
$\eta (G \setminus 1\d 5 ) \cup P$ yields a \he of Fourthapex),
and (3, 8, 14) (when $\eta(G) \cup P$ yields a \he of Extrapex),
in each case contradicting the minimality of the spine.
If $G$ is Thirdapex, the possibilities for $(a,b,c)$ are:
(1, 5, 13) or (2, 3, 14) (when $\eta(G \setminus 8\d 9) \cup P$
yields a \he of Fourthapex), (6, 10, 14) (when
$\eta (G \setminus 1\d 11) \cup P$ yields a \he of Thirdapex),
and (9, 10, 14) (when $\eta(G \setminus 2\d 7 ) \cup P$
yields a \he of Firstapex), in each case a contradiction.

If $G$ is Fourthapex, the possibilities are: (1, 5, 13) (when
$\eta (G \setminus 4\d 9) \cup P$ yields a \he of Thirdapex),
(6, 10, 13) (when $\eta (G) \cup P$ yields a \he of Extrapex),
(1, 2, 14) (when $\eta(G \setminus 4\d 9) \cup P$ yields a
\he of Secondapex), and (1, 11, 14) (when 
$\eta (G \setminus 10\d 12) \cup P$ yields a \he of Thirdapex),
in each case a contradiction. (We haved used a symmetry of Fourthapex 
not evident from the drawing, exchanging 13 with 15 and 1 with 9.)

If $G$ is Extrapex, the possibilities are: (1, 2, 13)
(when $\eta (G \setminus \{7\d 13,1\d 6\}) \cup P$
yields a \he of Secondapex) and (2, 7, 14) (when
$\eta ( G \setminus \{2\d 3,10\d 11\}) \cup P$
yields a \he of Thirdapex), in each case a contradiction.

Hence the spine has only three edges.
This proves \ref{spine}.~\bbox

\bigskip
\noindent{\bf Proof of \ref{apex}.}

``Only if'' is easy, and we omit it.
For ``if'', let $H$ be dodecahedrally-connected, and not contain
any of Petersen, Jaws, Starfish, Diamond, Concertina, Bigdrum.
By \ref{square} it contains none of Square(1),\ldots, Square(7).
We may assume that $H$ is not arched or doubly-apex, for such graphs
are apex; and so by \ref{apexs} $H$ contains an Apex-selector.
By \ref{spine}, there is a \he $\eta$ of some Apex-forcer $G$ in $H$
such that its spine has only three edges.
Let $F$ be the subgraph of $G$ induced on $\{13, 14, 15, 16\}$, and let $\eta_F$
be the restriction of $\eta$ to $F$.
Let ${\cal C}$ be as in \ref{forcer}; then $(G,F, {\cal C})$
is a framework, and $H, \eta_{F}$ satisfy (E1).
We claim they satisfy (E2)--(E7).
(E2) follows from \ref{forcer}, and (E3), (E7) are vacuously true.
For (E4), (E5) and (E6) a large amount of case-checking
is required, for $G$ = Firstapex, Secondapex, Thirdapex,
Fourthapex and Extrapex, separately.
(In the case-checking we use that $H$ contains none of 
Petersen, Jaws, Starfish, Diamond, Concertina, Bigdrum, and we could also use 
that it contains none of Square(1)--Square(7). In fact we find that
we don't need to use all of the latter; we just
need that $H$ does not contain Square(2).)
The details are in the Appendix~\cite{RSTappendix}.
From \ref{mainthm}, this proves \ref{apex}.~\bbox

\section{Die-connected non-apex graphs}

\begin{figure} [h!]
\centering
\includegraphics{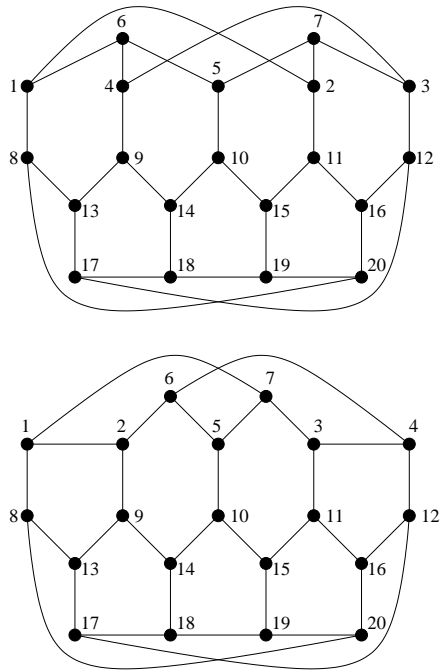}
\caption{Antilog and Log.}
\label{logfig}
\end{figure}
Our next real objective in this paper is modify \ref{apex} to find all
the cubic graphs $G$ minimal with the properties that they are
non-apex and dodecahedrally-connected, and $| \delta(X) | \geq 6$ for
all $X \subseteq V(G)$ with $|X|,|V(G) \setminus X| \geq 7$.
(There are only three such graphs, namely Petersen, Jaws and Starfish,
as we shall see in the next section.)
Diamond, Concertina and Bigdrum all have subsets $X$ with
$|\delta(X) | =5$ and $|X|,|V(G) \setminus X| \geq 9$,
so they are rather far from having the property we require; and
a convenient half-way stage is afforded by ``die-connectivity''.
We recall that a graph $G$ is {\em die-connected} if it is dodecahedrally-connected 
(and hence cubic and cyclically five-connected) and
$|\delta(X) | \geq 6$ for all $X \subseteq V(G)$ with
$|X|,|V(G) \setminus X| \geq 9$.
In this section we find all minimal graphs that are non-apex and
die-connected.
\begin{figure} [h!]
\centering
\includegraphics{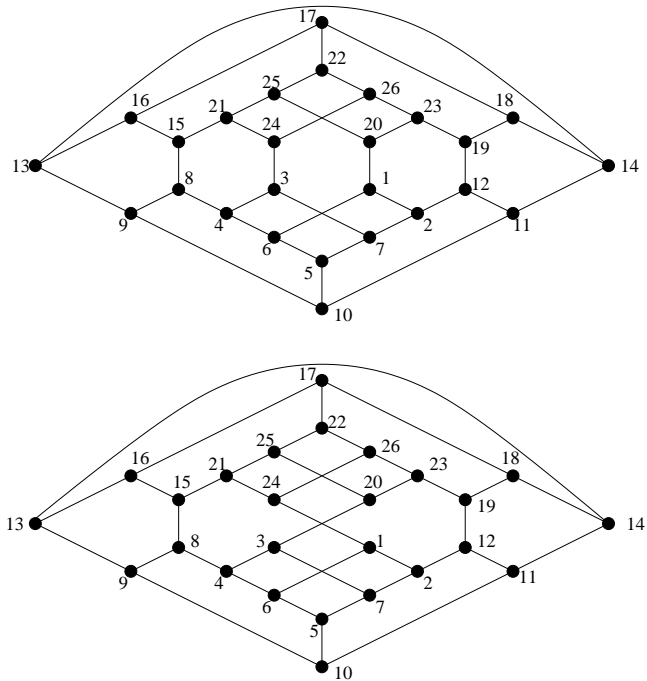}
\caption{Dice(1) and Dice(3).}
\label{dice1fig}
\end{figure}
\begin{figure} [h!]
\centering
\includegraphics{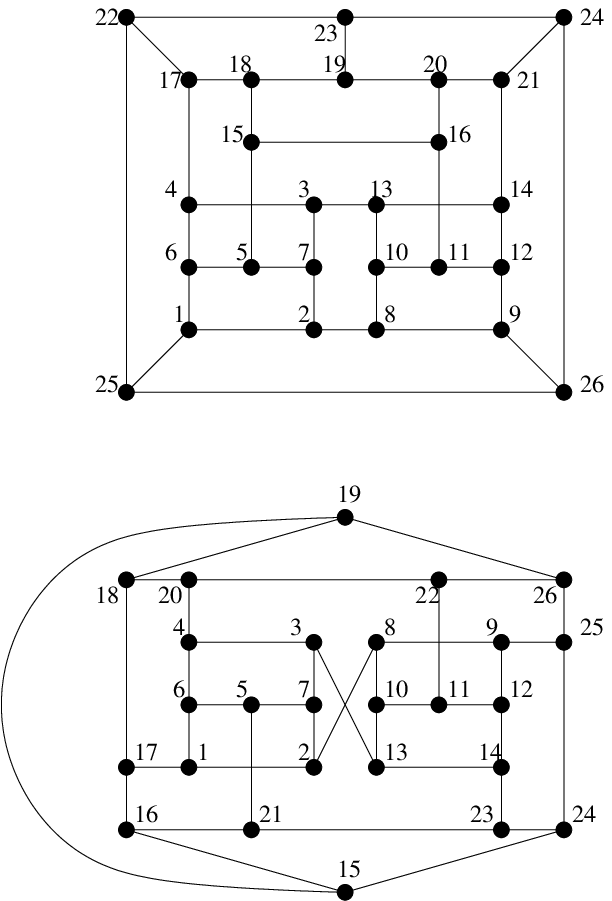}
\caption{Dice(2) and Dice(4).}
\label{dice3fig}
\end{figure}
The graphs Log, Antilog, and Dice(1),..., Dice(4) are defined in Figures 14,
15 and 16. We shall show the following.

\begin{thm}\label{dieapex}
Let $H$ be die-connected.
Then $H$ is apex if and only if H contains none of Petersen,
Jaws, Starfish, Log, Antilog, Dice(1), Dice(2), Dice(3),
Dice(4).
\end{thm}

We begin with the following.

\begin{thm}\label{diamond}
Any die-connected graph that contains Diamond also contains one of
Petersen, Antilog, Dice(4).
\end{thm}
\Proof
Let $H$ be die-connected, and contain no member of
${\cal L} =$ \{Petersen, Antilog, Dice(4)\}.
We claim first that
\\
\\
(1) {\em H does not contain Diamond $+ (1\d 2,10\d 11)$.}
\\
\\
\Subproof
Let $C$ be the quadrangle of $G =$ Diamond +$(1\d 2,10\d 11)$.
Then all $A$-extensions are killed by ${\cal L}$, and
\[
{\cal P} (C, {\cal L} ) = \{(2\d 19,4\d 5),
(11\d 20,10\d 13)\}^{\ast}.
\]
We verify the hypotheses of \ref{extn}
(the $E$-extension is isomorphic to Dice(4)). This proves (1).

\bigskip

Now let ${\cal L'} =$ \{Petersen, Antilog, Diamond +
$(1\d 2,10\d 11)$\}, and $X = \{1,\ldots,9\}$.
\\
\\
(2) {\em Every $X$-augmentation of Diamond contains a member of $ {\cal L}'$.}
\\
\\
\Subproof
Let $(e_{1},f_{1}),\ldots,(e_{n},f_{n})$ be an $X$-augmenting sequence,
and suppose the corresponding $X$-augmentation contains no member
of ${\cal L}'$.
In particular, Diamond + $(e_{1},f_{1})$ contains no member of
${\cal L}'$, and so (by checking all possibilities) it follows
that $f_{1}$ is 6-10 and $e_{1}$ is one of
1-2, 1-7, 4-9.
In particular, $n \geq 2$.
Since $f_{1}=6\d 10$ it follows that
$e_{2} = 6\d 20$.
If $e_{1}$ is 1-7 or 4-9 there is
no possibility for $f_{2}$.
Thus $e_{1}$ is 1-2, and then $f_{2}$ is
9-12, and $n \geq 3$,
and $e_{3}$ is 9-22.
Again by checking cases it follows that $f_{3}$ is 7-17, and hence
$n \geq 4$ and $e_{4}$ is 7-24; and there is no possibility
for $f_{4}$, a contradiction.
This proves (2).

\bigskip

From (1), (2) and \ref{augment}, the result follows since $H$ is die-connected.
This proves \ref{diamond}.~\bbox

\begin{thm}\label{bigdrum}
Every die-connected graph that contains Bigdrum also contains
one of Petersen, Diamond or Dice$(2)$.
\end{thm}
\Proof
Let $H$ be die-connected, and contain no member of
${\cal L} = $ \{Petersen, Diamond, Dice(2)\}.
We claim first
\\
\\
(1) {\em $H$ does not contain Bigdrum $+ (3\d 8,10\d 11)$.}
\\
\\
\Subproof
Let $G$ = Bigdrum + $(3\d 8,10\d 11)$, and let $C$ be
the quadrangle of $G$.
Then all $A$-extensions are killed by ${\cal L}$, and
\[
{\cal P} (C, {\cal L}) = \{ (8\d 11,9\d 13),
(19\d 20,10\d 14)\}^{\ast}.
\]
The result follows from \ref{extn} by checking all the various
extensions (in particular,
\[
G+(8\d 19, 5\d 9) + (11\d 20, 10\d 14 ) +
( 8\d 21, 20\d 23)
\]
is isomorphic to Dice(2)).
This proves (1).

\bigskip

Now let ${\cal L}'$ = \{Petersen, Diamond,
Bigdrum + $(3\d 8,10\d 11)\}$ and
$X = \{1,\ldots,9\}$.
We claim that
\\
\\
(2) {\em Every $X$-augmentation of Bigdrum contains a member of ${\cal L}'$.}
\\
\\
\Subproof
Let $(e_{1},f_{1}),\ldots,(e_{n},f_{n})$ be an $X$-augmenting sequence,
such that the corresponding  $X$-augmentation contains no member
of ${\cal L}'$.
Then by checking cases it follows that $(e_{1},f_{1})$ is one of
$(3\d 8,6\d 10), (4\d 7,9\d 13)$, and by the
symmetry we may assume the first.
Then $n \geq 2$, and $e_{2}$ is 6-20; and there is no possibility
for $f_{2}$, a contradiction.
This proves (2).

\bigskip

From (1), (2) and \ref{augment}, this proves \ref{bigdrum}.~\bbox

\begin{thm}\label{concertina}
Any die-connected graph that contains  Concertina also contains one of
Petersen, Log, Diamond, Bigdrum, Dice$(1)$, Dice$(3)$.
\end{thm}
\Proof
Let $H$ be a die-connected graph that contains no member of
${\cal L}$ = \{Petersen, Log, Diamond, Bigdrum, Dice(1), Dice(3)\}.
Let Conc(1), Conc(2), Conc(3) be Concertina + $(e,f)$ where
$(e,f)$ is $(4\d 8,10\d 11)$,
$(6\d 7,17\d 18)$, $(8\d 9,16\d 17)$; and let
Conc(4) be Concertina + $(2\d 3,8\d 11) + (8\d 20,16\d 17)$.
\\
\\
(1) {\em H does not contain Conc(1).}
\\
\\
\Subproof
Let $C$ be the quadrangle of $G$ = Conc(1).
All $A$-extensions are killed by ${\cal L}$, and
\[
{\cal P} (C, {\cal L}) = \{(8\d 11,9\d 17),
(19\d 20,2\d 10)\}^{\ast};
\]
and the result follows by verifying the other hypotheses of \ref{extn}.
(The $E$-extension is isomorphic to Dice(1).) This proves (1).

\bigskip

Let Conc(21) be Conc(2) + $(7\d 19,1\d 5)$, let
Conc(211) be Conc(21) + $(1\d 2,3\d 4)$, and let Conc(212)
be Conc(21) + $(1\d 2,3\d 7)$. 
\\
\\
(2) {\em $H$ does not contain  Conc$(211)$ or Conc$(212)$.}
\\
\\
\Subproof
Let $G$ = Conc(211) and let $C$ be its quadrangle.
Then all $A$-extensions are killed by ${\cal L}$, and
\[
{\cal P} (C, {\cal L}) = \{(2\d 23,1\d 12)\}^{\ast},
\]
and the result for Conc(211) follows by verifying the other hypotheses
of \ref{extn}.

Now let $G$ = Conc(212) and let $C$ be its quadrangle.
Again all $A$-extensions are killed by ${\cal L}$, and again
\[
{\cal P} (C,{\cal L}) = \{(2\d 23,1\d 12)\}^{\ast}
\]
and again the result follows from \ref{extn}.
(Conc(212) + $(3\d 24,1\d 22)$ is isomorphic to Dice(3).)
This proves (2).
\\
\\
(3) {\em $H$ does not contain Conc(21).}
\\
\\
\Subproof
Let ${\cal L}_{1} = {\cal L} \cup$ \{Conc(211), Conc(212)\}.
Let $X = \{1, 2, 10, 11, 12, 13, 14, 15, 16\}$; we claim that
every $X$-augmentation of Conc(21) contains a member of ${\cal L}_{1}$.
For suppose not, and let the corresponding sequence be
$(e_{1},f_{1}),\ldots,(e_{n},f_{n})$.
By checking cases, $e_{1}$ is 12-16 and
$f_{1}$ is 14-18; and so $n \geq 2$, and
$e_{2}$ is 14-20, and there is no possibility for $f_{2}$.
Hence (3) follows from \ref{augment} and (2).
\\
\\
(4) {\em $H$ does not contain Conc(2).}
\\
\\
\Subproof
Let ${\cal L}_{2} = {\cal L} \cup$ \{Conc(21)\},
$G$ = Conc(2), and $C$ the quadrangle of $G$.
Then all $A$-extensions are killed by ${\cal L}_{2}$, and
\[
{\cal P} (C,{\cal L}_{2} ) = \{(19\d 20,6\d 9),
(19\d 20, 9\d 17)\}^{\ast}
\]
and the result follows by verifying the hypotheses of \ref{extn}.
This proves (4).
\\
\\
(5) {\em H does not contain Conc(3).}
\\
\\
\Subproof
Let ${\cal L}_{3} = {\cal L} \cup$ \{Conc(2)\},
$G$ = Conc(3), and $C$ the quadrangle of $G$.
Then all $A$-extensions are killed by ${\cal L}_{3}$, and
\[
{\cal P}(C, {\cal L}_{3} ) = \{(9\d 19,4\d 8) \}^{\ast},
\]
and the result follows by verifying the hypotheses of \ref{extn}.
This proves (5).
\\
\\
(6) {\em $H$ does not contain Conc(4).}
\\
\\
\Subproof
Let ${\cal L}_{4} = {\cal L} \cup$ \{Conc(2), Conc(3)\}, and
$X$ = \{3, 4, 5, 6, 7, 8, 9, 17, 18\}.
We claim that every $X$-augmentation of $G$ = Conc(4) contains
a member of ${\cal L}_{4}$.
Suppose not, and let the
corresponding sequence be $(e_{1},f_{1}),\ldots,(e_{n},f_{n})$.
By checking cases, $e_{1}$ is 3-7 and
$f_{1}$ is 1-5; so $n \geq 2$, and $e_{2}$ is
5-24, and there is no possibility for $f_{2}$,
a contradiction.
Hence (6) follows from \ref{augment}.

\bigskip

Let ${\cal L}_{5} = {\cal L}_{4} \cup$ \{Conc(1), Conc(4)\}, and
$X$ = \{1,\ldots,9\}.
We claim that every $X$-augmentation of $G$ = Concertina contains
a member of ${\cal L}$.
Suppose not, and let $(e_{1},f_{1}),\ldots,(e_{n},f_{n})$ be the
corresponding sequence.
By checking cases $(e_{1},f_{1})$ is one of
$(2\d 3,8\d 11)$, $(4\d 8,2\d 10)$;
so $n \geq 2$, and in either case there is no possibility for $f_{2}$.
Hence the result follows from (1), (4), (5), (6) and \ref{augment}.
This proves \ref{concertina}.~\bbox

\bigskip

\noindent{\bf Proof of \ref{dieapex}.}
``Only if'' is easy, and we omit it.
For ``if'', let $H$ contain none of the given graphs.
By \ref{diamond}, \ref{bigdrum}, \ref{concertina} it contains none of Diamond,
Bigdrum, Concertina; and so by \ref{apex} it is apex.
This proves \ref{dieapex}.~\bbox



\section{Theta-connected non-apex graphs}

We recall that $G$ is {\em theta-connected} if it is cubic
and cyclically five-connected, and $|\delta (X) | \geq 6$ for all
$X \subseteq V(G)$ with $|X|, |V(G)  \setminus  X| \geq 7$
(and hence it is dodecahedrally-connected).
None of the graphs of Figures 14--16  are theta-connected, and our next
objective is to make a version of \ref{dieapex} for theta-connected graphs.
It becomes much simpler:

\begin{thm}\label{thetamain}
Let $H$ be theta-connected.
Then $H$ is apex if and only if it contains none of Petersen,
Jaws and Starfish.
\end{thm}

For the proof we use \ref{domino} below.
A {\em domino} in a cubic graph $H$ is a subgraph $D$ with
$|V(D) | = 7$, consisting of the union of three paths
$P_{1},P_{2},P_{3}$ of lengths two, three and three respectively, which have
common ends and otherwise are disjoint.
The middle vertex of $P_{1}$ is called the {\em centre} of the domino,
and the other four vertices of degree two are its {\em corners};
an {\em attachment sequence}
\begin{figure} [h!]
\centering
\includegraphics{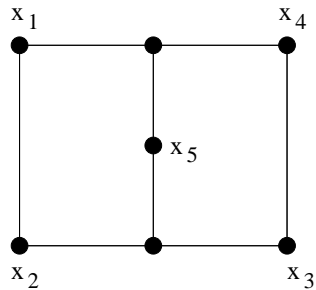}
\caption{A domino.}
\label{dominofig}
\end{figure}
is some sequence $(x_{1},\ldots,x_{5})$ where $x_{1},\ldots,x_{4}$
are the corners, $x_{5}$ is the centre,
$x_{1}x_{2}$ is an edge, and $x_{2},x_{3}$ have a common neighbour.
(See Figure 17.)

A domino $D$ in $G$ with attachment sequence $(x_{1},\ldots,x_{5})$
is said to be {\em crossed} if
\begin{itemize}
\item
there are two disjoint connected subgraphs
$P,Q$ of $G$, both edge-disjoint from $D$, with
$V(P \cap D) = \{x_{1},x_{3}\}$ and
$V(Q \cap D) = \{x_{2},x_{4},x_{5}\}$, and
\item
there are two disjoint connected subgraphs $P,Q$ of $G$, both
edge-disjoint from $D$, with
$V(P \cap D) = \{x_{1},x_{3},x_{5}\}$ and
$V(Q \cap D) = \{x_{2},x_{4}\}$.

\end{itemize}

\begin{thm}\label{domino}
Let $D$ be a crossed domino with attachment sequence
$x_{1},\ldots,x_{5}$, in a cyclically five-connected cubic graph $G$ with
$|V(G)| \geq 14$.
Let $x_{5}$ be incident with $g \not\in E(D)$.
Let $H$ be a cubic graph, cyclically five-connected, 
and let $\eta$ be a \he of $G$ in $H$.
Then either
\begin{itemize}
\item
there exists $X \subseteq V(H)$ with $| \delta_{H} (X) | =5$, such that for
all $v \in V(G)$, $\eta (v) \in X$ if and only if $v \in V(D)$, or
\item
$H$ contains Petersen, or
\item
for some $e \in E(D)$ and $f \in E(G \setminus V(D)) $
there is a \he $\eta '$ of $G+(e,f)$ in $H$, or
\item
for some $e \in \{x_{1}x_{2},x_{3}x_{4}\}$, and for some
$f \in E(G \setminus V(D)) $ such that $f, g$ are diverse in $G$, 
there is a \he $\eta'$ of
\[
G+(e,g) + (yx_{5},f)
\]
in $H$, where $x, y$ are the new vertices of $G+(e,g)$.
\end{itemize}
\end{thm}
\Proof
Let $X=V(D)$.
We assume that (i) and (ii) are false.
Since $|V(G)| \geq 14$ and $|\delta_{G} (X) | =5$, and since (i) is false,
it follows from \ref{augment} that there is an $X$-augmentation $G'$ of $G$, and
a \he $\eta'$ of $G'$ in $G$.
Let $(e_{1},f_{1}),\ldots,(e_{n},f_{n})$ be the corresponding sequence.
If $n=1$ then (iii) is true, so we assume that $n \geq 2$.
For $1 \leq i \leq 5$, let $x_{i}$ be adjacent in $G$ to
$y_{i} \in V(G) \setminus V(D)$.
Let the neighbours of $x_{5}$ in $G$ be $y_{5},x_{6},x_{7}$, where
$x_{6}$ is adjacent to $x_{1}$.
Let $G_{1} = G+(e_{1},f_{1})$ with new vertices $s_{1},t_{1}$, and
let $D_{1}$ be the subgraph of $G_{1}$ induced on
$V(D) \cup \{s_{1},t_{1}\}$.

Suppose first that $f_{1} = x_{1}y_{1}$.
Then since $e_{1}$ and $f_{1}$ are diverse in $G$, it follows that
$e_{1}=a_{1}b_{1}$ say where
$a_{1},b_{1} \in \{x_{3},x_{4},x_{5},x_{7}\}$, that is,
$e_{1}$ is one of $x_{3}x_{4},x_{3}x_{7},x_{5}x_{7}$.
If $f_{1}$ is 3-4 or 3-7, let $P,Q$ be disjoint
paths of $G_{1}$ from $x_{2}$ to $x_{4}$ and from $t_{1}$
to $x_{5}$, with no vertices or edges in $D_{1}$ except
their ends; and let $R$ be a path of $G \setminus V(D)$ between
$V(P)$ and $V(Q)$ with no internal vertex or edge in $P$ or $Q$.
Then $D_{1} \cup P \cup Q \cup R$ is homeomorphic to Petersen,
and so $G_{1}$ and hence $H$ contains Petersen, and (ii) is true,
a contradiction.
So $e_{1}=x_{5}x_{7}$.
Let $P,Q$ be disjoint paths of $G_{1}$ from $t_{1}$ to $x_{3}$ and from
$x_{2}$ to $x_{5}$, with no vertices or edges in $D_{1}$ except their ends,
and let $R$ be as before.
Then $D_{1} \cup P \cup Q \cup R$ again is homeomorphic to Petersen,
a contradiction.

Hence $f_{1} \neq x_{1}y_{1}$, and so by symmetry
$f_{1} \neq x_{2}y_{2}, x_{3}y_{3}, x_{4}y_{4}$;
and hence $f=x_{5}y_{5}$.
Hence $e_{1}$ is 1-2 or 3-4, and by symmetry
we may assume the first.
Also, $e_{2}=x_{5} t_{1}$, and there are
(up to the symmetry) three possibilities for $f_{2}$, namely
$f_{2}=x_{1}y_{1},f_{2}=x_{4}y_{4}$, and
$f_{4} \in E(G \setminus V(D))$.
In the third case the theorem is true, so we assume for a contradiction
that one of the first two cases hold.
Let $G_{2} = G_{1} +(e_{2},f_{2})$, with new vertices
$s_{2},t_{2}$, and let $D_{2}$ be the subgraph of $G_{2}$ induced on
$V(D) \cup \{s_{1},t_{1},s_{2},t_{2}\}$.

If $f_{2}=x_{1}y_{1}$, let $P,Q$ be disjoint paths  of $G_{2}$
from $t_{2}$ to $x_{3}$ and from $t_{1}$ to $x_{4}$
with no vertices or edges in $D_{2}$ except their ends; then
$D_{2} \cup P \cup Q$ is homeomorphic to Petersen, a contradiction.
But if $f_{2} =x_{4}y_{4}$, let $P,Q$ be disjoint paths of $G_{2}$
from $x_{2}$ to $t_{2}$ and $t_{1}$ to $x_{3}$, with no vertices
or edges in $D_{2}$ except their ends; then
$D_{2} \cup P \cup Q$ is homeomorphic to Petersen, a contradiction.
This proves \ref{domino}.~\bbox

\bigskip

\noindent{\bf Proof of \ref{thetamain}.}
``Only if'' is easy and we omit it.
For ``if'', let $H$ be theta-connected and not contain Petersen,
Jaws or Starfish.
\\
\\
(1) {\em $H$ does not contain Antilog.}
\\
\\
\Subproof
Let $G$ be Antilog, let $X= \{1,\ldots,7\}$, and let
$D= G|X$.
Then $D$ is a crossed domino of $G$.
But the following all contain Petersen:
\begin{itemize}
\item[(i)]
$G+(e,f)$ for all $e \in E(D)$ and $f \in E(G \setminus X)$
\item[(ii)]
$G+(1\d 6,5\d 10) +(5\d 22,xy)$ for all
$xy \in E(G \setminus X)$ with $x,y \neq 10,14,15$.
\end{itemize}
From \ref{domino}, this proves (1).

\bigskip

Let ${\cal L}$ = \{Petersen, Jaws\}.
\\
\\
(2) {\em $H$ does not contain Log.}
\\
\\
\Subproof
Let Log(1) be Log + $(1\d 2,8\d 13)$, let $C$ be its
quadrangle, and let ${\cal L}_{1} = {\cal L} \cup$ \{Antilog\}.
All $A$-extensions are killed by ${\cal L}_{1}$, and
\[
{\cal P} (C, {\cal L}_{1} ) =
\{ (21\d 22,2\d 9), (21\d 22,13\d 9) \}^{\ast},
\]
and it follows by verifying the hypotheses of \ref{extn} that
$H$ does not contain Log(1).

Let Log(2) be Log + $(1\d 2,9\d 13)$, let $C$ be
its quadrangle, and ${\cal L}_{2} = {\cal L}_{1} \cup$ \{Log(1)\}.
All $A$-extensions are killed by ${\cal L}_{2}$, and
${\cal P}(C,{\cal L}_{2}) = \emptyset$, and so by \ref{Aextn}
$H$ does not contain Log(2).

Now let $G$ = Log, $X = {1,\ldots,7}$, and
${\cal L}_{3} = {\cal L}_{2} \cup$ \{Log(2)\}.
For any edge $e$ of $G|X$ and edge $f$ of $G$ not in $G|X$
(we permit $f$ to have one end in $X$), if $e,f$ are diverse then
$G+(e,f)$ contains a member of ${\cal L}_{3}$; and so $H$ does
not contain Log, by (1) and \ref{augment}.
This proves (2).
\\
\\
(3) {\em H does not contain Dice(1).}
\\
\\
\Subproof
Let Dice(11) = Dice(1) + $(1\d 2,20\d 23)$,
let $C$ be its quadrangle, and
${\cal L}_{4}$ = \{Petersen, Jaws, Log, Antilog\}.
All $A$-extensions are killed by ${\cal L}_{4}$, and
${\cal P}(C, {\cal L}_{4}) = \emptyset$, so by \ref{Aextn}
$H$ does not contain Dice(11).

Now let ${\cal L}_{5} = {\cal L}_{4} \cup$
\{Dice(11)\}, let $G$ = Dice(1), $X$ = \{1,\ldots,7\} and
$D=G|X$; then $D$ is a crossed domino in $G$.
For all $e \in E(D)$ and $f \in E(G \setminus X)$,
$G+(e,f)$ contains a member of ${\cal L}_{4}$; and for all
$xy \in E(G \setminus X)$ with
$x,y \neq 9,10,11$,
$G+(1\d 2,5\d 10) + (5\d 28,xy)$
contains Petersen.
Hence the result follows from \ref{domino}.
This proves (3).
\\
\\
(4) {\em $H$ does not contain Dice(2).}
\\
\\
\Subproof
Let $G$ = Dice(2), $X$ = \{1,\ldots,7\} and ${\cal L}_{6}$
= \{Petersen, Antilog, Dice(1)\}.
For all $e \in E(G|X)$ and $f \in E(G) \setminus E(G|X)$, if
$e,f$ have no common end then $G+(e,f)$ contains a member of
${\cal L}_{6}$; so (4) follows from (1), (3) and \ref{augment}.
\\
\\
(5) {\em $H$ does not contain Dice(3).}
\\
\\
\Subproof
Let Dice(31) = Dice(3) + $(3\d 4,13\d 14)$,
let $C$ be its quadrangle, and ${\cal L}_{4}$ as before.
All $A$-extensions are killed by ${\cal L}_{4}$, and
${\cal P}(C, {\cal L}_{4}) = \emptyset$, so by \ref{Aextn}
$H$ does not contain Dice(31).

Let ${\cal L}_{7} = {\cal L}_{4} \cup$ \{Dice(31)\}.
Let $G$ = Dice(3), $X = \{1,\ldots,7\}$, and
$D = G|X$.
Then $D$ is a crossed domino in $G$.
For all $e \in E(D)$ and $f \in E(G \setminus X)$,
$G+(e,f)$ contains a member of ${\cal L}_{7}$.
Moreover, for all $xy \in E(G\setminus X)$ with
$x,y \neq 15,16,18$,
\[
G+(1\d 2,5\d 15) + (5\d 28,xy) 
\]
\[
G+(3\d 4,5\d 15) + (5\d 28,xy)
\]
both contain Petersen or Log.
From (1)--(3) and \ref{domino}, this proves (5).
\\
\\
(6) {\em $H$ does not contain Dice(4).}
\\
\\
\Subproof
Let $G$ = Dice(4), $X$ = \{1,\ldots,7\} and $D=G|X$.
Then $D$ is a crossed domino in $G$.
But for all $e \in E(D)$ and $f \in E(G\setminus X)$,
$G+(e,f)$ contains Petersen or Log; and for all
$xy \in E(G \setminus X)$ with $x,y \neq 16,21,23$,
\[
G+(1\d 2,5\d 21) + (5\d 28,xy) 
\]
\[
G+(3\d 4,5\d 21) + (5\d 28,xy)
\]
both contain Petersen or Log.
The result follows from (2) and \ref{domino}. This proves (6).

\bigskip

From (1)--(6) and \ref{diamond}, this proves \ref{thetamain}.~\bbox

The reader may have noticed that Starfish hardly ever is needed
for anything.
There is an explanation, the following (previously stated as \ref{starfish0}).

\begin{thm}\label{starfish3}
Every dodecahedrally-connected graph $H$ containing Starfish
either is isomorphic to Starfish or contains Petersen.
\end{thm}
\Proof
If $H$ ``properly'' contains $G$ = Starfish, then by \ref{splitter} $H$
contains a graph $G' = G+(e,f)$ for some choice of
diverse edges $e,f$ of $G$.
But every such graph $G'$ contains Petersen.
This proves \ref{starfish3}.~\bbox

From \ref{starfish3} we obtain a slightly stronger reformulation of \ref{thetamain},
previously stated as \ref{main0}.

\begin{thm}\label{thetastar}
Let $H$ be theta-connected, and not isomorphic to Starfish.
Then $H$ is apex if and only if it contains neither of Petersen, Jaws.
\end{thm}

The proof is clear.

\section{Excluding Petersen}

In this section we prove \ref{main0}, thereby completing the proof of \ref{main}.
We restate it:

\begin{thm}\label{doublecross}
Let $H$ be theta-connected, and contain Jaws but not Petersen.
Then $H$ is doublecross.
\end{thm}
\Proof
Let Jaws(1) be Jaws $+ (1\d 2,3\d 4)$,
let Jaws(11) be Jaws(1) $+ (3\d 22,1\d 6)$, and let
Jaws(12) be Jaws(1) $+ (21\d 22,1\d 6)$.
\\
\\
(1) {\em H does not contain Jaws(11) or Jaws(12).}
\\
\\
\Subproof
Let $G$ be Jaws(11), and let
$X=V(G)  \setminus $ \{1, 2, 3, 21, 22, 23, 24\}.
If $ab \in E(G|X)$ and $cd \in E(G)  \setminus  E(G|X)$, with
$a,b \neq c,d$ and with $a,b$ non-adjacent to any of $c,d$
that are in $X$, then
$G+(ab,cd)$ contains Petersen.
Hence the result follows from \ref{augment} when $G$ is Jaws(11). 

When $G$ is Jaws(12), the argument is not so simple. Again we apply \ref{augment} to the same set
$X$. Let $(e_1,f_1)\l (e_k,f_k)$ be an augmenting sequence. By checking cases, we find that $f_1$ is not an
edge of $G\setminus X$ (because every choice of $e_1\in E(G|X)$ and $f_1\in E(G\setminus X)$ gives a Petersen), and
so $k\ge 2$; and having fixed $(e_1,f_1)$, we try all the possibilities for $(e_2,f_2)$. Again, there is no case with
$f_2\in E(G\setminus X)$, and so $k\ge 3$, and for each surviving choice of $(e_2,f_2)$ we try the possibilities for $(e_3,f_3)$. We
find in every case that there is no choice of $(e_3,f_3)$. (See the Appendix~\cite{RSTappendix}.)
This proves (1).
\\
\\
(2) {\em $H$ does not contain Jaws(1).}
\\
\\
\Subproof
Let $C$ be the quadrangle of $G$ = Jaws(1), and let
${\cal L}$ = \{Petersen, Jaws(11), Jaws(12)\}.
Then all $A$-extensions are killed by ${\cal L}$, and
${\cal P}(C, {\cal L}) = \emptyset$, so (2) follows from \ref{Aextn}.

\bigskip

Let Jaws(2) be Jaws $+(8,3,5,6)+(21,3,22,6)$, let Jaws(21) be Jaws(2) $+(6,7,11,12)$, and let Jaws(22) be Jaws(2) $+(7,8,19,10)$.
\\
\\
(3) {\em $H$ does not contain Jaws(21).}
\\
\\
We apply \ref{extn} to the quadrangle $\{25, 26, 12, 7 \}$, taking ${\cal L}$ to be \{Petersen, Jaws1\}. 
Again, see the Appendix for details. (Note that Jaws(21) has two circuits of length four, 
but it is quad-connected; this was the reason we extended
\ref{extn} to quad-connected graphs instead of graphs $G$ that were 
cyclically five-connected except for one circuit of length four.)
\\
\\
(4) {\em $H$ does not contain Jaws(22).}
\\
\\
This is easier; we apply \ref{Aextn} to the quadrangle $\{8, 20, 26, 25\}$, taking ${\cal L}$ to be \{Petersen, Jaws1, Jaws(21)\}.
\\
\\
(5) {\em $H$ does not contain Jaws(2).}
\\
\\
Let $X = \{ 6, 7, 8, 21, 22, 23, 24\}$. We apply \ref{augment} to $X$, and try all possibilities for the first three terms
of the augmenting sequence; and find in each case contains one of Petersen, Jaws(1), Jaws(21), Jaws(22). (See the Appendix.)

\bigskip

Now let ${\cal C}_{1}$ be the set of the seven circuits of Jaws
that bound regions in the drawing in Figure 2, not containing
1-6, 3-8, 13-18 or 15-20.
Let ${\cal C}_{2}$ be the set of paths of Jaws induced on
the following sets:
\begin{eqnarray*}
&6, 1, 2, 3, 8;&\\&
 8, 3, 4, 5, 6, 1;&\\& 
1, 6, 7, 8, 3; &\\&
3, 8, 20, 15;&\\&
15, 20, 19, 18, 13;&\\&
13, 18, 17, 16, 15, 20;&\\& 
20, 15, 14, 13, 18;&\\& 
18, 13, 1, 6.&
\end{eqnarray*}
Let $G$ = Jaws, let $F$ and $\eta_F$ be null, and let
${\cal C} = {\cal C}_{1} \cup {\cal C}_{2}$; then
$(G,F, {\cal C})$ is a framework.
By hypotheses, there is a \he $\eta$ of $G$ in $H$.
We claim that (E1)--(E7) hold.

Since $F$ is null, 
(E4), (E5) are vacuously true, and
(E1), (E3) are obvious.
It remains to check (E2), (E6) and (E7).
For (E2) we check that if $e,f \in E(G)$, not both in some
member of ${\cal C}$, then $G+(e,f)$ contains either Petersen
or Jaws(1); so (E2) follows from (2).
For (E6) it is only necessary to check cross extensions on the
circuit with vertex set \{4, 5, 11, 17, 16, 10\} and the path
with vertex set \{1, 6, 5, 4, 3, 8\}, since all the other
circuits and paths are too short or are equivalent by symmetry.
Hence we must check
\begin{eqnarray*}
&G+(4\d 5,16\d 17) + (4\d 21,17\d 22) &\\&
G+(4\d 5,16\d 17) + (4\d 10,11\d 17) &\\&
G+(4\d 10,11\d 17) + (4\d 21,17\d 22) &\\&
G+(4\d 10,11\d 17) + (10\d 16,5\d 11 ) &\\&
G+(3\d 8,5\d 6) + (3\d 4,1\d 6)&\\&
G+(3\d 8,5\d 6) + (3\d 21,6\d 22 ); &
\end{eqnarray*}
but they all contain Petersen, except the last which contains Jaws(2).
Hence (E6) holds.

For (E7) we must check
\[
G+(3\d 8,5\d 6 ) + (3\d 21,1\d 6) +
(8\d 21,1\d 24);
\]
but this contains Petersen.
Hence (E7) holds.
From \ref{mainthm}, this proves \ref{doublecross}.~\bbox

\end{document}